\documentclass[a4paper,10pt]{article}

\usepackage{setspace}
\usepackage[bottom]{footmisc}
\usepackage{color}
\usepackage{subfig}
\usepackage{hyperref}
\usepackage{multirow}
\usepackage{cite}
\usepackage{makecell}
\usepackage{amsmath} % easier math formulae, align, subequations \ldots
\usepackage{mathtools}
\usepackage{cleveref}
\usepackage{graphicx}
\usepackage{epsfig}
\usepackage{epstopdf}
\usepackage{float}
\usepackage{titling}
\usepackage[margin=1.0in]{geometry}
\usepackage{etoolbox}
 
\title{A review of troubled cell indicators for discontinuous Galerkin method on structured and unstructured meshes}
\makeatletter
\renewcommand\@date{{%
  \vspace{-\baselineskip}%
  \large\centering
  \begin{tabular}{@{}c@{}}
    S R Siva Prasad Kochi\textsuperscript{1} \\
    \normalsize siva.ksr@gmail.com
  \end{tabular}%
  \quad and\quad
  \begin{tabular}{@{}c@{}}
    M Ramakrishna\textsuperscript{2} \\
    \normalsize krishna@ae.iitm.ac.in
  \end{tabular}

  \bigskip

  \textsuperscript{1}Dept. of Aerospace Engg., IIT Madras.\par
  \textsuperscript{2}Professor, Dept. of Aerospace Engg., IIT Madras.

  \bigskip

  \today
}}
\makeatother
\begin{document}

\maketitle

\begin{abstract}
 In this paper, eight different troubled cell indicators (shock detectors) are reviewed for the solution of nonlinear hyperbolic conservation laws using discontinuous Galerkin (DG) method and a WENO limiter on both structured and unstructured meshes. Extensive simulations using one-dimensional and two-dimensional problems (2D Riemann problem and the double Mach reflection) for various orders on the hyperbolic system of Euler equations are used to compare these troubled cell indicators. They are evaluated based on the percentage of cells flagged as troubled cells for various orders and various grid sizes. CPU time taken to test a single cell for discontinuity is also compared. For one-dimensional problems, the performance of Fu and Shu indicator \cite{fs1} and the modified KXRCF indicator \cite{lpr} is better than other indicators. For two-dimensional problems, the performance of the artificial neural network (ANN) indicator of Ray and Hesthaven \cite{rh1} is quite good and the Fu and Shu and the modified KXRCF indicators are also good. These three indicators are suitable candidates for applications of DGM using WENO limiters though it should be noted that the ANN indicator is quite expensive and requires a lot of training.
 
 {{\bf Keywords:} discontinuous Galerkin method, troubled cell indicators, shocks, contact discontinuities}
\end{abstract}

\section{Introduction}
\label{intro}
\noindent In this paper, we review different troubled cell indicators using discontinuous Galerkin method (DGM). In 2005, Qiu and Shu \cite{qs1} evaluated various troubled cell indicators based on the percentage of cells flagged as troubled cells on structured meshes using various one-dimensional test cases and the double Mach reflection problem. We review troubled cell indicators developed after their review which are applicable for both structured and unstructured grids. We also use several one-dimensional test cases. In 2D, we use the 2D Riemann problem and the double Mach reflection problem on both structured and unstructured meshes. We also show the time history of the detected troubled cells for one-dimensional problems. CPU time taken to test a single cell for discontinuity is also compared for the various troubled cell indicators. We have identified the Fu and Shu indicator \cite{fs1}, the modified KXRCF indicator \cite{lpr} and the artificial neural network (ANN) indicator of Ray and Hesthaven \cite{rh1} as suitable candidates for applications of DGM using WENO limiters.
\\
\\
For solutions of nonlinear conservation laws using the discontinuous Galerkin method \cite{cs6}, Qiu and Shu \cite{qs1} adopted a new methodology for limiting. They first identify \textit{troubled cells}, namely, those cells where limiting is needed and use a limiter only in those troubled cells. This method is extremely popular and is usually the default process for limiting using the discontinuous Galerkin method. This makes the troubled cell indicator (or shock detector) which identifies the troubled cells very important. Qiu and Shu \cite{qs2} reviewed many troubled cell indicators available at the time and concluded that the minmod-based TVB indicator (with a constant $M$ chosen carefully) \cite{cs2}, the KXRCF indicator by Krivodonova et al \cite{kxrcf}, and an indicator based on Harten's subcell resolution idea \cite{har2} are better than others. Later, Zhu and Qiu \cite{zq2}, \cite{zq3} investigated these troubled cell indicators for adaptive discontinuous Galerkin methods.
\\
\\
Here, we consider only those troubled cell indicators which are applicable for both structured and unstructured grids. For this reason, we didn't consider the troubled cell indicators developed by Balsara et al \cite{bamd}, Don et al \cite{dglw}, Fu \cite{fu1}, and Huang and Fu \cite{hf1}. We also don't consider the various \textit{a posteriori} troubled cell indicators \cite{jyy}, \cite{bat} some of them based on the Multi-dimensional Optimal Order Detection (MOOD) approach \cite{dzls}, as they require special study.
\\
\\
\noindent Discontinuous Galerkin method was first introduced by Reed and Hill \cite{rh} and later was developed for nonlinear conservation laws as the Runge-Kutta discontinuous Galerkin (RKDG) method by Cockburn et al \cite{cs6}. For solutions containing strong shocks, RKDG method uses a nonlinear limiter which is applied to detect discontinuities and control spurious oscillations near such discontinuities. Many such limiting strategies were developed in literature \cite{zcq} and WENO (weighted essentially non oscillatory) limiters are preferred as they maintain the order of the scheme. Some of the earlier limiting strategies have their own troubled cell indicators and most of them were reviewed by Qiu and Shu \cite{qs2}. We use the compact subcell WENO limiting strategy \cite{srspkmr1}, \cite{srspkmr2} developed by us for our review of troubled cell indicators.
\\
\\
\noindent In recent years, troubled cell indicators are attracting a lot of attention. There are troubled cell indicators coupled with a artificial viscosity addition \cite{pp},\cite{lodato1} and those developed based on edge detection techniques for spectral data \cite{sj}. We also have shock detectors based on various derivatives for higher order methods \cite{olll}. Hierro et al \cite{hbk} developed a novel troubled-cell indicator and hp-adaptivity scheme that exploits the adaptivity history. Vuik and Ryan \cite{vr2} addressed the issue of problem-dependent parameters that exist in troubled cell indicators and proposed a novel automatic parameter selection strategy. Later, Fu and Shu \cite{fs1} proposed a new troubled cell indicator that contains a parameter depending only on the DG polynomial degree. Zhu et al \cite{zhw} generalized their troubled cell indicator to h-adaptive meshes. Maltsev et al \cite{myjst} used this as a shock detector for a hybrid discontinuous Galerkin - finite volume method to great effect. The much praised KXRCF troubled cell indicator was modified for DG spectral element methods by Li et al \cite{lpr}. We have also modified the troubled cell indicator of Fu and Shu \cite{fs1} and used it to obtain results similar to that of the parent indicator. We have included these two indicators for review in this paper.
\\
\\
\noindent Recently, neural networks in deep learning were introduced to design new troubled cell indicators. The works of Ray and Hesthaven \cite{rh1}, \cite{rh2} and Feng et al \cite{flw} use feed-forward neural networks and the works of Sun et al \cite{swcxx} and Wang et al \cite{wzcx} are based on convolutional neural networks. Beck et al \cite{bzsf} used a neural network based on the so called holistically nested edge detection (HED) algorithm for shock detection. Xiao et al \cite{xgs} developed an edge detection algorithm based on Bayesian learning. Recently, Zhu et al \cite{zwg} proposed a new framework using K-means clustering and improved it further in \cite{zwwzg} and \cite{wgwzz}. It is notable that some of these frameworks were designed for only structured grids and not yet been extended to unstructured meshes. Design of a robust and accurate troubled cell indicator is a very important problem and we review some of the recent works on this topic which are simple and are applicable on both structured and unstructured grids.
\\
\\
\noindent We evaluate the performance of the selected troubled cell indicators using a process similar to that followed in Qiu and Shu \cite{qs2}. We look at the percentage of cells flagged as troubled cells for various orders and various grid sizes and also show the time history of troubled cells for some specific test cases. We have observed that for each of the selected troubled cell indicators, the CPU time taken to test a single cell for discontinuity is approximately the same (excluding the time required to train a neural network for the artificial neural network indicator of Ray and Hesthaven \cite{rh1}). We have also observed that as long as we use an indicator which identifies all the troubled cells, if we use a good limiter, the accuracy of the solution remains the same. But the additional troubled cells add to the computational time and reduce the efficiency of the scheme. We found the solution to be non-oscillatory using all the troubled cell indicators used.
\\
\\
\noindent The paper is organized as follows. We describe the formulation of the discontinuous Galerkin method used for all our results in section \ref{sec:1}, the discontinuities detection methods reviewed are described in section \ref{sec:2} and the results are described in section \ref{sec:3} and we conclude the paper in section \ref{sec:4}.

\section{Formulation of discontinuous Galerkin method}\label{sec:1}

\noindent Consider a two-dimensional conservation law of the form

\begin{equation}\label{2dConsLaw}
\frac{\partial Q}{\partial t} + \frac{\partial F(Q)}{\partial x} + \frac{\partial G(Q)}{\partial y} = 0 \quad \text{in the domain} \quad \Omega
\end{equation}

\noindent with the initial condition

\begin{displaymath}
    Q(x,y,0) = Q_{0}(x,y),
\end{displaymath}

\noindent We approximate the domain $\Omega$ by $K$ non-overlapping elements given by $\Omega_{k}$, $k=1,2,\ldots K$. 
\\
\\
\noindent We look at solving \eqref{2dConsLaw} using the discontinuous Galerkin method. We approximate the local solution as:

\begin{equation}\label{modalApprox}
 Q_{h}(r,s) = \sum_{i=0}^{N_{p}-1} \hat{Q}^{i} \psi_{i}(r,s)
\end{equation}

\noindent where $r$ and $s$ are the local coordinates, $h$ the grid size, and $\hat{Q}^{i}$ are the degrees of freedom representing the approximate solution $Q_{h}$ and $N_{p}=(N+1)(N+1)$ are the number of degrees of freedom ($N$ is the degree of the one-dimensional representation polynomial) for quadrilaterals and $N_{p}=(N+1)(N+2)/2$ for triangles. The basis chosen is the tensor product orthonormalized Legendre polynomials for quadrilaterals and polynomial basis given in \cite{hestha1} for triangles.  Now, using $\psi_{j}(r,s)$ as the test function, the weak form of the equation \eqref{2dConsLaw} is obtained as

\begin{eqnarray}\label{weakFormScheme}
 \sum_{i=0}^{N_{p}-1} \frac{\partial \hat{Q}^{i}}{\partial t} \int_{\Omega_{k}} \psi_{i} \psi_{j} d\Omega + \int_{\partial \Omega_{k}} \hat{F} \phi_{j} ds - \int_{\Omega_{k}} \vec{F} \cdot \nabla \phi_{j} d\Omega = 0
\end{eqnarray}

\noindent where $\partial \Omega_{k}$ is the boundary of $\Omega_{k}$, $\vec{F} = (F(Q),G(Q))$ and $\hat{F} = \bar{F^{*}}\cdot \hat{n}$ where $\bar{F^{*}}$ is the monotone numerical flux at the interface which is calculated using an exact or approximate Riemann solver and $\hat{n}$ is the unit outward normal. We use the Lax-Friedrichs numerical flux unless otherwise specified. This formulation is termed to be $P^{N}$ based discontinuous Galerkin method. Equation \eqref{weakFormScheme} is integrated using an appropriate Gauss Legendre quadrature and is discretized in time by using an appropriate Runge-Kutta time discretization given in \cite{shu} unless otherwise specified.
\\
\\
\noindent To control spurious oscillations that occur near a discontinuity, either a nonlinear limiter is applied or artificial viscosity is added near such discontinuities. This is known as shock capturing in literature. The general strategy for shock capturing is
\\
\\
\noindent \textbf{1)} Identify the cells which need to be limited, known as troubled cells. \\
\noindent \textbf{2)} Limit the solution polynomial or add artificial viscosity in troubled cells.
\\
\\
\noindent The first step requires a shock detector or a troubled cell indicator (as it is called in DGM terminology). In this paper, we investigate various troubled cell indicators and their properties. For step 2, we use the compact subcell weighted essentially non oscillatory limiting (CSWENO) method proposed in \cite{srspkmr1} and \cite{srspkmr2}.

\section{Discontinuity detection methods}\label{sec:2}

\noindent The discontinuity detection methods we intend to review are listed out below:
\\
\\
\noindent \textbf{1) Subcell indicator of Persson and Peraire \cite{pp} - abbreviated as PP indicator}: This indicator requires the solution to be expressed in terms of an orthogonal basis in each element:

\begin{equation}\label{modalApprox2}
 Q = \sum_{i=1}^{N_{p}(N)} \hat{Q}^{i} \psi_{i}
\end{equation}

\noindent where $N_{p}(N)=(N+1)(N+1)$ for quadrilaterals and $N_{p}(N)=(N+1)(N+2)/2$ for triangles. We now consider a truncated expansion of the same solution with terms only up to order $N-1$. It is written as:

\begin{equation}\label{modalApprox3}
 \overline{Q} = \sum_{i=1}^{N_{p}(N-1)} \hat{Q}^{i} \psi_{i}
\end{equation}

\noindent Within each element $\Omega_{k}$, the troubled cell indicator is defined as

\begin{equation}\label{ppDetector}
 T_{k} = \frac{(Q-\overline{Q})_{L^{2}}^{2}}{(Q)_{L^{2}}^{2}}
\end{equation}

\noindent where $(.)_{L^{2}}$ represents the $L^{2}$-norm of the quantity in brackets. Assuming that the polynomial expansion has a similar behavior to the Fourier expansion, Persson and Peraire concluded that the value of $T_{k}$ will scale like $1/N^{4}$ for continuous functions. Based on this, the element $\Omega_{k}$ is a troubled cell if $T_{k}>1/N^{4}$. This works quite well as a troubled cell indicator and Persson and Peraire \cite{pp} used the quantity $T_{k}$ to even determine the artificial viscosity to be added. This was improved further to give quite good results by Klockner et al \cite{kwh}.
\\
\\
\noindent \textbf{2) Concentration method of Sheshadri and Jameson \cite{sj} - abbreviated as SJ indicator}: The concentration method, proposed by Gelb, Cates and Tadmor \cite{gt1},\cite{gt2},\cite{gc} is a general framework for detecting jump discontinuities in piecewise smooth functions using their spectral data. Their approach is based on localization using appropriate concentration kernels and separation of scales by nonlinear enhancement. Sheshadri and Jameson used the extended framework for polynomial modes (instead of Fourier modes) \cite{gt2} as a discontinuity detector. This only works for polynomial basis derived from Chebyshev and Jacobi polynomials and hence will work for a Legendre polynomial basis. The step by step process for discontinuity detection in 1D is given below:
\\
\\
\noindent Step 1: At each solution point $x$ (typically a Gauss-Legendre quadrature point), compute the quantity

\begin{displaymath}
 \frac{\pi}{N_{p}}\sqrt{1-x^{2}}\sigma\left(\frac{|k|}{N_{p}}\right)\psi'_{k}(x) \quad \text{for} \quad k = 1,2,\ldots N_{p}
\end{displaymath}

\noindent where $N_{p}=N+1$ ($N$ is the degree of the polynomial basis) and $\sigma(.)$ are the concentration factors. A variety of concentration factors are available for this purpose and we have used the polynomial concentration factors given by Gelb and Tadmor in \cite{gt2} for all our calculations. Since the solution points are fixed, this gives us a $N_{p}\times N_{p}$ matrix for each element. This is known as the concentration matrix.
\\
\\
\noindent Step 2: Find the modal matrix of each element (again a $N_{p}\times N_{p}$ matrix) or convert the local nodal solution to modal form using the Vandermonde matrix.
\\
\\
\noindent Step 3: Find the concentration kernel $K$ by a matrix multiplication of the concentration matrix and the modal matrix.
\\
\\
\noindent Step 4: Find the enhanced kernel $K^{p}_{N_{p},J}$ by using the following function for each element of $K$:
\begin{equation}
[K^{p}_{N_{p},J}]_{ij} = 
\begin{cases}
    K_{ij},& \text{if } N_{p}^{p/2}|K_{ij}|^{p}>J\\
    0,              & \text{otherwise}
\end{cases}
\end{equation}

\noindent where $p$ ($>1$) is the enhancement exponent typically chosen to be 2 or 3 and $J$ is an appropriately chosen threshold.
\\
\\
\noindent Step 5: Identify the points for which the enhanced kernel is zero as a point of discontinuity. If an element contains a point of discontinuity, it is identified as a troubled cell. This completes the procedure for shock detection.

\begin{figure}[htbp]
\begin{center}
\includegraphics[scale=1.0]{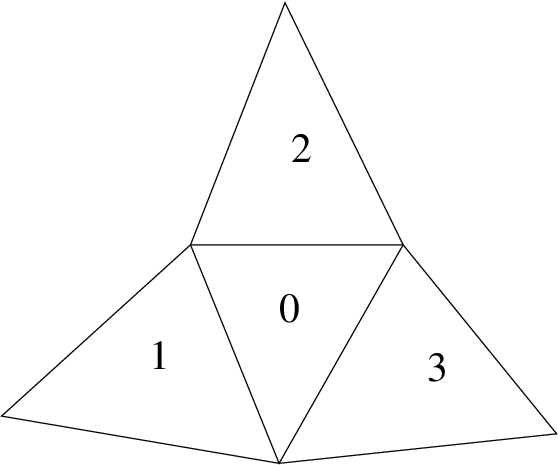}
\caption{The stencil S = $\{ \Omega_{0}, \Omega_{1}, \Omega_{2}, \Omega_{3}\}$ for Fu and Shu troubled cell indicator \cite{fs1}}
\label{fig1:stencilFS}
\end{center}
\end{figure}

\noindent This procedure can be extended to 2D problems in a dimension by dimension fashion as given in \cite{sj}. This shock detector is a little expensive computationally but works really well for some class of problems as shown by \cite{sj}.
\\
\\
\noindent \textbf{3) Troubled cell indicator of Fu and Shu \cite{fs1} - abbreviated as FS1 indicator}: To check whether a cell is troubled, Fu and Shu considered a stencil which contains the target cell and all its immediate neighbors. For example, consider a target cell $\Omega_{0}$ as shown in Figure \ref{fig1:stencilFS}. The indicator stencil is S = $\{ \Omega_{0}, \Omega_{1}, \Omega_{2}, \Omega_{3}\}$. Now consider the DG polynomial in each of these four cells as $Q_{l}(x,y)$, $l=0,1,2,3$.
\\

\noindent We define 

\begin{equation}\label{cellAvgDef}
 \overline{\overline{Q}}_{j} = \frac{1}{\Omega_{0}} \int_{\Omega_{0}} Q_{j}(x,y) dxdy, \quad \overline{Q} = \frac{1}{\Omega_{j}} \int_{\Omega_{j}} Q_{j}(x,y) dxdy
\end{equation}

\noindent where $\overline{\overline{.}}$ denotes the cell average of the function $.$ on the target cell while $\overline{.}$ denotes the cell average of the function $.$ on its own associated cell. In order to evaluate $\overline{\overline{Q}}_{l}$ for $l=1,2,3$, the polynomial data are naturally extrapolated to the target cell $\Omega_{0}$. We now define

\begin{equation}\label{FSTC1}
 T_{\Omega_{0}} = \frac{\sum_{j=1}^{3} |\overline{\overline{Q}}_{0}-\overline{\overline{Q}}_{j}|}{\text{max}_{j\in \{0,1,2,3\}} \{|\overline{Q}_{j}|\}}
\end{equation}

\noindent The target cell $\Omega_{0}$ is treated as a troubled cell if

\begin{equation}\label{FSTC2}
 T_{\Omega_{0}} > C_{N}
\end{equation}

\noindent for a constant $C_{N}$ that depends only on the degree of the polynomial basis $N$. For 1D problems, the values of $C_{N}$ are $0.05, 0.1, 0.25, 0.5$ for $N=1,2,3,4$ respectively \cite{fs1}. Fu and Shu compared this troubled cell indicator with one of the best shock detectors - the so called KXRCF troubled cell indicator \cite{kxrcf} and obtained very favourable results for many test problems. This indicator has also been used as a shock detector for a hybrid discontinuous Galerkin - finite volume method by Maltsev et al \cite{myjst} to great effect.
\\
\\
\noindent \textbf{4) Modified Fu and Shu troubled cell indicator - abbreviated as FS2 indicator}: We have modified the FS1 troubled cell indicator described above and used it to obtain similar results. We use the same stencil S = $\{ \Omega_{0}, \Omega_{1}, \Omega_{2}, \Omega_{3}\}$ and define:

\begin{equation}\label{cellAvgDef2}
 \overline{\overline{Q}}_{j} = \frac{1}{\Omega_{j}} \int_{\Omega_{j}} Q_{0}(x,y) dxdy, \quad \overline{Q} = \frac{1}{\Omega_{j}} \int_{\Omega_{j}} Q_{j}(x,y) dxdy
\end{equation}

\noindent Here, in contrast to FS1, while evaluating $\overline{\overline{Q}}_{l}$ for $l=1,2,3$, instead of extrapolating the polynomial data to the target cell $\Omega_{0}$, we extrapolate the polynomial data of the target cell $Q_{0}(x,y)$ to the neighboring cells. Now we use the same equations \eqref{FSTC1}, \eqref{FSTC2} to determine whether $\Omega_{0}$ is a troubled cell or not. This indicator works quite well as the tests below show.
\\
\\
\noindent \textbf{5) Modified KXRCF troubled cell indicator of Li et al \cite{lpr} - abbreviated as LPR indicator}: Li et al modified the well regarded KXRCF troubled cell indicator \cite{kxrcf} and constructed their own limiter based on their modification. They first used the method employed in KXRCF troubled cell indicator and found:

\begin{equation}
 T_{1} = \frac{|\int_{\partial \Omega_{i}^{-}} (Q_{i}-Q_{j})|}{h |\Omega_{i}^{-}|\overline{Q}_{i}}
\end{equation}

\noindent where $h$ is the grid size, $\Omega_{i}^{-}$ represents the inflow surface of element $i$, and element $j$ is adjacent to element $i$ with $\Omega_{i}^{-}$ as the interface. $\overline{Q}_{i}$ is the cell average of element $i$. The cell is marked as a troubled cell if $T_{1}>C_{1}$, where $C_{1}$ is 1 for 1D problems and 2D triangles and 2 for quadrilaterals.
\\
\\
\noindent To improve this indicator, they constructed linear candidate polynomials in physical space in the following fashion. $Q(r,s)$ is the higher order polynomial representation of $Q$ in a given element with local coordinates $r$ and $s$. These polynomials in the candidate troubled cell (identified using $T_{1}$) and its immediate neighbors as shown in Figure \ref{fig1:stencilFS} (the stencil S = $\{ \Omega_{0}, \Omega_{1}, \Omega_{2}, \Omega_{3}\}$) are projected as linear polynomials in physical space. These polynomials are written as

\begin{equation}
 Q1_{k}(x,y) = \overline{Q}_{i} + Q1_{x}(x-x_{i}) + Q1_{y}(y-y_{i}), \quad k = 0,1,2,3
\end{equation}

\noindent where $k=0$ is for the candidate troubled cell and $k=1,2,3$ for its immediate neighbors. $(x_{i},y_{i})$ is the center of gravity of the element $i$. The polynomials $Q1_{k}(x,y)$ are obtained by minimizing the function

\begin{equation}
 L = \int_{\Omega_{j}} (Q1_{k}(x,y)-Q_{j}(r,s))^{2} d\Omega, \quad j\in \text{S}
\end{equation}

\noindent The coefficients $Q1_{x}$ and $Q1_{y}$ are obtained by using

\begin{equation}
 \frac{\partial L}{\partial Q1_{x}} = 0, \quad \frac{\partial L}{\partial Q1_{y}} = 0
\end{equation}

\noindent After obtaining $Q1_{x}$ and $Q1_{y}$, we calculate

\begin{equation}
 T_{2} = \frac{\sum_{k=1}^{3}|Q1_{k}(x_{k},y_{k})-\overline{Q}_{k}|}{h \times \text{max}_{l\in \{0,1,2,3\}} \{|\overline{Q}_{l}|\}}
\end{equation}

\noindent If $T_{2}>C_{2}$ with $C_{2}$ again being 1 for 1D problems and triangles and 2 for quadrilaterals, the element $i$ is marked as a troubled cell, otherwise it is considered a smooth cell. This indicator has been used by Li et al for various test problems to great effect \cite{lpr}.
\\
\\
\noindent \textbf{6) Troubled cell indicator using an artificial neural network \cite{rh1}, \cite{rh2} - abbreviated as RH indicator}: Ray and Hesthaven constructed a troubled cell indicator using an artificial neural network for one-dimensional conservation laws in \cite{rh1} and extended it to two-dimensional unstructured grids in \cite{rh2}. We briefly describe their procedure. Using an artificial neural network (ANN), our aim is to find an approximation to an unknown function $\mathbf{G}$ given a large data set $\mathbf{T}$. For our problem, this corresponds to the underlying true indicator function that correctly flags the genuine troubled-cells.
\\
\\
\noindent An artificial neural network is described mathematically by the triplet ($A,V,W$). $A$ is the set of all artificial neurons in the network, $V$ represents the set of all directed connections $(i,j)$, $i,j\in A$, where $i$ is the sending neuron and $j$ is the receiving neuron, and $W$ represents set of weights $w_{i,j}$ used for the neural connections $(i,j)$. Each artificial neuron receives a lot of signals whose weighted accumulation is stored in terms of a propagation function which is usually chosen to be linear. Once the accumulation crosses a certain threshold (called a bias), the neuron transmits a scalar-valued signal, which is modeled using an activation function. Several architectures which describe the arrangement and connectivity of neurons were proposed for the modeling of ANN's \cite{schmidhuber1}. We follow Ray and Hesthaven and use the architecture known as the multi-layer perceptron (MLP), where the neurons are arranged in multiple layers (an input layer, an output layer and several in-between layers called hidden layers).
\\
\\
\noindent The appropriate weights and biases are obtained by training the network. As with Ray and Hesthaven, we use the training paradigm known as supervised learning \cite{kriesel1}, where the training aims to minimize the error between prediction and truth. In this case, the true responses of the training set are known a priori. For a candidate troubled cell $i$, we use a input layer with $5$ neurons corresponding to the vector $(\overline{Q}_{i-1},\overline{Q}_{i},\overline{Q}_{i+1},\overline{Q}_{i-1/2}^{+},\overline{Q}_{i+1/2}^{-})$, where the $\overline{Q}_{i}$ represents the cell average of $Q$ in cell $i$. We use 5 hidden layers and an output layer of 2 neurons (which corresponds to a probability of $i$ being a troubled cell). We train the network using a stochastic optimization algorithm as in \cite{rh1}. The training is very costly, but it needs to be done just once, and can be used as a troubled cell indicator for all one-dimensional conservation laws.
\\
\\
\noindent For two-dimensional unstructured grids, we follow \cite{rh2} and again train an MLP network to perform as a troubled cell indicator. In this case, we use an input layer with 12 neurons with quantities selected from the cell and its immediate neighbors. For an element shown in Figure \ref{fig1:stencilFS}, we use the modal coefficients $(\hat{Q}^{0},\hat{Q}^{1},\hat{Q}^{N+1})$ for the element and its immediate neighbors, when a basis of order $N$ is used. This gives the 12 quantities required for the input layer. We again use 5 hidden layers and a stochastic optimization algorithm as in \cite{rh2} to train the network. Again, the training is done just once and is used for all two-dimensional test cases. Using this troubled-cell indicator, we obtained results quite similar to that of Ray and Hesthaven \cite{rh1},\cite{rh2}.
\\
\\
\noindent \textbf{7) Characteristic modal shock detector of Lodato \cite{lodato1} - abbreviated as PPL indicator}: This troubled cell detection strategy is quite simple and is based on the subcell detector of Persson and Peraire \cite{pp} (and hence the label PPL for the detector). For a given conservative variable vector ($Q$) in a given cell $i$, we compute its values in the characteristic space by premultiplying with the left eigenvectors of the conservative Jacobian. The characteristic direction is chosen to coincide with the physical direction corresponding to the cell average of the velocity vector as explained in \cite{lodato1}. We convert these nodal quantities in characteristic space to modal quantities by multiplying with the appropriate Vandermonde matrix. Now these modal quantities in characteristic space are used to find $T_{k}$ as in \eqref{ppDetector}. Now, this quantity is converted back to the physical space to be used as a troubled cell indicator as done by Persson and Peraire. Lodato \cite{lodato1} used this to calculate the artificial viscosity to be added quite effectively. We only use the shock detecting feature of this process for our test cases.
\\
\\
\noindent \textbf{8) Shock detector of Maeda and Heijima \cite{mh1} - abbreviated as MH indicator}: Maeda and Heijima used the pressure-change indicator and a compressible sonic-point condition to detect shocks accurately in a given element. We describe their procedure for two-dimensional problems. The pressure-change indicator \cite{klr}, \cite{zlcz} is defined as:

\begin{equation}\label{pCInd}
 f = \text{min}\left(\frac{p_{L}}{p_{R}},\frac{p_{R}}{p_{L}}\right)^{3}
\end{equation}

\noindent where $p_{L}$ and $p_{R}$ are static pressures in adjacent elements. When $f<f_{c}$ (where $f_{c}=0.5$ \cite{zlcz}), the cell is a candidate for the existence of a shock. Next, a local angle is calculated from the coordinates of two cell centers ($B$ and $C$), where one is the candidate cell and the other - one of its neighbors that is supposed to be a part of shock:

\begin{equation}
 \beta = \tan^{-1}\left(\frac{y_{C}-y_{B}}{x_{C}-x_{B}}\right)
\end{equation}

\noindent Now, the sonic point condition is considered between two cells striding across the line segment $BC$ with

\begin{equation}
 V_{n} = u\sin \beta - v\cos \beta 
\end{equation}

\noindent where $u$ and $v$ are the $x$ and $y$ components of the velocity. The sonic point condition now is 

\begin{equation}
 V_{nL}-c_{L}>0 \quad \text{and} \quad V_{nR}-c_{R}<0
\end{equation}

\noindent or

\begin{equation}
 V_{nL}+c_{L}>0 \quad \text{and} \quad V_{nR}+c_{R}<0
\end{equation}

\noindent where $c$ is the sonic velocity. The subscripts $L$ and $R$ represent adjacent cells being considered. If this condition is satisfied, we mark the cell as a troubled cell. This has to be repeated for all the neighbors of the candidate cell (even neighbors sharing nodes). To detect a contact discontinuity, we use density instead of pressure in \eqref{pCInd}. For one-dimensional problems, we simply use equation \eqref{pCInd} and consider an element as a troubled cell if $f<f_{c}$. Maeda and Heijima used this shock detector for a variety of two-dimensional problems with good results.

\section{Results}\label{sec:3}

\noindent In this section, we use different numerical test cases to compare the eight different troubled-cell indicators outlined in the previous section with the discontinuous Galerkin method while using a WENO reconstruction limiter \cite{srspkmr1}, \cite{srspkmr2}.

The CPU time taken to test a single cell for discontinuity is approximately the same for all the troubled cell indicators (excluding the training time for the RH indicator). The normalised CPU time for each troubled cell indicator, normalised with respect to the maximum value is given in Table \ref{table:0New}. From this table, we can clearly see that the CPU time taken to test a single cell for discontinuity is not an important factor for the chosen troubled cell indicators. Hence the number of troubled cells identified near a discontinuity is the most important factor affecting the solution scheme. We therefore look at the percentage of cells flagged as troubled cells for evaluating the troubled cell indicators for both structured and unstructured meshes. The time history of number of troubled cells flagged using each troubled-cell indicator is also shown.

\begin{figure}[htbp]
  \centering
  \subfloat[Solution of 1D Euler equations for the single moving contact discontinuity with $P^{1}$, $P^{2}$ and $P^{3}$ based DGM and SJ indicator]{\label{fig:CDSolution}\includegraphics[width=0.48\textwidth]{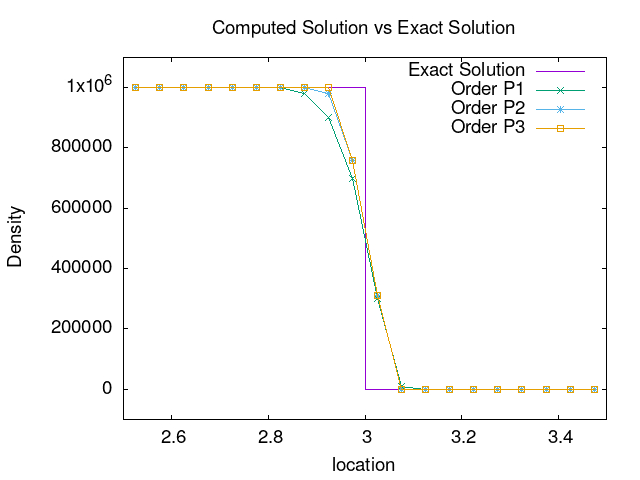}}
  \subfloat[Density error ($|\rho-\rho_{exact}|$) for $P^{1}$, $P^{2}$ and $P^{3}$ based DGM]{\label{fig:CDError}\includegraphics[width=0.48\textwidth]{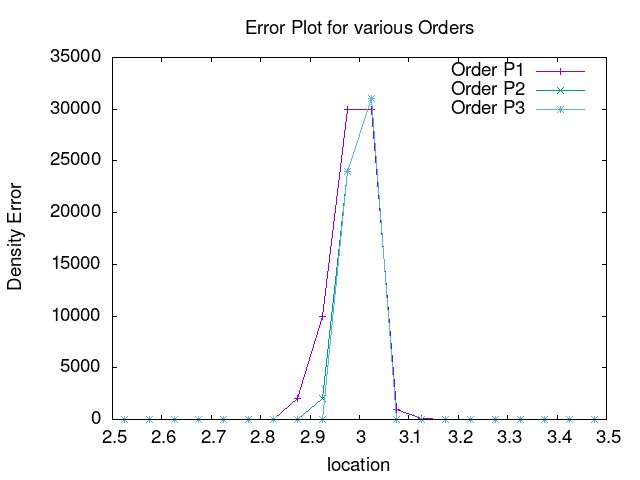}}\hfill
  \caption{Comparison of density solutions of single moving contact discontinuity at $t=3.0$ using 200 elements obtained with the $P^{1}$, $P^{2}$ and $P^{3}$ based DGM and the exact solution on [2.5,3.5]. Density error ($|\rho-\rho_{exact}|$) for $P^{1}$, $P^{2}$ and $P^{3}$ based DGM is also plotted}
  \label{fig:CDCSWENO}
\end{figure}

\begin{figure}[htbp]
  \centering
  \subfloat[PP Indicator]{\label{fig:SCDPP}\includegraphics[width=0.32\textwidth]{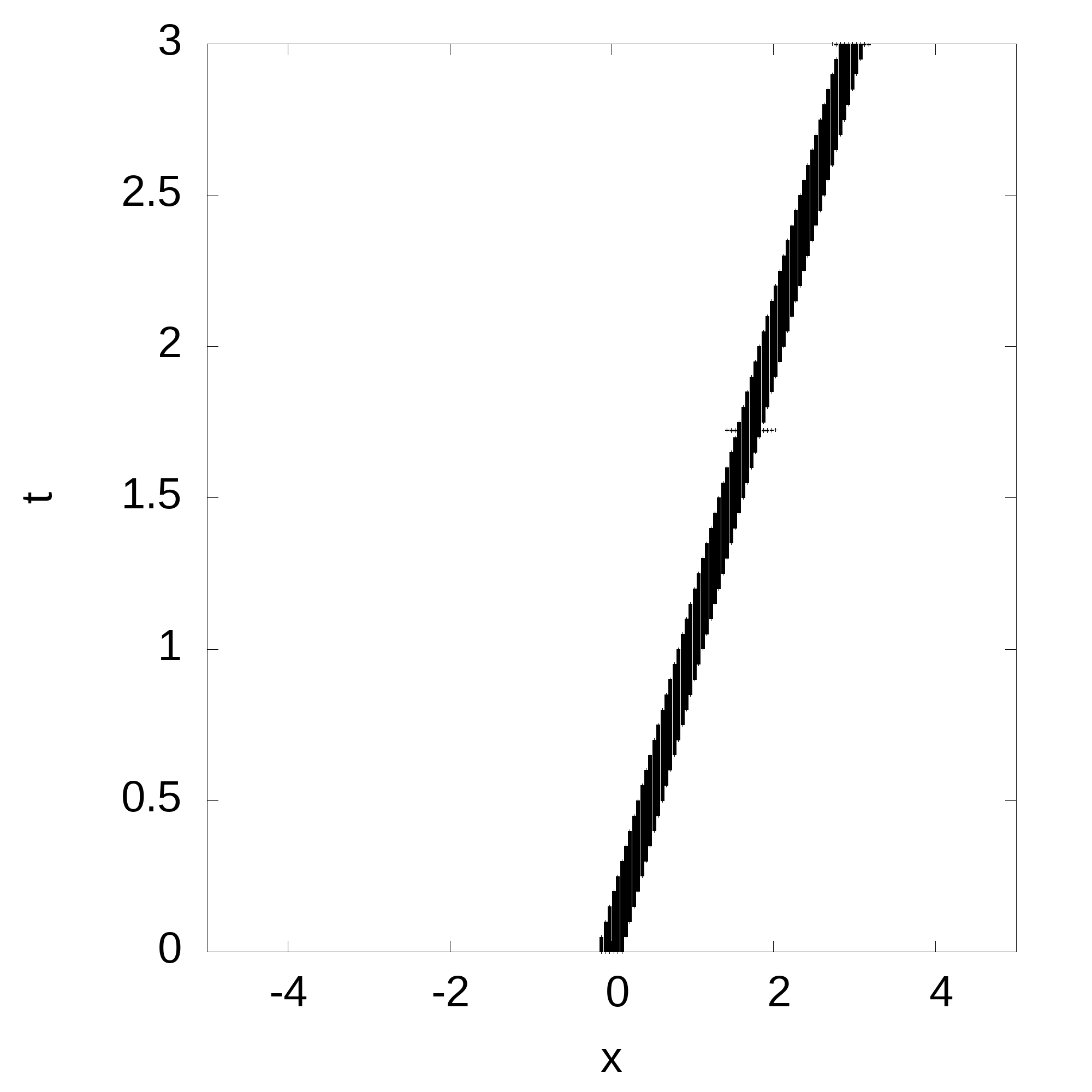}}
  \subfloat[SJ Indicator]{\label{fig:SCDSJ}\includegraphics[width=0.32\textwidth]{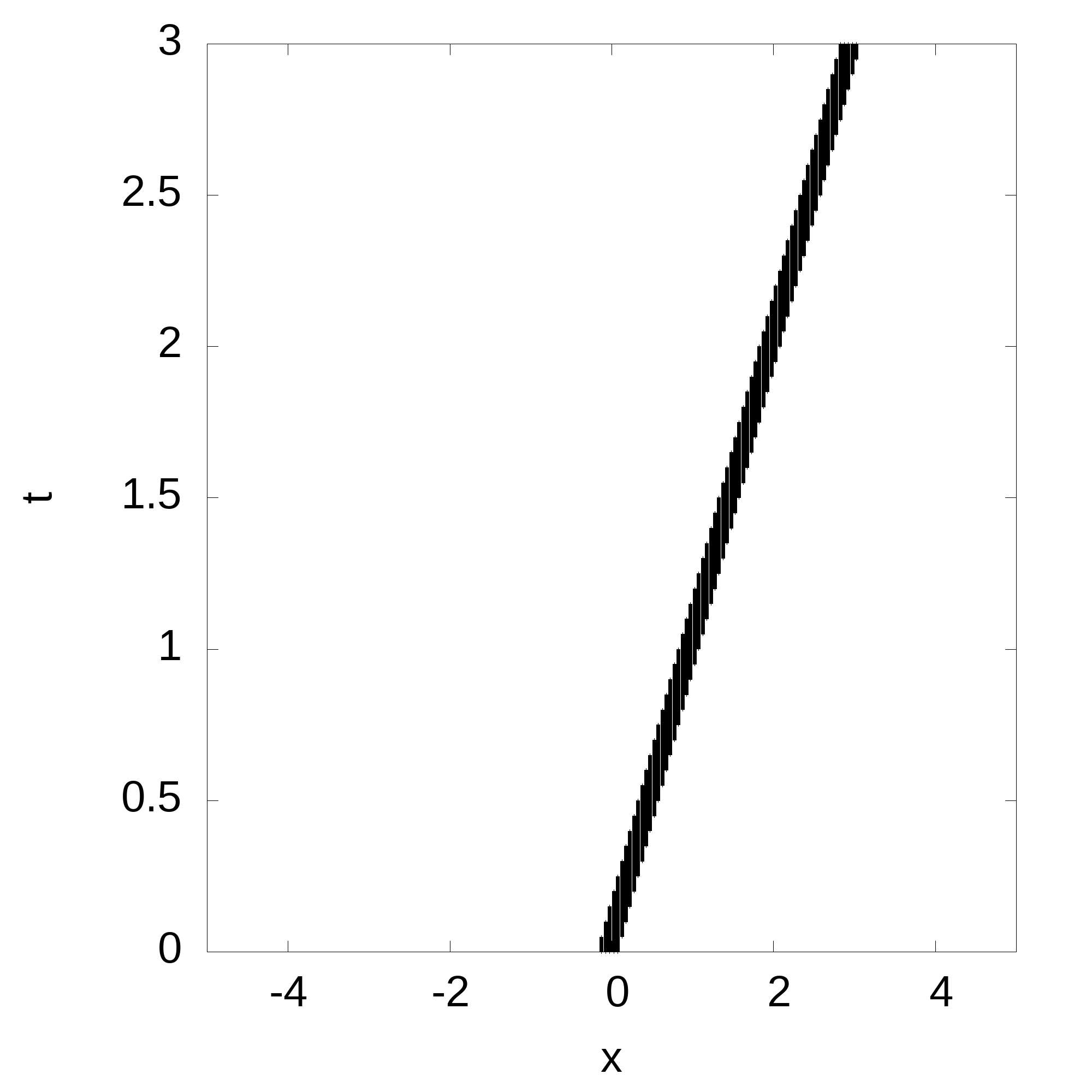}}
  \subfloat[FS1 Indicator]{\label{fig:SCDFS1}\includegraphics[width=0.32\textwidth]{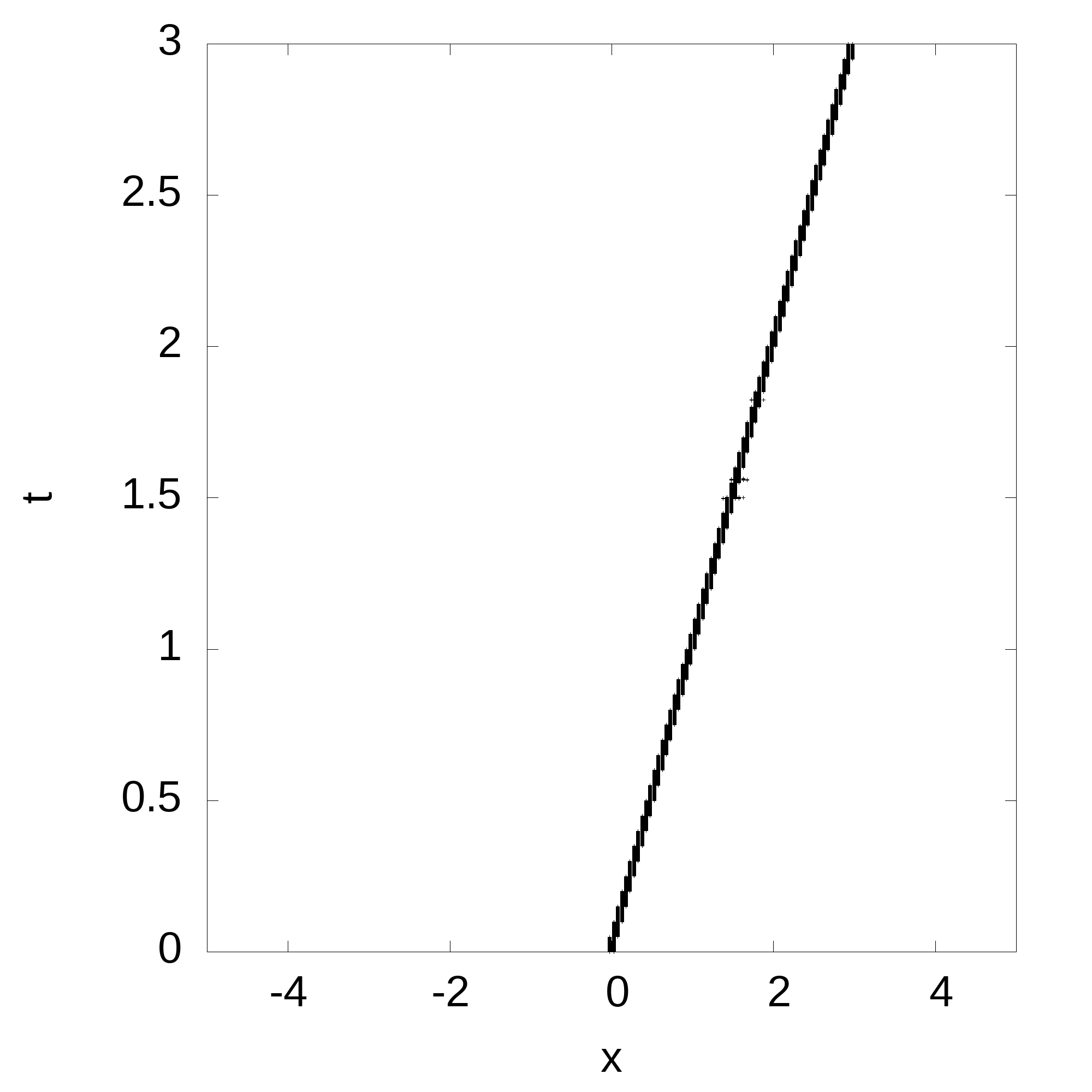}}\hfill
  \subfloat[FS2 Indicator]{\label{fig:SCDFS2}\includegraphics[width=0.32\textwidth]{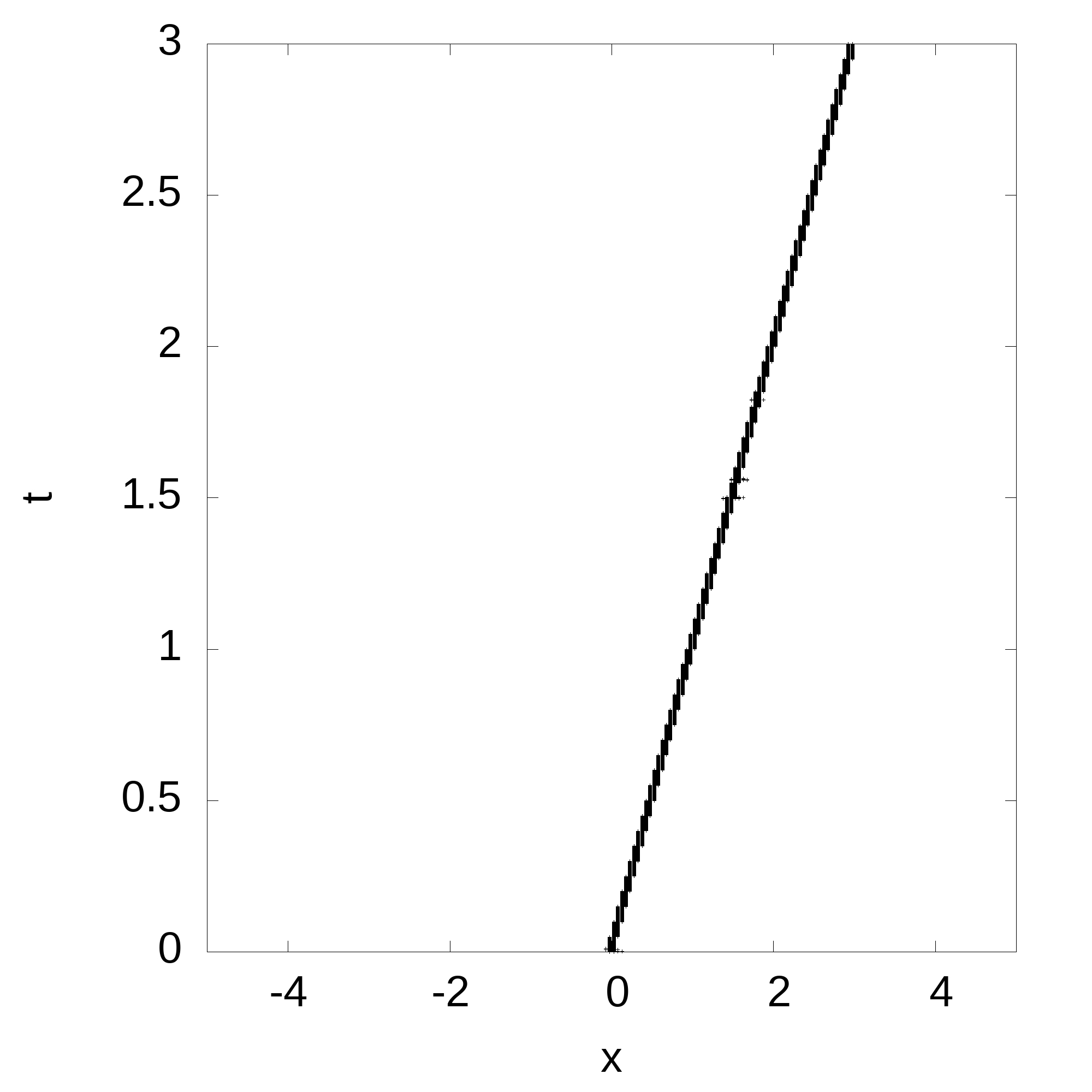}}
  \subfloat[LPR Indicator]{\label{fig:SCDLPR}\includegraphics[width=0.32\textwidth]{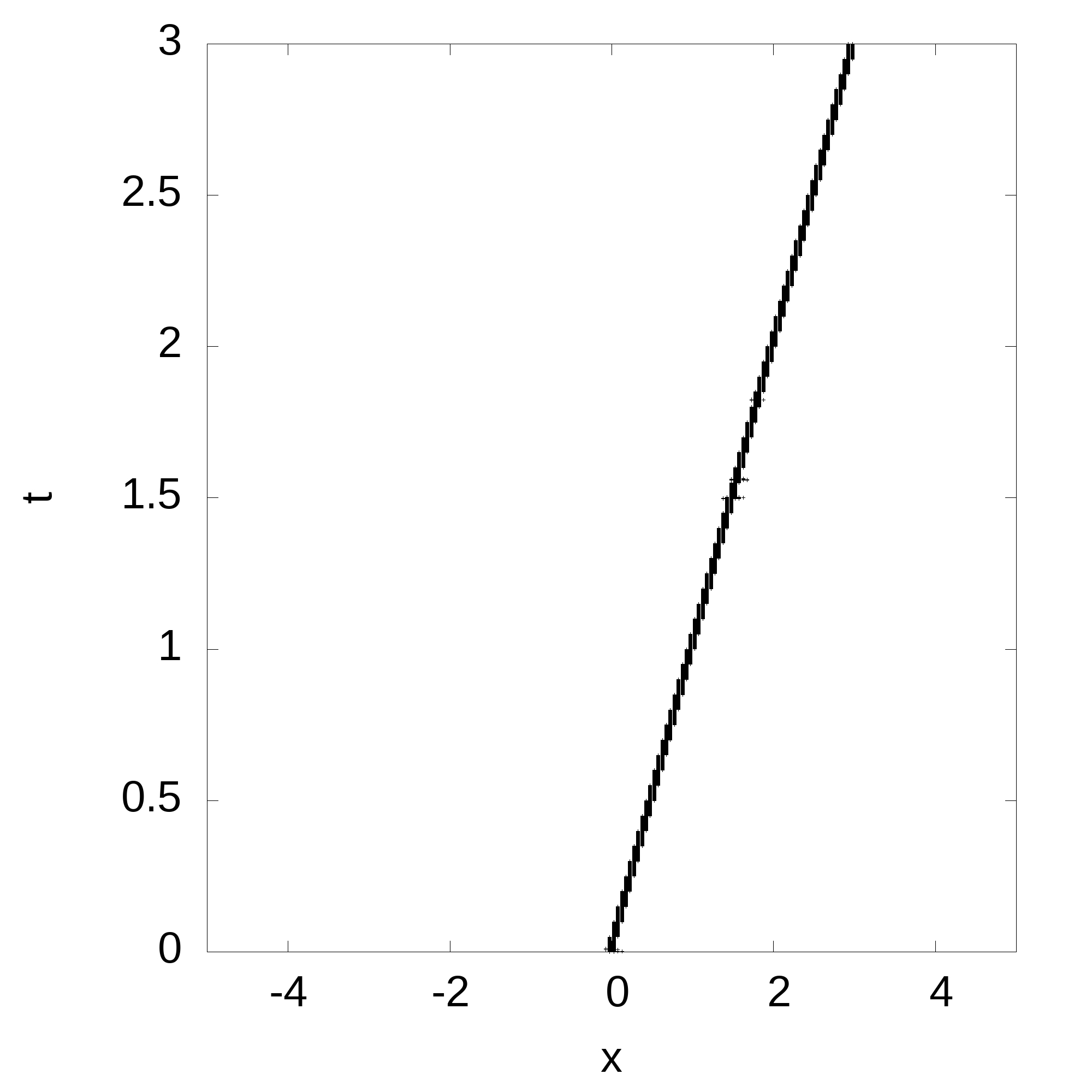}}
  \subfloat[RH Indicator]{\label{fig:SCDRH}\includegraphics[width=0.32\textwidth]{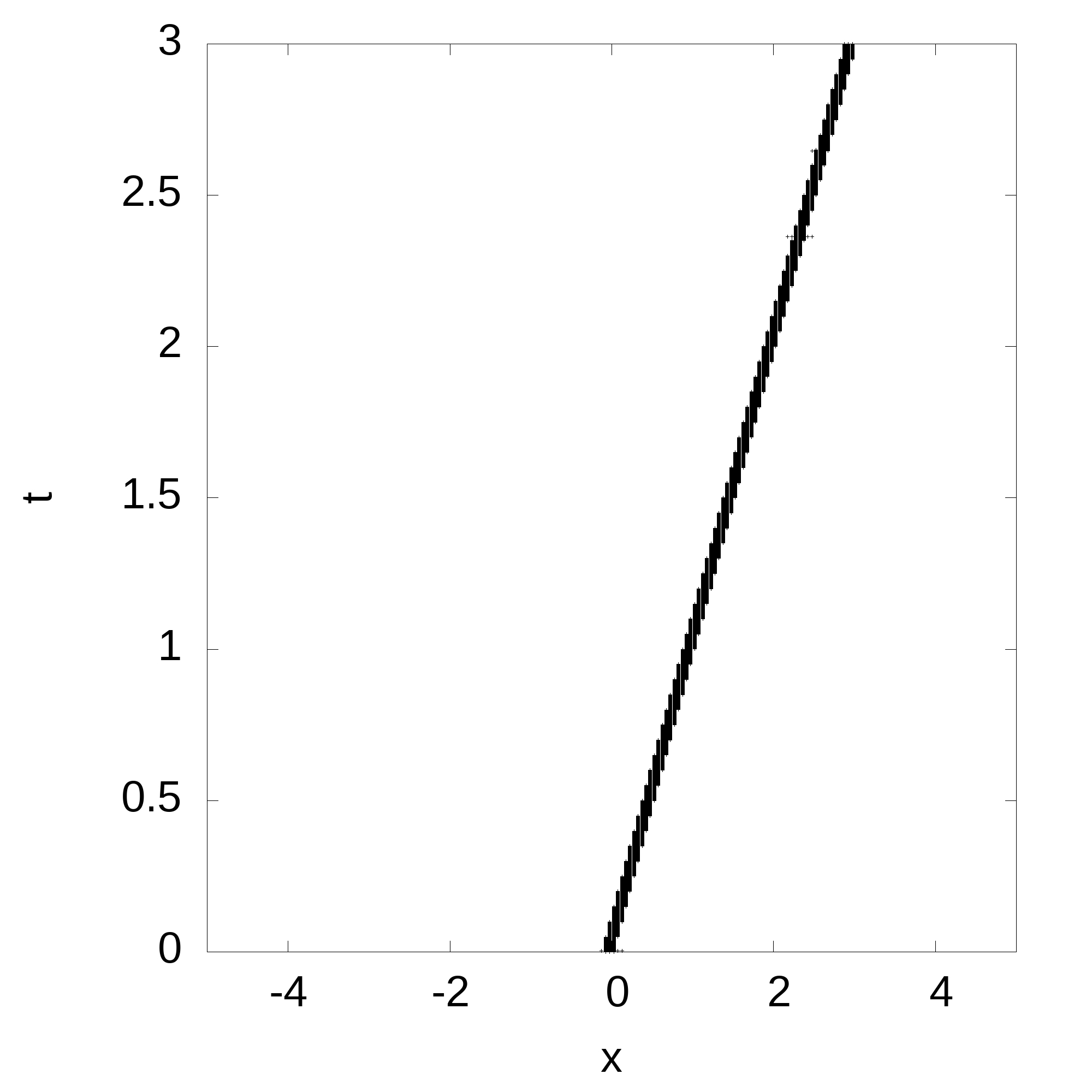}}\hfill
  \subfloat[PPL Indicator]{\label{fig:SCDPPL}\includegraphics[width=0.32\textwidth]{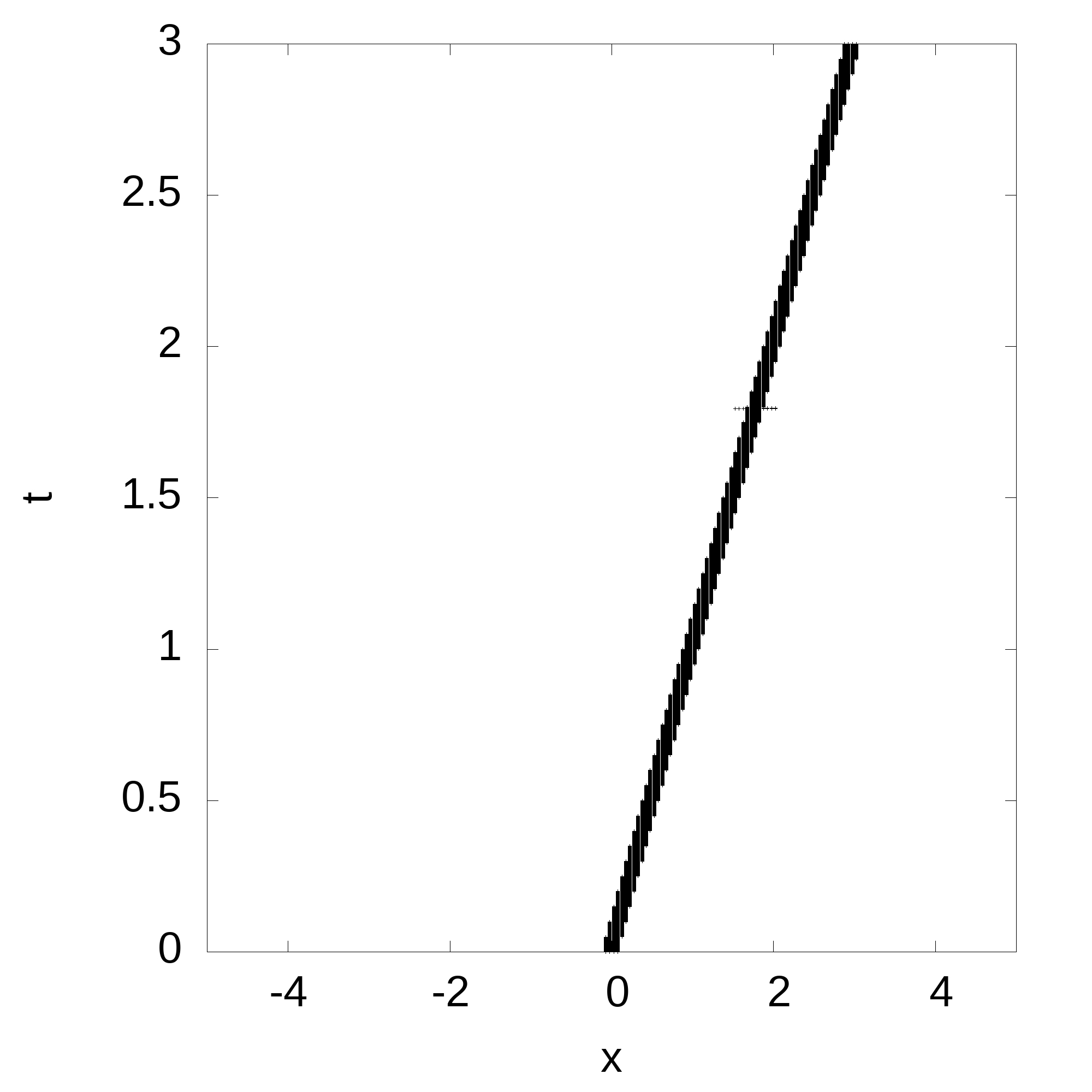}}
  \subfloat[MH Indicator]{\label{fig:SCDMH}\includegraphics[width=0.32\textwidth]{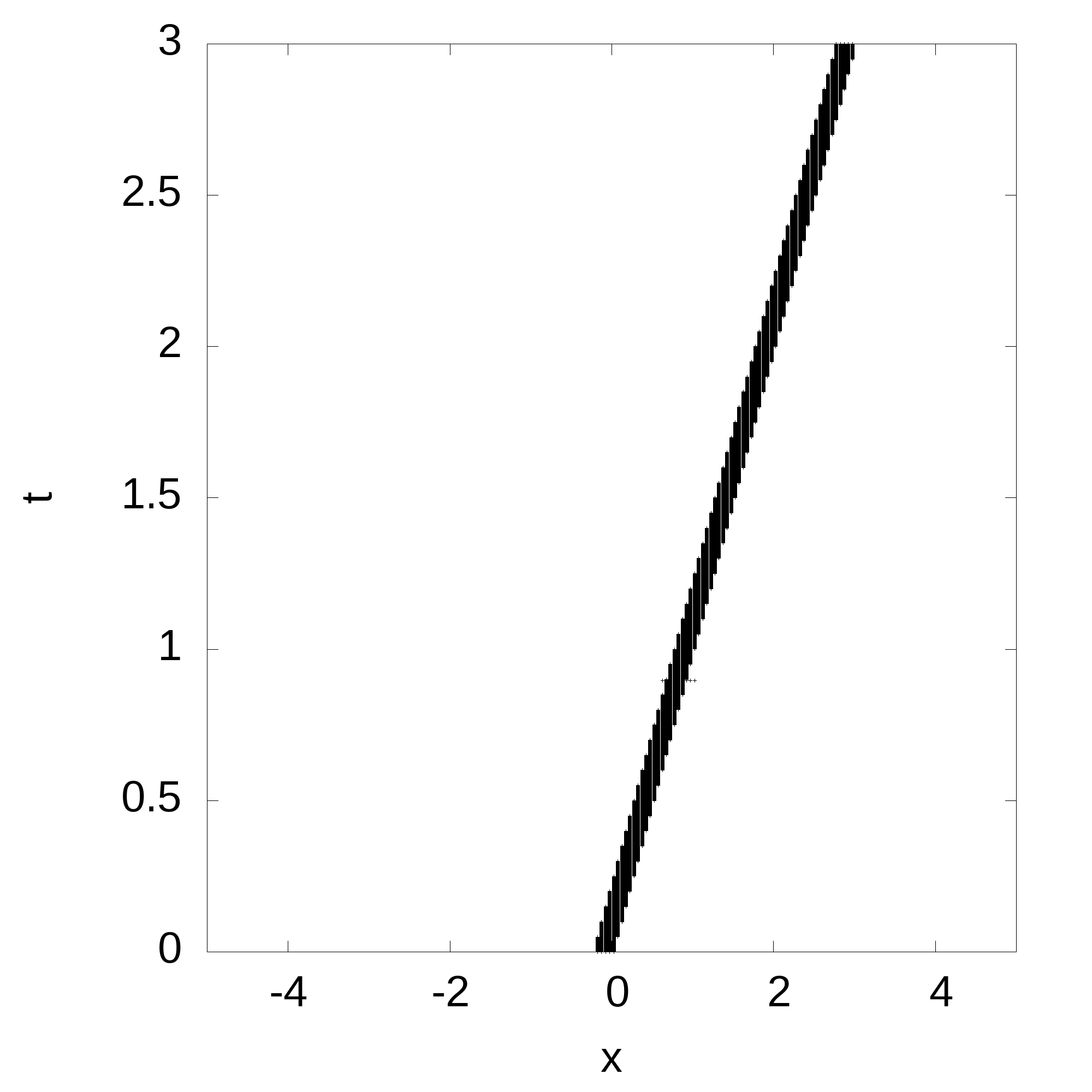}}\hfill
  \caption{The time history of flagged troubled cells of the single contact discontinuity problem for the one-dimensional Euler equations, simulated until $t=3.0$ with 200 elements and $P^{1}$ based DGM}
  \label{fig:SCD}
\end{figure}

\begin{table}[htbp]
%\small
\centering
\begin{tabular}{|c|c|c|c|}
%\hline
%\multicolumn{6}{|c|}{$L_{2}$ error for the Riemann problem configuration-10} \\
\hline
Scheme Indicator & \multicolumn{3}{|c|}{\makecell{Normalised CPU time, normalised\\ with respect to the maximum value}} \\ \cline{2-4}
  & 1D & 2D (Structured) & 2D (Unstructured) \\ 
\hline
PP & 0.974 & 0.967 & 0.965 \\
\hline
SJ & 1.0 & 1.0 & 1.0 \\
\hline
FS1 & 0.966 & 0.962 & 0.965 \\
\hline
FS2 & 0.966 & 0.962 & 0.965 \\
\hline
LPR & 0.985 & 0.983 & 0.991 \\
\hline
RH$^{*}$ & 0.951 & 0.958 & 0.962 \\
\hline
PPL & 0.973 & 0.984 & 0.988 \\
\hline
MH & 0.954 & 0.977 & 0.981 \\
\hline
\end{tabular}
\caption{Normalised CPU time, normalised with respect to the maximum value for various troubled cell indicators. $^{*}$Note that, for the RH indicator, we are excluding the time taken to train the neural network}
\label{table:0New}
\end{table}
%We don't show the solution for each problem as it is quite similar for each case (since we use the same limiter). 
We found the solution to be non-oscillatory using all the troubled cell indicators used. We show the computed solution using a single indicator as the solutions obtained using other indicators are quite similar to each other and stay within an $L^{2}$ norm of $\sim  10^{-14}$ as shown in Tables \ref{table:1New}, \ref{table:2New}, \ref{table:3New}, \ref{table:4New}, \ref{table:5New}, \ref{table:7New}, \ref{table:9New} and \ref{table:11New}. This shows that as long as we use an indicator which identifies all the troubled cells, if we use a good limiter, the accuracy of the solution remains the same. But the additional troubled cells add to the computational time and reduce the efficiency of the scheme.
\\
\\
\noindent \textbf{Test Problem 1 (Single contact discontinuity)\cite{zcq}:} We solve the one-dimensional Euler equations for an ideal gas given by

\begin{equation}\label{1dEulerEquations}
\frac{\partial Q}{\partial t} + \frac{\partial F(Q)}{\partial x} = 0
\end{equation}

\noindent where $Q = (\rho, \rho u, E)^{T}$, $F(Q) = uQ + (0, p, pu)^{T}$, $p = (\gamma -1)(E-\frac{1}{2}\rho u^{2})$ and $\gamma = 1.4$ with the initial condition

\begin{equation}
(\rho, u, p) = 
\begin{cases}
    (\eta_{p},1,1),& \text{if } x\geq 0\\
    (1,1,1),              & \text{otherwise}
\end{cases}
\end{equation}

\noindent in the domain $[-5,5]$. Here, $\rho$ is the density, $u$ is the velocity, $E$ is the total energy and $p$ is the pressure. We choose $\eta_{p}=10^{n}$ where $n=1,\ldots,6$. We only show the results for $\eta_{p}=10^{6}$ as the results for other cases are quite similar. For all the troubled cell indicators, we use the density as the detection variable. Average (over all time steps) and maximum percentages of cells being flagged as troubled cells, for the different troubled-cell indicators, are summarized in Table \ref{table:1} for two different grid sizes and various orders. The computed solution for density obtained at $t=3.0$ using 200 elements while using the SJ indicator and CSWENO limiter for $P^{1}$, $P^{2}$ and $P^{3}$ based DGM is compared and plotted against the exact solution in Figure \ref{fig:CDSolution} on [2.5,3.5]. The error ($|\rho-\rho_{exact}|$) obtained for $P^{1}$, $P^{2}$ and $P^{3}$ based DGM is also plotted in Figure \ref{fig:CDError}. We are not showing the solutions obtained using other indicators as they are quite similar to the solution obtained with the SJ indicator and stay within an $L^{2}$ norm of $\sim  10^{-14}$. We show the $L^{2}$ difference in density between those solutions in table \ref{table:1New}. This clearly shows that as long as we use an indicator which identifies all the troubled cells, if we use a good limiter, the accuracy of the solution remains the same. But the additional troubled cells add to the computational time and reduce the efficiency of the scheme. \\
We also show the time history of the flagged troubled cells using the eight different indicators in Figure \ref{fig:SCD} using $P^{1}$ based DGM for 200 elements. In general, for all the indicators, the number of flagged troubled cells slightly increase when the order increases. From the tabulated results and the figure, we can say that, in this case the troubled cell indicators FS1, FS2 and LPR out perform the other indicators.
\\
\\
\begin{table}[htbp]
\small
\centering
\begin{tabular}{|c|c|c|c|c|c|c|c|c|c|}
%\hline
%\multicolumn{6}{|c|}{$L_{2}$ error for the Riemann problem configuration-10} \\
\hline
 \makecell{No. of \\ Cells} & \makecell{Scheme \\ Indicator} & \multicolumn{2}{|c|}{$P^{1}$} & \multicolumn{2}{|c|}{$P^{2}$} & \multicolumn{2}{|c|}{$P^{3}$} & \multicolumn{2}{|c|}{$P^{4}$} \\
\cline{3-10}
 & & Ave & Max & Ave & Max & Ave & Max & Ave & Max \\
\hline 
\multirow{8}{*}{200} & PP & 3.26 & 6.5 & 4.10 & 6.5 & 4.25 & 7.0 & 4.42 & 7.5 \\
\cline{2-10}
 & SJ & 2.46 & 4.5 & 2.60 & 4.5 & 2.65 & 5.0 & 2.74 & 5.0 \\
\cline{2-10}
 & FS1 & 1.05 & 2.5 & 1.16 & 3.0 & 1.31 & 3.0 & 1.33 & 3.0 \\
\cline{2-10}
 & FS2 & 1.06 & 2.5 & 1.17 & 3.0 & 1.30 & 3.0 & 1.34 & 3.0 \\
\cline{2-10}
 & LPR & 1.02 & 2.0 & 1.18 & 2.5 & 1.27 & 2.5 & 1.35 & 3.0\\
\cline{2-10}
 & RH & 1.45 & 2.5 & 1.56 & 3.0 & 1.67 & 3.5 & 1.76 & 3.5 \\
 \cline{2-10}
 & PPL & 1.95 & 2.5 & 2.05 & 2.5 & 2.10 & 2.5 & 2.15 & 3.0 \\
\cline{2-10}
 & MH & 2.55 & 4.5 & 2.67 & 4.5 & 2.75 & 5.5 & 2.85 & 5.5\\
\hline
\multirow{8}{*}{400} & PP & 3.06 & 5.75 & 4.01 & 6.25 & 4.05 & 6.5 & 3.94 & 7.0 \\
\cline{2-10}
 & SJ & 2.36 & 4.0 & 2.52 & 4.5 & 2.55 & 4.75 & 2.61 & 4.5 \\
\cline{2-10}
 & FS1 & 1.01 & 2.75 & 1.12 & 2.75 & 1.21 & 3.0 & 1.28 & 3.0 \\
\cline{2-10}
 & FS2 & 1.03 & 2.75 & 1.14 & 2.75 & 1.31 & 3.0 & 1.29 & 3.0 \\
\cline{2-10}
 & LPR & 0.96 & 2.25 & 1.15 & 2.5 & 1.22 & 2.5 & 1.32 & 3.0\\
\cline{2-10}
 & RH & 1.42 & 2.5 & 1.52 & 2.75 & 1.63 & 3.0 & 1.77 & 3.25 \\
 \cline{2-10}
 & PPL & 1.91 & 2.25 & 2.04 & 2.75 & 2.11 & 2.75 & 2.22 & 3.0 \\
\cline{2-10}
 & MH & 2.51 & 4.25 & 2.61 & 4.75 & 2.72 & 5.0 & 2.84 & 5.25\\
\hline
\end{tabular}
\caption{Average (marked as Ave) and maximum (marked as Max) percentages of cells flagged as troubled cells subject to different troubled-cell indicators for the single contact discontinuity for various orders and two grid sizes.}
\label{table:1}
\end{table}

\begin{table}[htbp]
%\small
\centering
\begin{tabular}{|c|c|c|c|c|}
%\hline
%\multicolumn{6}{|c|}{$L_{2}$ error for the Riemann problem configuration-10} \\
\hline
Scheme Indicator & \multicolumn{4}{|c|}{\makecell{$L^{2}$ difference in density in comparison with \\ SJ indicator solution for various orders}} \\ \cline{2-5}
  & $P^{1}$ & $P^{2}$ & $P^{3}$ & $P^{4}$ \\ 
\hline
PP & 1.6E-14 & 1.1E-15 & 2.3E-15 & 3.5E-15\\
\hline
FS1 & 2.1E-14 & 3.7E-15 & 5.4E-15 & 3.6E-15\\
\hline
FS2 & 2.1E-14 & 3.7E-15 & 5.4E-15 & 3.6E-15\\
\hline
LPR & 1.9E-15 & 2.4E-15 & 3.2E-15 & 2.6E-15\\
\hline
RH & 2.5E-14 & 9.6E-15 & 2.7E-15 & 5.4E-15\\
\hline
PPL & 4.8E-14 & 8.2E-15 & 1.9E-15 & 8.3E-15\\
\hline
MH & 5.3E-15 & 6.5E-15 & 3.9E-15 & 1.2E-15\\
\hline
\end{tabular}
\caption{$L^{2}$ difference in density using CSWENO limiter and SJ indicator (shown in Figure \ref{fig:CDCSWENO}) and the solution obtained using other indicators for the moving contact discontinuity problem using 200 elements between solution obtained}
\label{table:1New}
\end{table}

\begin{figure}[htbp]
  \centering
  \subfloat[Solution of 1D Euler equations for the Sod problem with $P^{1}$, $P^{2}$ and $P^{3}$ based DGM and SJ indicator]{\label{fig:SodSolution}\includegraphics[width=0.48\textwidth]{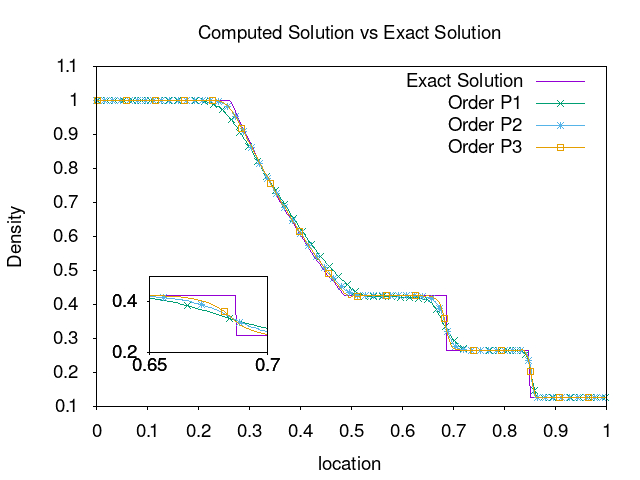}}
  \subfloat[Density error ($|\rho-\rho_{exact}|$) for $P^{1}$, $P^{2}$ and $P^{3}$ based DGM]{\label{fig:SodError}\includegraphics[width=0.48\textwidth]{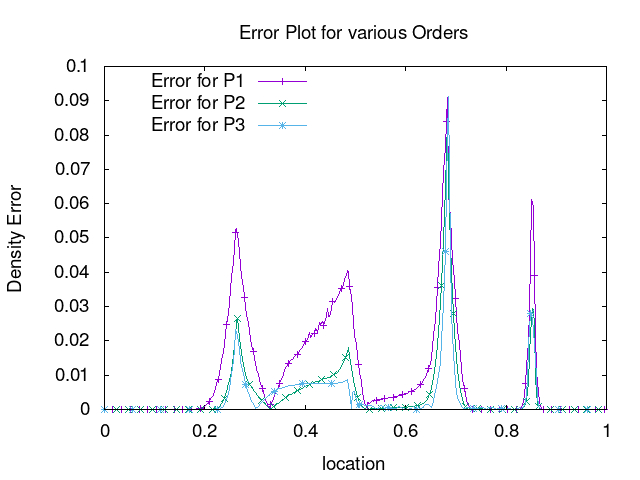}}\hfill
  \caption{Comparison of density solutions of Sod Problem at $t=0.2$ using 200 elements obtained with the $P^{1}$, $P^{2}$ and $P^{3}$ based DGM and the exact solution. Density error ($|\rho-\rho_{exact}|$) for $P^{1}$, $P^{2}$ and $P^{3}$ based DGM is also plotted. Figure \ref{fig:SodSolution} also includes a zoomed in portion of the solution for better comparison}
  \label{fig:SodCSWENO}
\end{figure}

\begin{figure}[htbp]
  \centering
  \subfloat[PP Indicator]{\label{fig:SodPP}\includegraphics[width=0.32\textwidth]{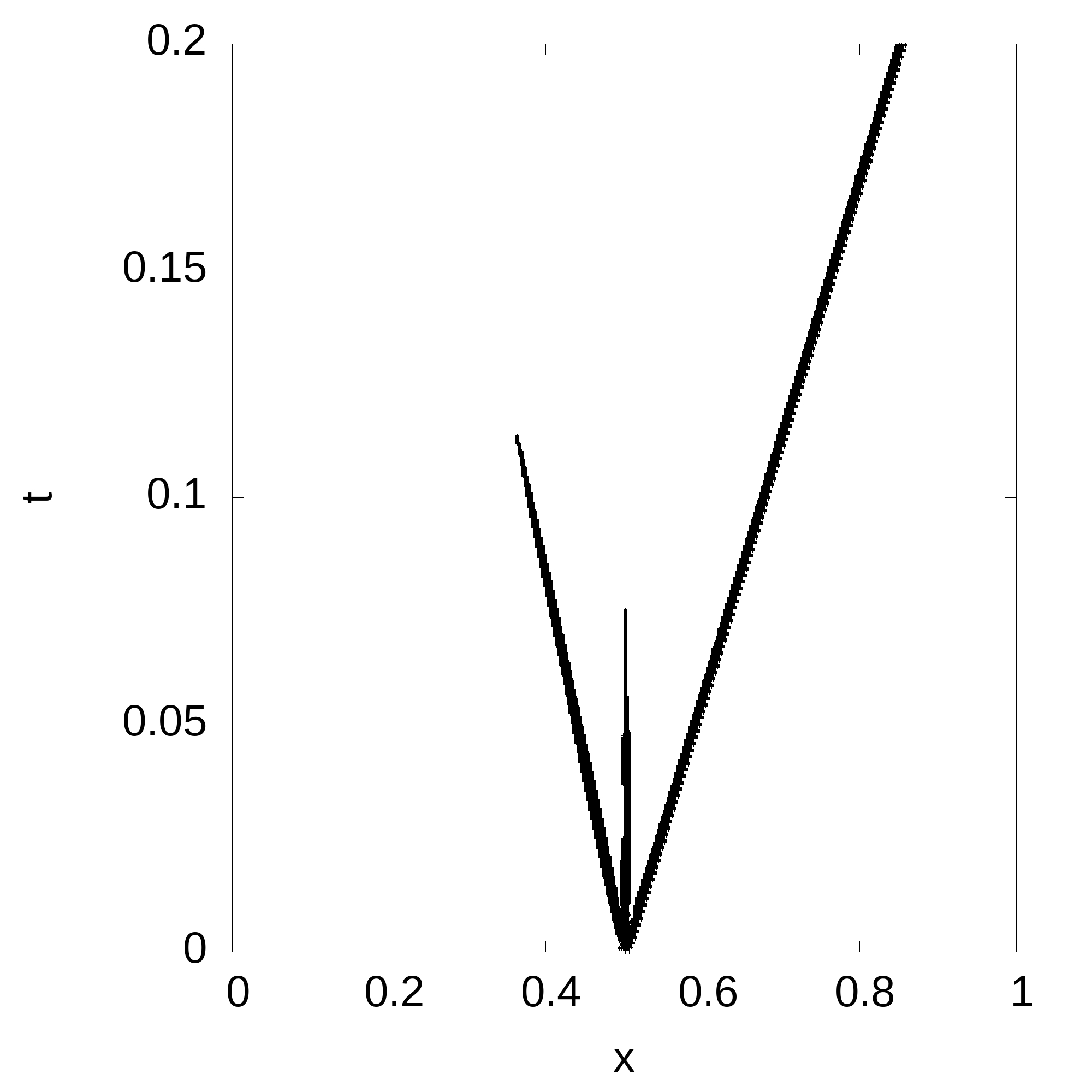}}
  \subfloat[SJ Indicator]{\label{fig:SodSJ}\includegraphics[width=0.32\textwidth]{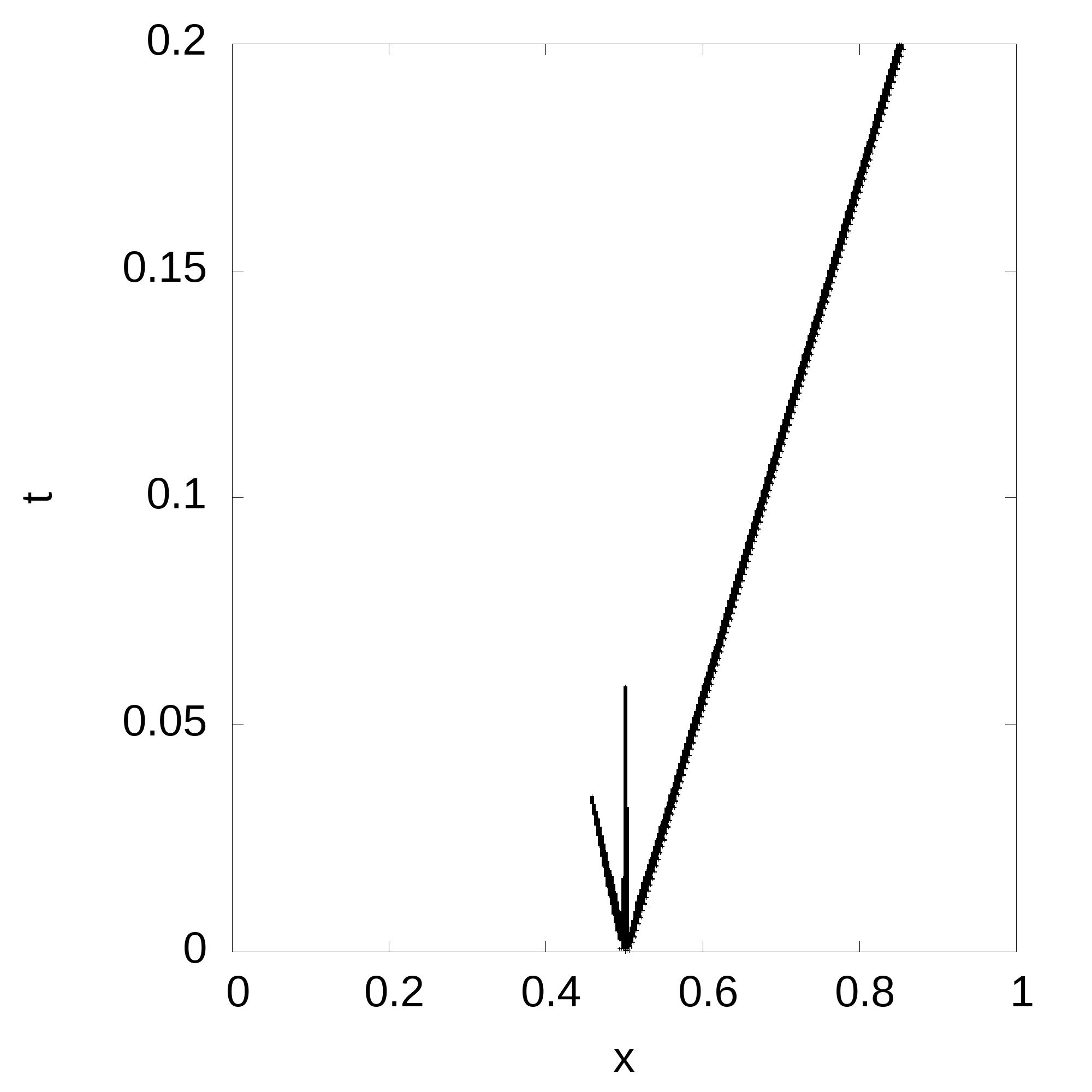}}
  \subfloat[FS1 Indicator]{\label{fig:SodFS1}\includegraphics[width=0.32\textwidth]{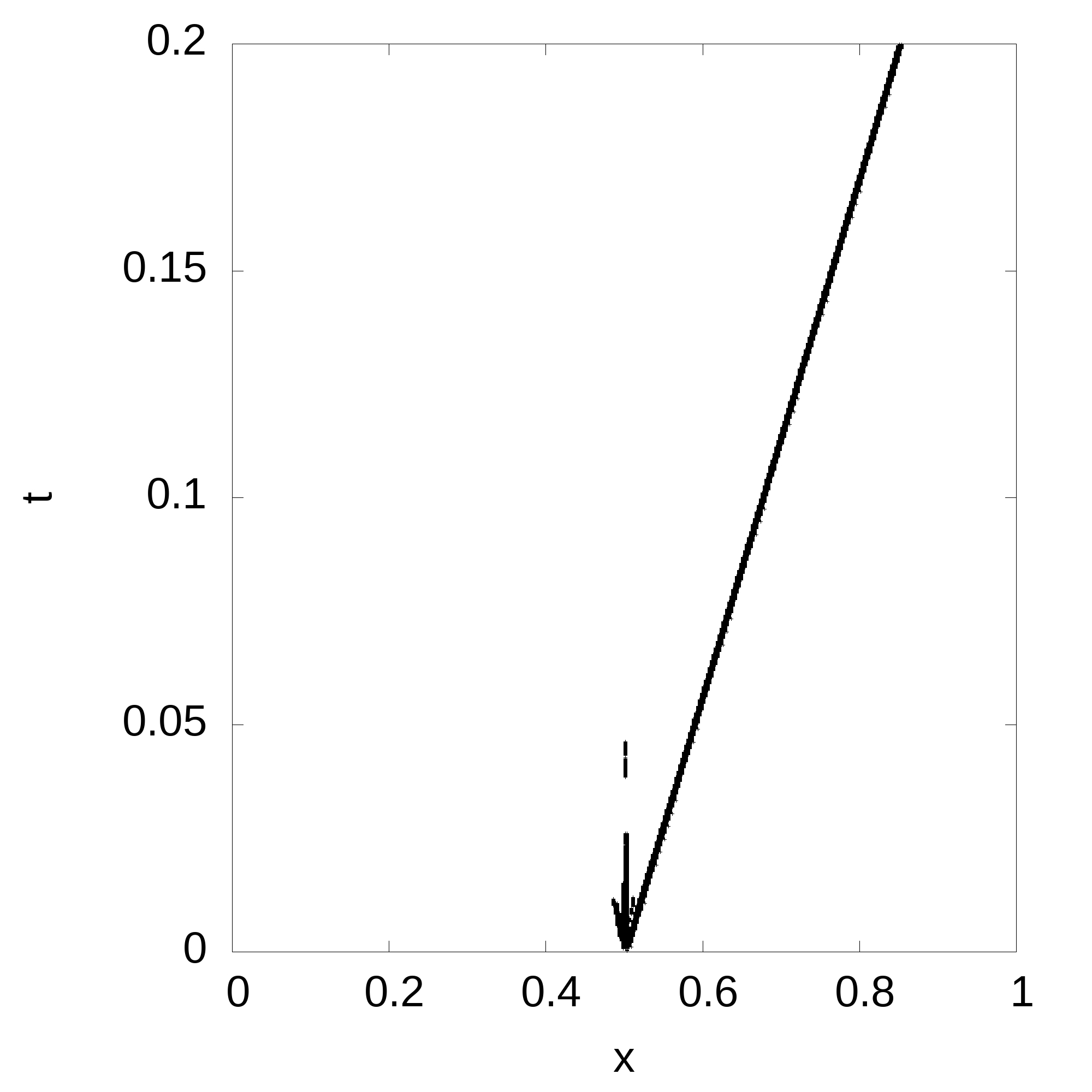}}\hfill
  \subfloat[FS2 Indicator]{\label{fig:SodFS2}\includegraphics[width=0.32\textwidth]{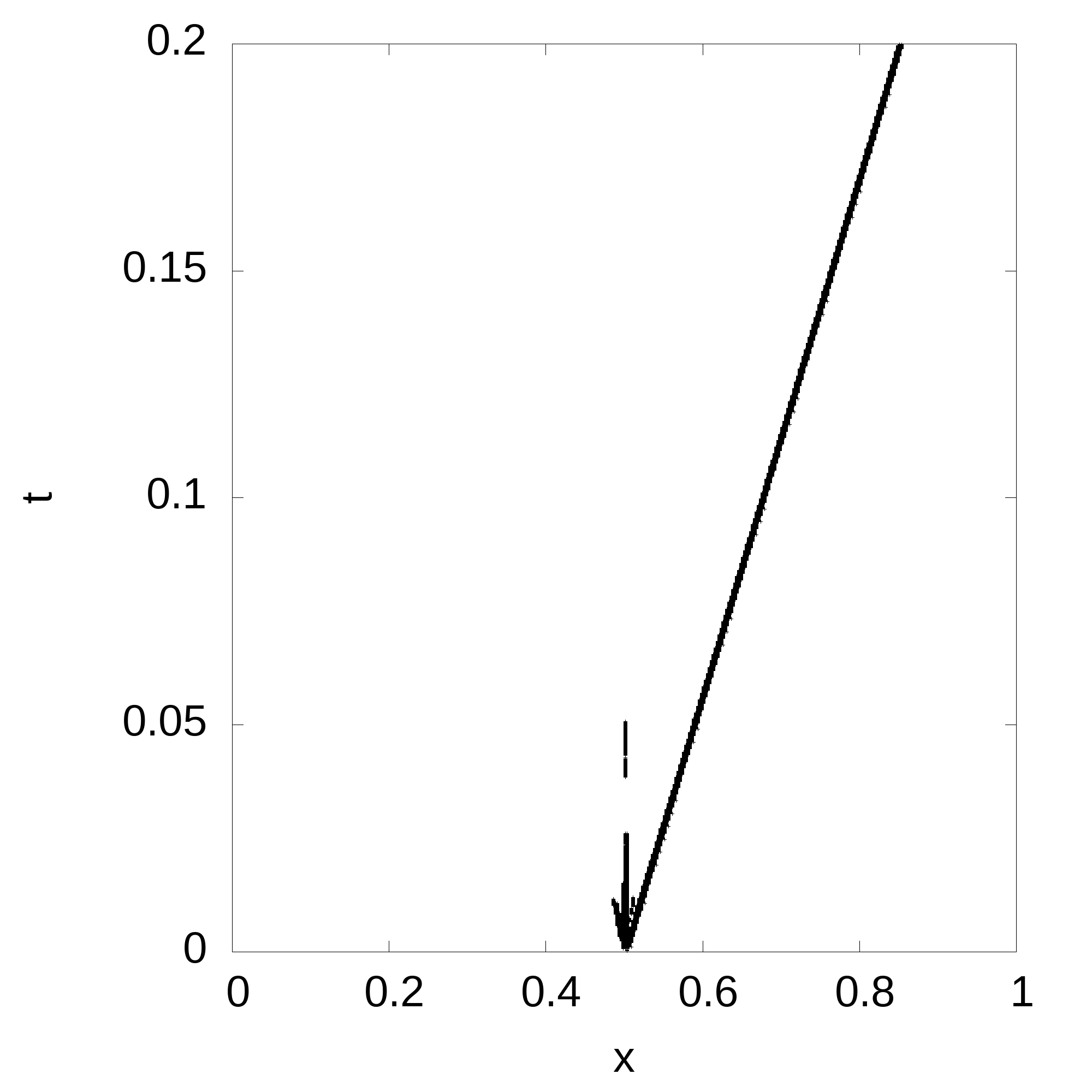}}
  \subfloat[LPR Indicator]{\label{fig:SodLPR}\includegraphics[width=0.32\textwidth]{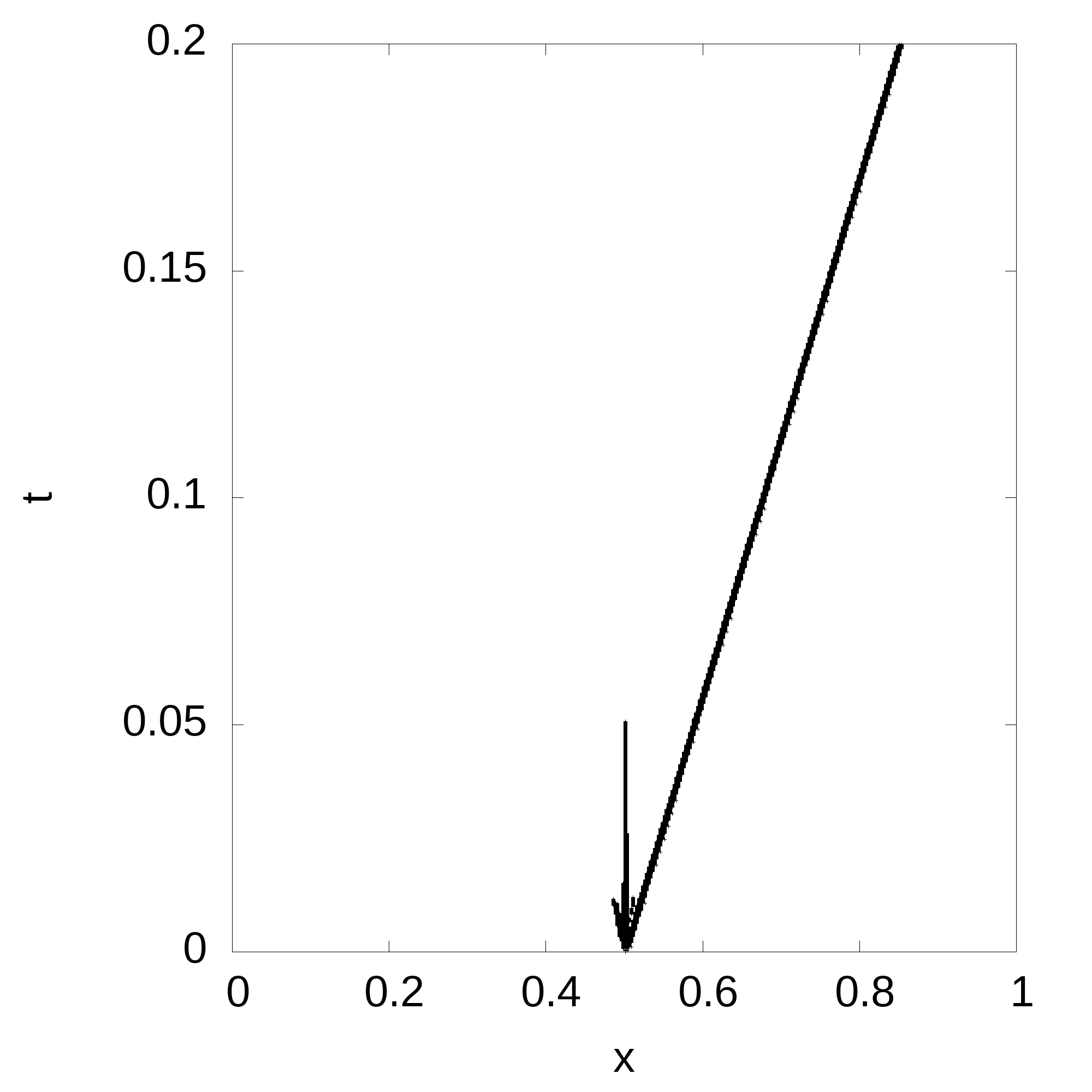}}
  \subfloat[RH Indicator]{\label{fig:SodRH}\includegraphics[width=0.32\textwidth]{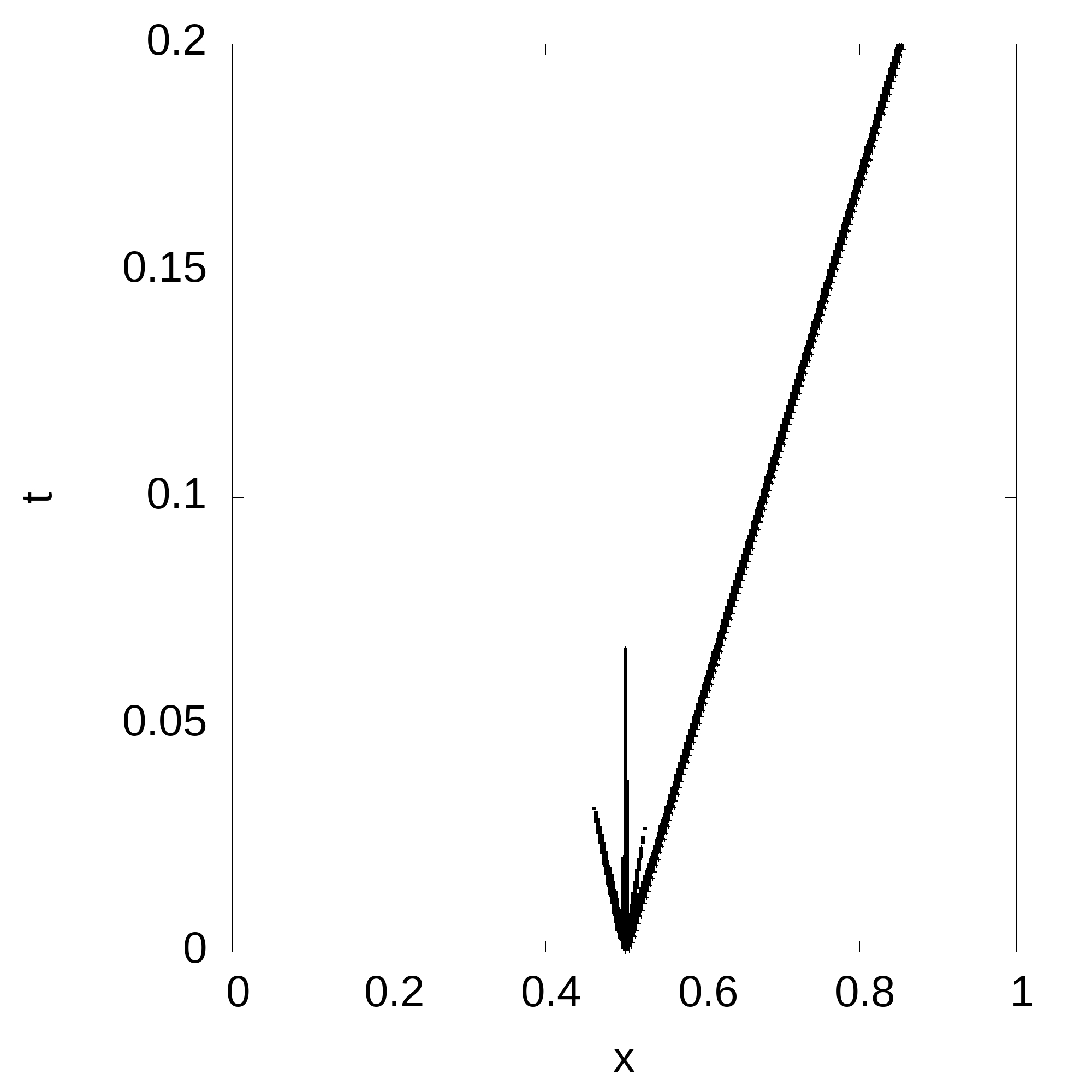}}\hfill
  \subfloat[PPL Indicator]{\label{fig:SodPPL}\includegraphics[width=0.32\textwidth]{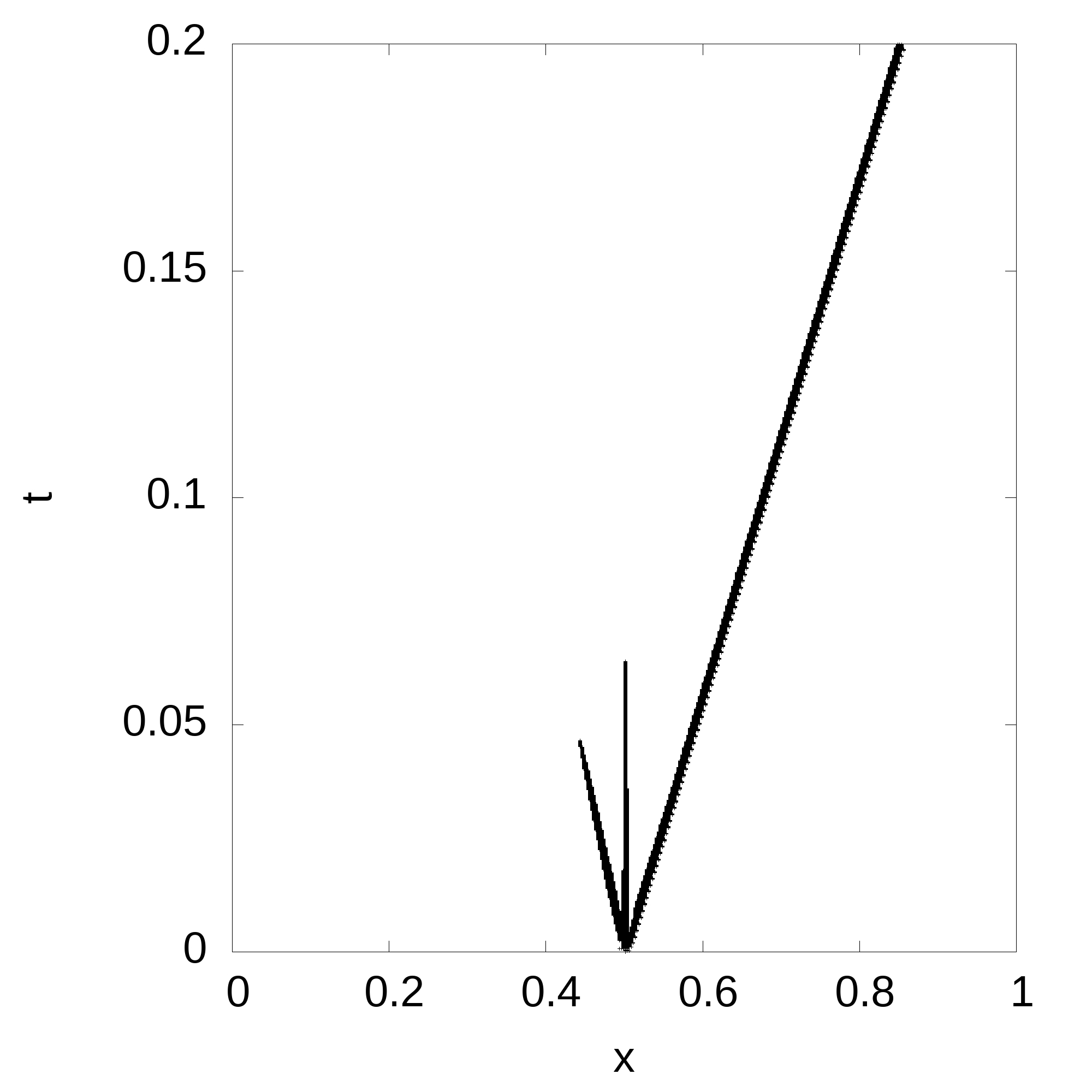}}
  \subfloat[MH Indicator]{\label{fig:SodMH}\includegraphics[width=0.32\textwidth]{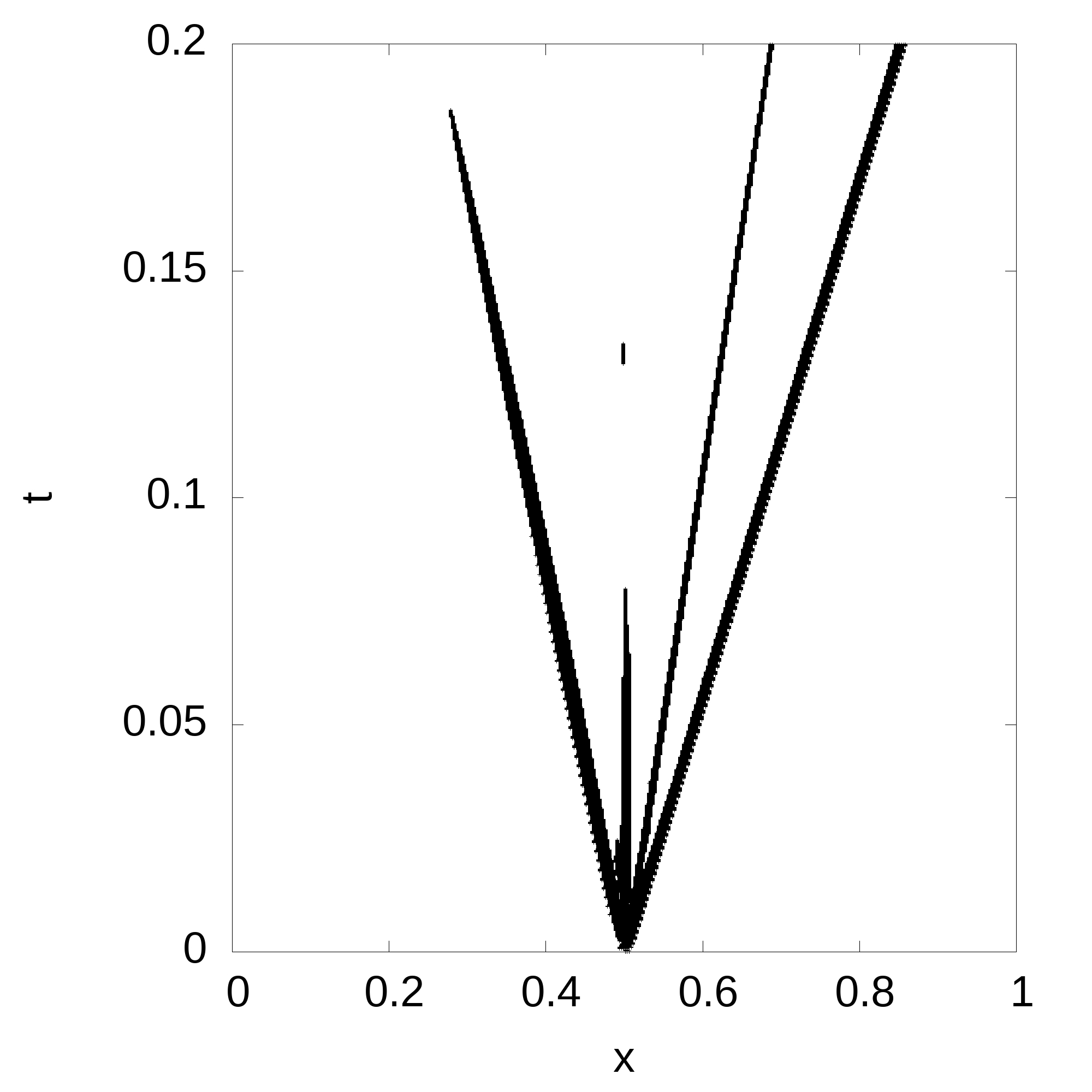}}\hfill
  \caption{The time history of flagged troubled cells of the Sod test problem for the one-dimensional Euler equations, simulated until $t=0.2$ with 200 elements and $P^{1}$ based DGM}
  \label{fig:Sod}
\end{figure}

\noindent \textbf{Test Problem 2 (Sod Problem)\cite{sod}:} We solve the one-dimensional Euler equations as given by \eqref{1dEulerEquations} with the initial condition

\begin{equation}
(\rho, u, p) = 
\begin{cases}
    (1,0,1),& \text{if } x< 0.5\\
    (0.125,0,0.1),              & \text{otherwise}
\end{cases}
\end{equation}

\noindent in the domain $[0,1]$. For all the troubled cell indicators, we again use the density as the detection variable. Average (over all time steps) and maximum percentages of cells being flagged as troubled cells, for the different troubled-cell indicators, are summarized in Table \ref{table:2} for two different grid sizes and various orders. The computed solution for density obtained at $t=0.2$ using 200 elements while using the SJ indicator and CSWENO limiter for $P^{1}$, $P^{2}$ and $P^{3}$ based DGM is compared and plotted against the exact solution in Figure \ref{fig:SodSolution}. The error ($|\rho-\rho_{exact}|$) obtained for $P^{1}$, $P^{2}$ and $P^{3}$ based DGM is also plotted in Figure \ref{fig:SodError}. We are not showing the solutions obtained using other indicators as they are quite similar to the solution obtained with the SJ indicator and stay within an $L^{2}$ norm of $\sim  10^{-14}$. We show the $L^{2}$ difference in density between those solutions in table \ref{table:2New}.
\\
We also show the time history of the flagged troubled cells using the eight different indicators in Figure \ref{fig:Sod} using $P^{1}$ based DGM for 200 elements. For all the indicators, again, the number of flagged troubled cells slightly increase when the order increases. From the tabulated results and the figure, the troubled cell indicators SJ, PP and PPL perform in a similar fashion, while the MH indicator is quite bad. Again, we can say that, even for this test problem, the troubled cell indicators FS1, FS2 and LPR out perform the other indicators.
\\
\\
\begin{table}[htbp]
\small
\centering
\begin{tabular}{|c|c|c|c|c|c|c|c|c|c|}
%\hline
%\multicolumn{6}{|c|}{$L_{2}$ error for the Riemann problem configuration-10} \\
\hline
 \makecell{No. of \\ Cells} & \makecell{Scheme \\ Indicator} & \multicolumn{2}{|c|}{$P^{1}$} & \multicolumn{2}{|c|}{$P^{2}$} & \multicolumn{2}{|c|}{$P^{3}$} & \multicolumn{2}{|c|}{$P^{4}$} \\
\cline{3-10}
 & & Ave & Max & Ave & Max & Ave & Max & Ave & Max \\
\hline 
\multirow{8}{*}{200} & PP & 4.74 & 5.5 & 4.85 & 5.5 & 5.35 & 6.5 & 6.31 & 7.5 \\
\cline{2-10}
 & SJ & 4.02 & 4.5 & 4.58 & 5.5 & 5.09 & 5.5 & 5.66 & 6.0 \\
\cline{2-10}
 & FS1 & 1.02 & 2.5 & 1.38 & 3.5 & 1.47 & 3.5 & 1.68 & 4.0 \\
\cline{2-10}
 & FS2 & 1.09 & 2.5 & 1.44 & 3.5 & 1.54 & 3.5 & 1.75 & 4.0 \\
\cline{2-10}
 & LPR & 1.10 & 2.5 & 1.88 & 3.0 & 2.12 & 4.5 & 2.55 & 5.0\\
\cline{2-10}
 & RH & 2.45 & 3.5 & 2.76 & 4.0 & 3.42 & 4.5 & 3.96 & 5.0 \\
 \cline{2-10}
 & PPL & 4.05 & 4.5 & 4.62 & 5.5 & 5.11 & 5.5 & 5.75 & 6.0 \\
\cline{2-10}
 & MH & 8.55 & 12.5 & 9.77 & 12.5 & 10.27 & 13.5 & 10.93 & 14.0\\
\hline
\multirow{8}{*}{400} & PP & 4.25 & 5.0 & 4.53 & 5.5 & 5.14 & 6.5 & 6.02 & 7.0 \\
\cline{2-10}
 & SJ & 3.88 & 4.25 & 4.22 & 4.5 & 4.71 & 5.25 & 5.16 & 5.25 \\
\cline{2-10}
 & FS1 & 0.64 & 1.25 & 0.69 & 1.25 & 0.75 & 1.5 & 0.81 & 2.0 \\
\cline{2-10}
 & FS2 & 0.65 & 1.25 & 0.72 & 1.25 & 0.81 & 1.5 & 0.85 & 2.25 \\
\cline{2-10}
 & LPR & 0.99 & 2.5 & 1.73 & 2.75 & 1.92 & 3.5 & 2.22 & 4.0\\
\cline{2-10}
 & RH & 2.25 & 3.5 & 2.55 & 4.25 & 3.13 & 4.75 & 3.85 & 5.0 \\
 \cline{2-10}
 & PPL & 3.94 & 4.25 & 4.33 & 4.75 & 4.84 & 5.25 & 5.27 & 5.5 \\
\cline{2-10}
 & MH & 8.25 & 12.25 & 9.51 & 12.5 & 10.04 & 13.25 & 10.81 & 13.75\\
\hline
\end{tabular}
\caption{Average (marked as Ave) and maximum (marked as Max) percentages of cells flagged as troubled cells subject to different troubled-cell indicators for the Sod test problem for various orders and two grid sizes.}
\label{table:2}
\end{table}

\begin{table}[htbp]
%\small
\centering
\begin{tabular}{|c|c|c|c|c|}
%\hline
%\multicolumn{6}{|c|}{$L_{2}$ error for the Riemann problem configuration-10} \\
\hline
Scheme Indicator & \multicolumn{4}{|c|}{\makecell{$L^{2}$ difference in density in comparison with \\ SJ indicator solution for various orders}} \\ \cline{2-5}
  & $P^{1}$ & $P^{2}$ & $P^{3}$ & $P^{4}$ \\ 
\hline
PP & 2.3E-14 & 4.1E-15 & 6.8E-15 & 6.2E-15\\
\hline
FS1 & 1.8E-14 & 5.3E-15 & 7.4E-15 & 7.9E-15\\
\hline
FS2 & 1.8E-14 & 5.3E-15 & 7.4E-15 & 7.9E-15\\
\hline
LPR & 7.6E-15 & 1.8E-15 & 4.9E-15 & 3.6E-15\\
\hline
RH & 1.3E-14 & 8.3E-15 & 1.8E-15 & 4.5E-15\\
\hline
PPL & 4.2E-15 & 9.1E-15 & 5.8E-15 & 5.7E-15\\
\hline
MH & 6.4E-15 & 2.3E-15 & 4.1E-15 & 2.9E-15\\
\hline
\end{tabular}
\caption{$L^{2}$ difference in density using CSWENO limiter and SJ indicator (shown in Figure \ref{fig:SodCSWENO}) and the solution obtained using other indicators for the Sod problem using 200 elements between solution obtained}
\label{table:2New}
\end{table}

\begin{figure}[htbp]
  \centering
  \subfloat[Solution of 1D Euler equations for the Lax problem with $P^{1}$, $P^{2}$ and $P^{3}$ based DGM]{\label{fig:LaxSolution}\includegraphics[width=0.48\textwidth]{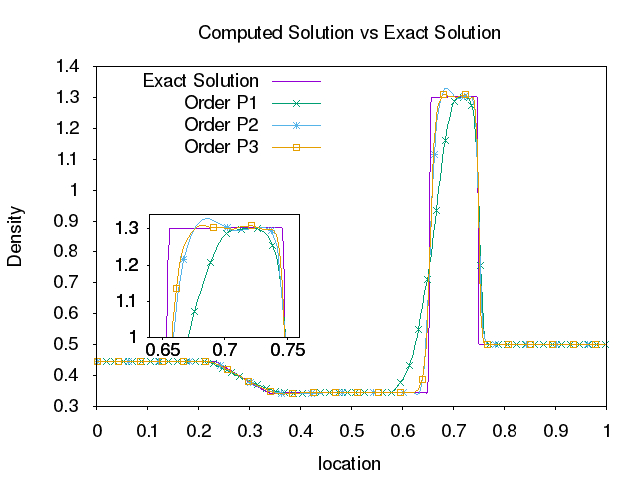}}
  \subfloat[Density error ($|\rho-\rho_{exact}|$) for $P^{1}$, $P^{2}$ and $P^{3}$ based DGM]{\label{fig:LaxError}\includegraphics[width=0.48\textwidth]{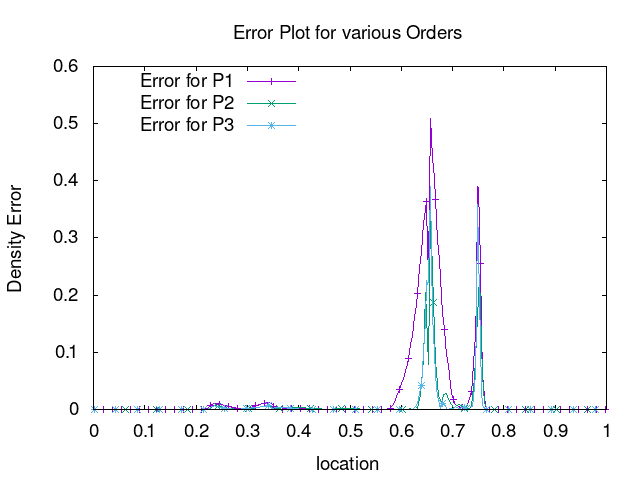}}\hfill
  \caption{Comparison of density solutions of Lax Problem at $t=0.1$ using 200 elements obtained with the $P^{1}$, $P^{2}$ and $P^{3}$ based DGM and the exact solution. Density error ($|\rho-\rho_{exact}|$) for $P^{1}$, $P^{2}$ and $P^{3}$ based DGM is also plotted. Figure \ref{fig:LaxSolution} also includes a zoomed in portion of the solution for better comparison}
  \label{fig:LaxCSWENO}
\end{figure}

\begin{figure}[htbp]
  \centering
  \subfloat[PP Indicator]{\label{fig:LaxPP}\includegraphics[width=0.32\textwidth]{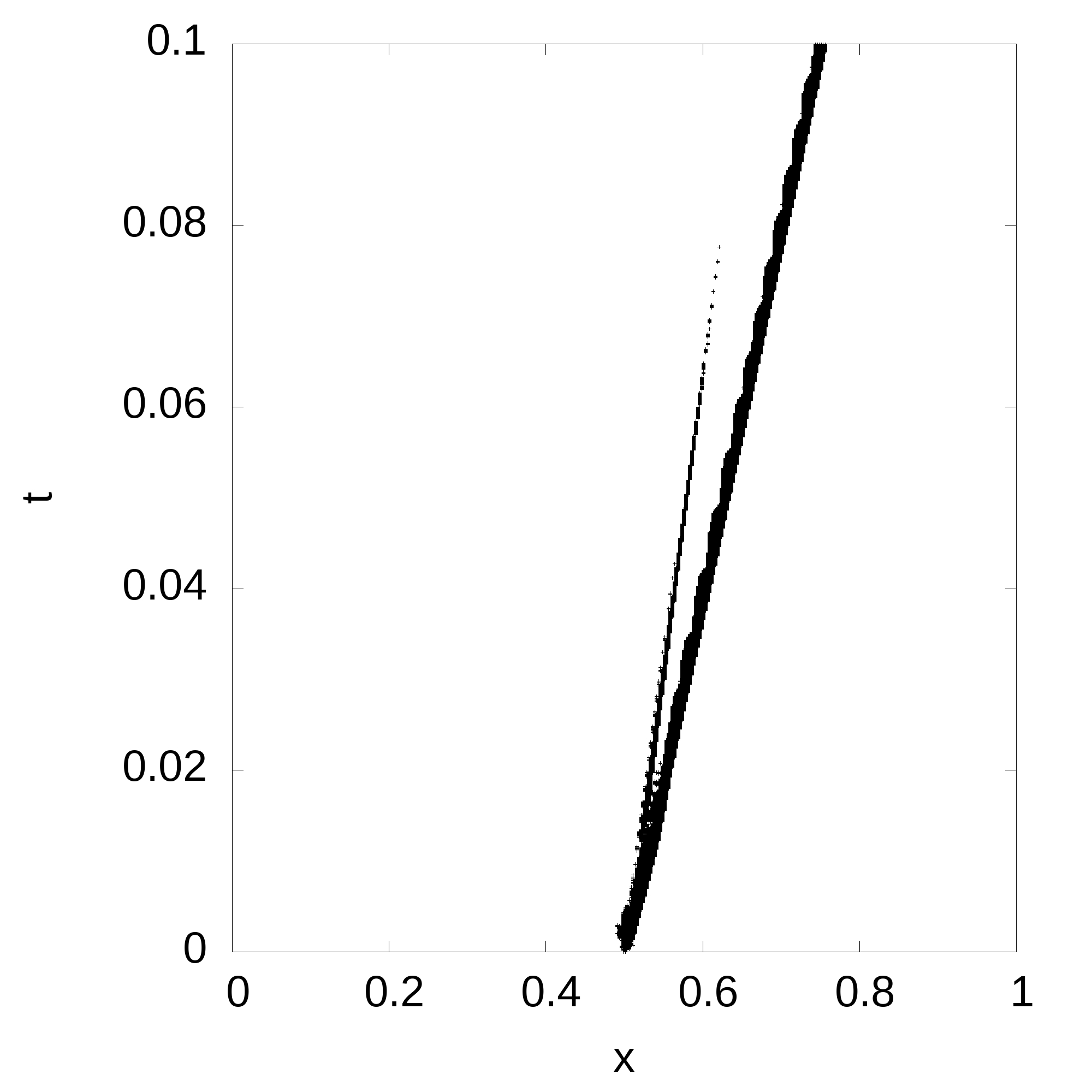}}
  \subfloat[SJ Indicator]{\label{fig:LaxSJ}\includegraphics[width=0.32\textwidth]{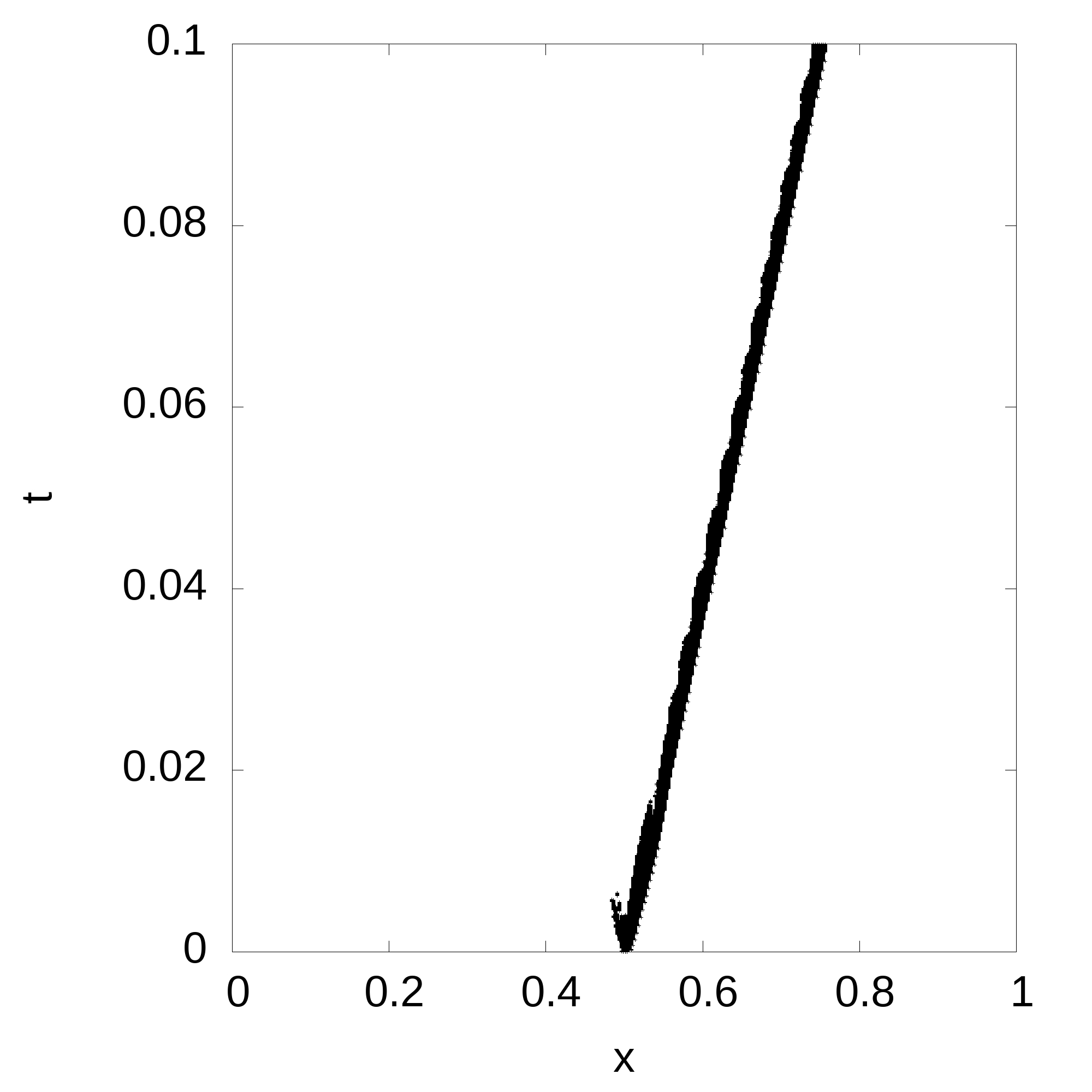}}
  \subfloat[FS1 Indicator]{\label{fig:LaxFS1}\includegraphics[width=0.32\textwidth]{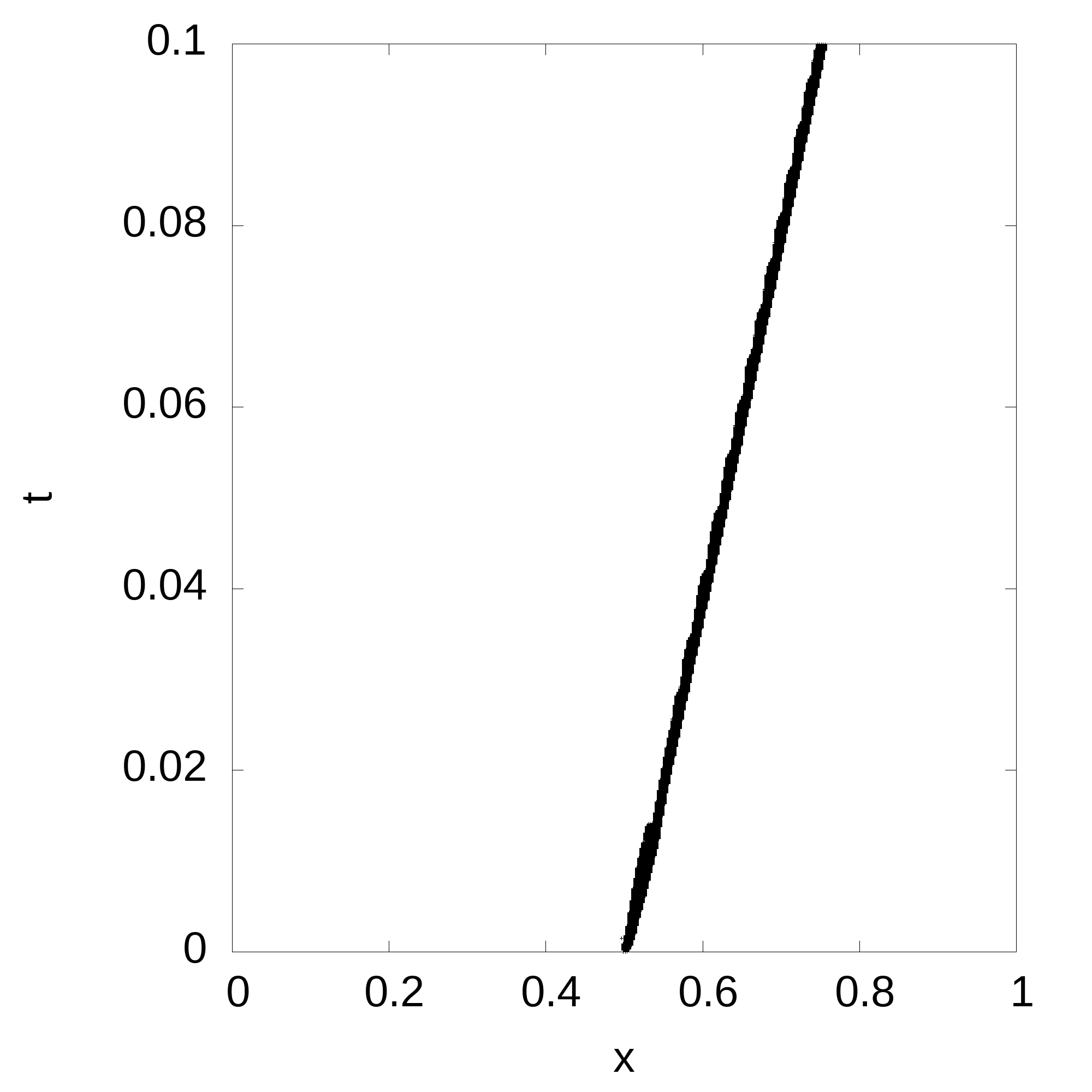}}\hfill
  \subfloat[FS2 Indicator]{\label{fig:LaxFS2}\includegraphics[width=0.32\textwidth]{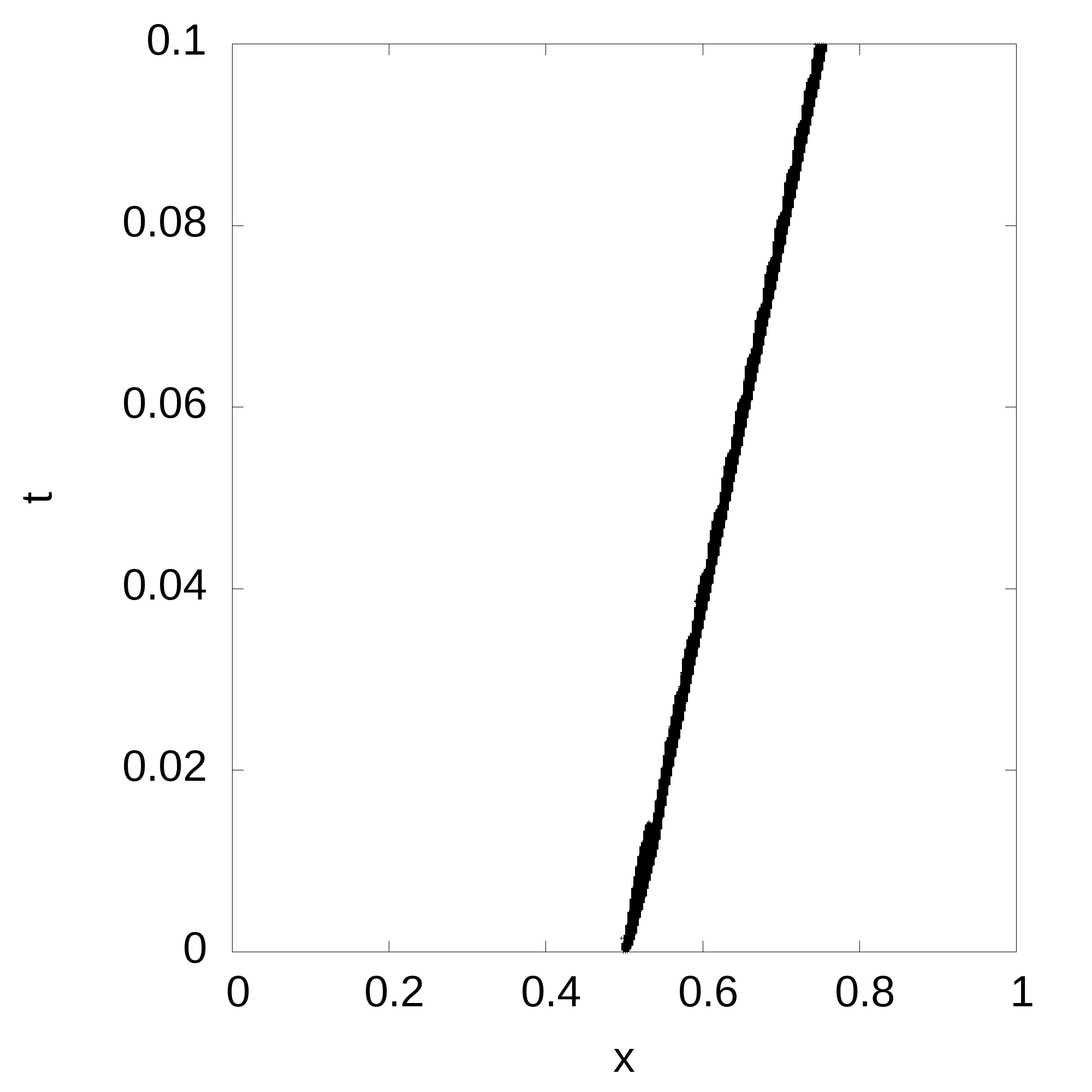}}
  \subfloat[LPR Indicator]{\label{fig:LaxLPR}\includegraphics[width=0.32\textwidth]{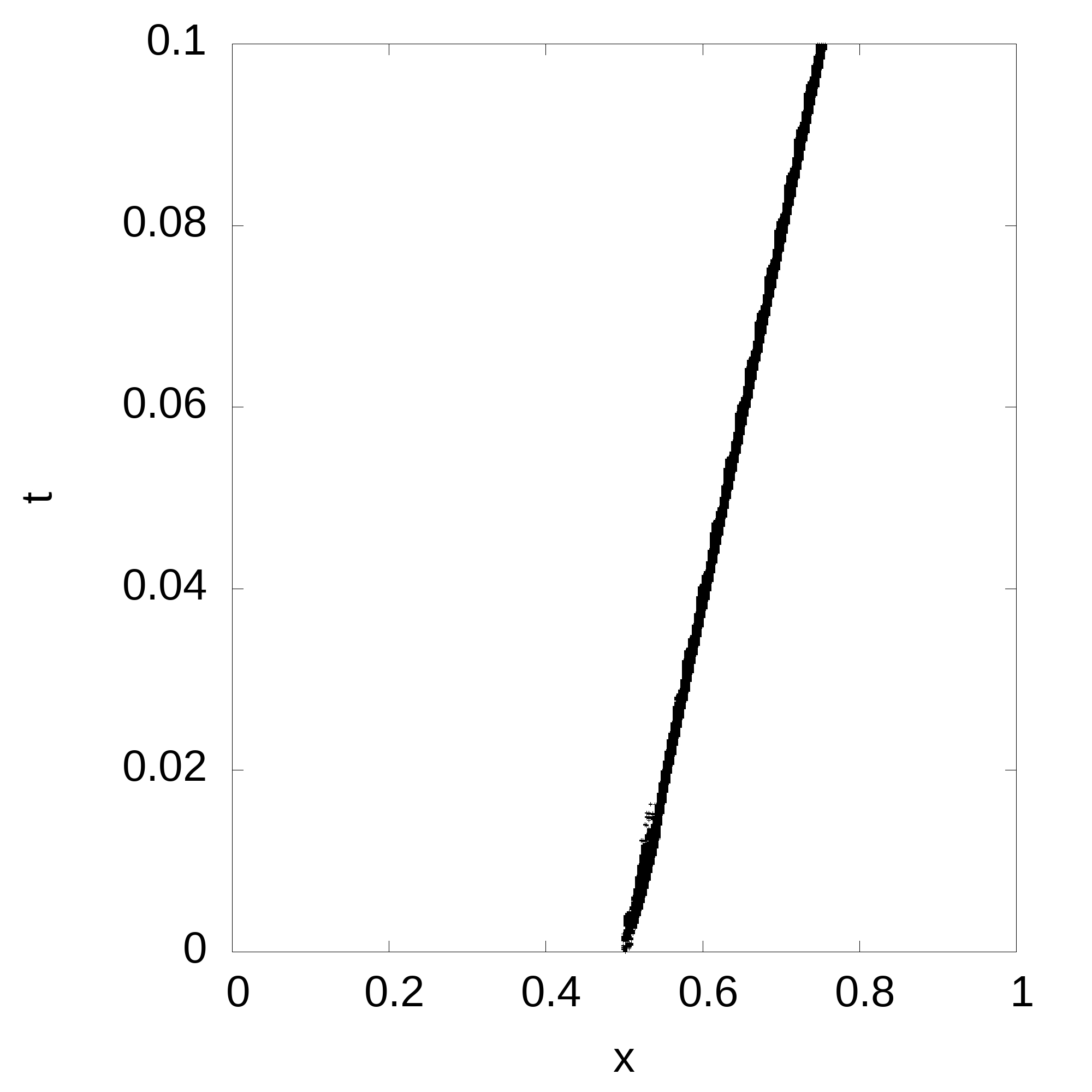}}
  \subfloat[RH Indicator]{\label{fig:LaxRH}\includegraphics[width=0.32\textwidth]{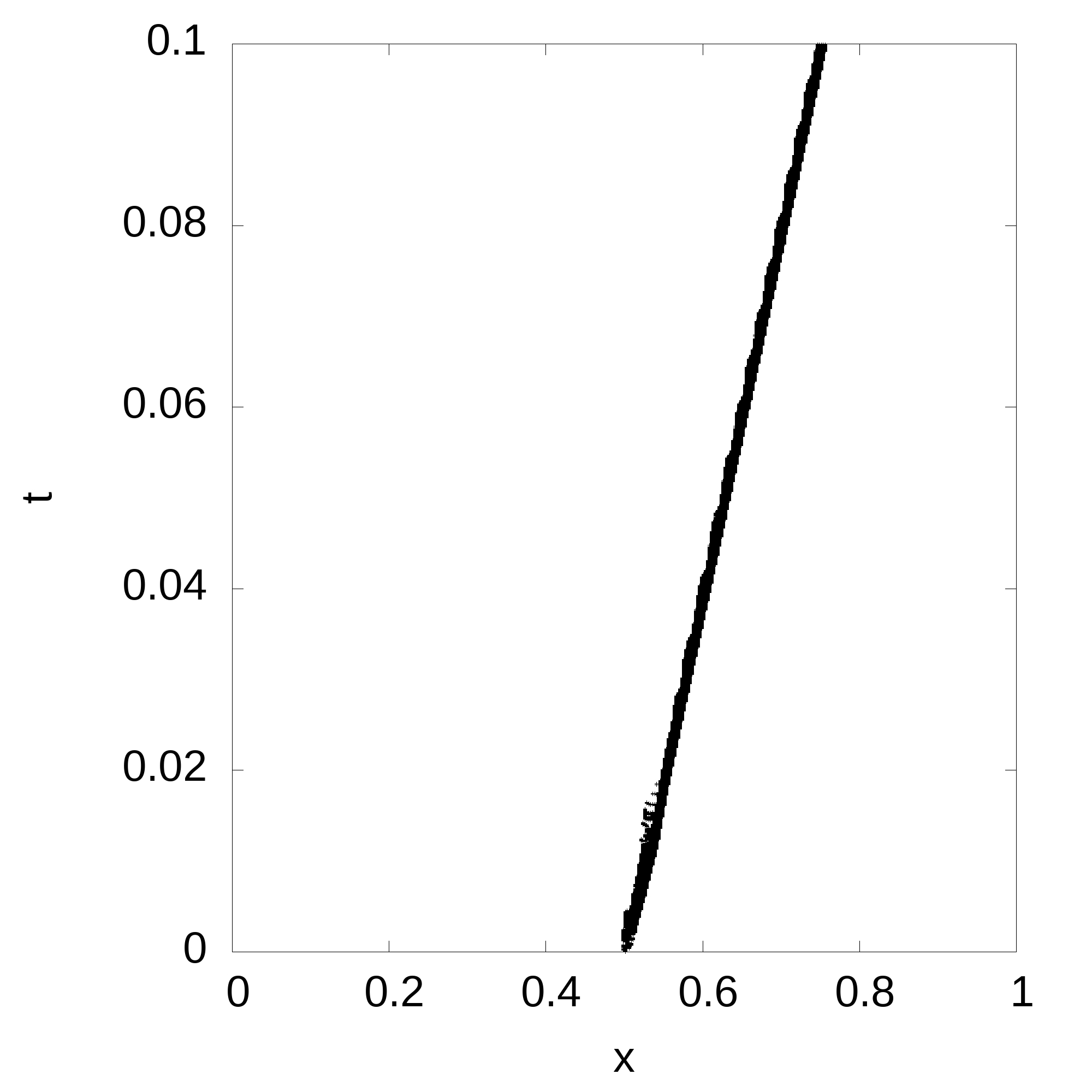}}\hfill
  \subfloat[PPL Indicator]{\label{fig:LaxPPL}\includegraphics[width=0.32\textwidth]{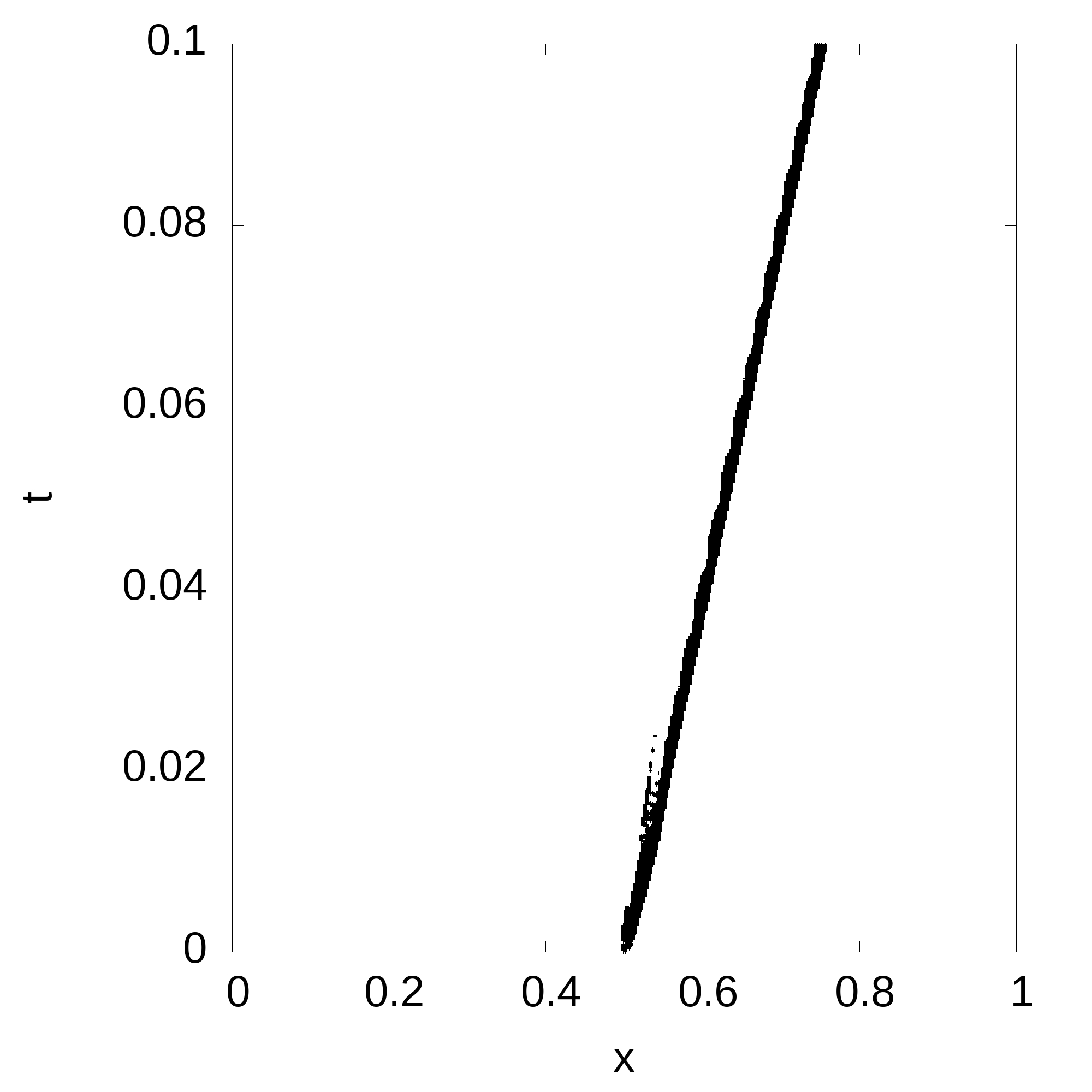}}
  \subfloat[MH Indicator]{\label{fig:LaxMH}\includegraphics[width=0.32\textwidth]{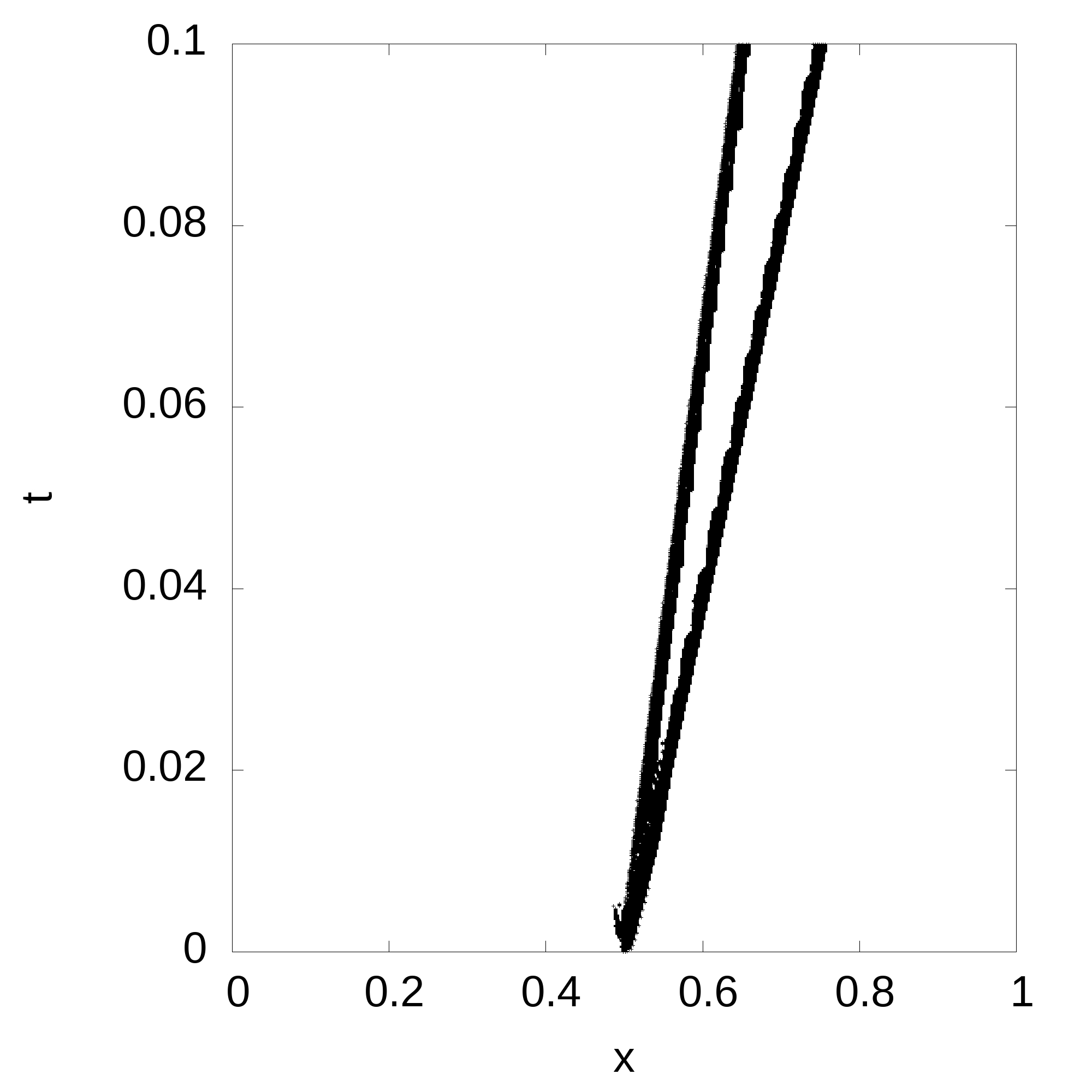}}\hfill
  \caption{The time history of flagged troubled cells of the Lax problem for the one-dimensional Euler equations, simulated until $t=0.1$ with 200 elements and $P^{1}$ based DGM}
  \label{fig:Lax}
\end{figure}

\noindent \textbf{Test Problem 3 (Lax Problem)\cite{lax1}:} We solve the one-dimensional Euler equations as given by \eqref{1dEulerEquations} with the initial condition

\begin{equation}
(\rho, u, p) = 
\begin{cases}
    (0.445,0.698,3.528),& \text{if } x< 0.5\\
    (0.5,0,0.571),              & \text{otherwise}
\end{cases}
\end{equation}

\noindent in the domain $[0,1]$. For all the troubled cell indicators, we again use the density as the detection variable. Average (over all time steps) and maximum percentages of cells being flagged as troubled cells, for the different troubled-cell indicators, are summarized in Table \ref{table:3} for two different grid sizes and various orders. The computed solution for density obtained at $t=0.1$ using 200 elements while using the SJ indicator and CSWENO limiter for $P^{1}$, $P^{2}$ and $P^{3}$ based DGM is compared and plotted against the exact solution in Figure \ref{fig:LaxSolution}. The error ($|\rho-\rho_{exact}|$) obtained for $P^{1}$, $P^{2}$ and $P^{3}$ based DGM is also plotted in Figure \ref{fig:LaxError}. We are not showing the solutions obtained using other indicators as they are quite similar to the solution obtained with the SJ indicator and stay within an $L^{2}$ norm of $\sim  10^{-14}$. We show the $L^{2}$ difference in density between those solutions in table \ref{table:3New}.
\\
We also show the time history of the flagged troubled cells using the eight different indicators in Figure \ref{fig:Lax} using $P^{1}$ based DGM for 200 elements. For all the indicators, again, the number of flagged troubled cells slightly increase when the order increases. From the tabulated results and the figure, the troubled cell indicators PP and PPL perform in a similar fashion, while the SJ indicator is a little better and the MH indicator is quite bad comparatively. We can say that, for this test problem, the troubled cell indicators FS1, FS2, RH, and LPR out perform the other indicators.
\\
\\
\\
\\
\begin{table}[htbp]
\small
\centering
\begin{tabular}{|c|c|c|c|c|c|c|c|c|c|}
%\hline
%\multicolumn{6}{|c|}{$L_{2}$ error for the Riemann problem configuration-10} \\
\hline
 \makecell{No. of \\ Cells} & \makecell{Scheme \\ Indicator} & \multicolumn{2}{|c|}{$P^{1}$} & \multicolumn{2}{|c|}{$P^{2}$} & \multicolumn{2}{|c|}{$P^{3}$} & \multicolumn{2}{|c|}{$P^{4}$} \\
\cline{3-10}
 & & Ave & Max & Ave & Max & Ave & Max & Ave & Max \\
\hline 
\multirow{8}{*}{200} & PP & 5.75 & 7.5 & 6.54 & 7.5 & 7.10 & 8.0 & 7.84 & 8.5 \\
\cline{2-10}
 & SJ & 2.84 & 4.5 & 2.98 & 4.5 & 3.34 & 5.0 & 4.12 & 5.5 \\
\cline{2-10}
 & FS1 & 2.32 & 4.0 & 2.43 & 4.5 & 2.55 & 4.5 & 2.74 & 5.0 \\
\cline{2-10}
 & FS2 & 2.35 & 4.0 & 2.44 & 4.5 & 2.55 & 4.5 & 2.75 & 5.0 \\
\cline{2-10}
 & LPR & 2.24 & 3.5 & 2.35 & 3.5 & 2.47 & 4.0 & 2.69 & 4.5\\
\cline{2-10}
 & RH & 2.45 & 4.5 & 2.74 & 4.5 & 3.02 & 5.0 & 3.45 & 5.5 \\
 \cline{2-10}
 & PPL & 4.86 & 7.0 & 5.32 & 7.5 & 6.11 & 7.5 & 7.13 & 8.5 \\
\cline{2-10}
 & MH & 9.55 & 13.5 & 10.22 & 14.5 & 10.59 & 14.5 & 11.02 & 15.0\\
\hline
\multirow{8}{*}{400} & PP & 4.94 & 7.0 & 5.73 & 7.5 & 6.33 & 7.75 & 7.24 & 8.5 \\
\cline{2-10}
 & SJ & 2.64 & 4.5 & 2.85 & 4.5 & 3.25 & 4.75 & 4.02 & 5.5 \\
\cline{2-10}
 & FS1 & 2.24 & 3.5 & 2.33 & 3.5 & 2.47 & 3.75 & 2.63 & 4.0 \\
\cline{2-10}
 & FS2 & 2.24 & 3.5 & 2.33 & 3.5 & 2.47 & 3.75 & 2.63 & 4.0 \\
\cline{2-10}
 & LPR & 2.20 & 3.25 & 2.31 & 3.25 & 2.44 & 3.75 & 2.64 & 4.25\\
\cline{2-10}
 & RH & 2.34 & 4.0 & 2.64 & 4.5 & 2.91 & 5.0 & 3.37 & 5.25 \\
 \cline{2-10}
 & PPL & 4.14 & 6.75 & 4.64 & 7.0 & 5.14 & 7.25 & 6.73 & 8.0 \\
\cline{2-10}
 & MH & 9.23 & 13.25 & 9.98 & 14.0 & 10.16 & 14.25 & 10.86 & 14.75\\
\hline
\end{tabular}
\caption{Average (marked as Ave) and maximum (marked as Max) percentages of cells flagged as troubled cells subject to different troubled-cell indicators for the Lax problem for various orders and two grid sizes.}
\label{table:3}
\end{table}

\begin{table}[htbp]
%\small
\centering
\begin{tabular}{|c|c|c|c|c|}
%\hline
%\multicolumn{6}{|c|}{$L_{2}$ error for the Riemann problem configuration-10} \\
\hline
Scheme Indicator & \multicolumn{4}{|c|}{\makecell{$L^{2}$ difference in density in comparison with \\ SJ indicator solution for various orders}} \\ \cline{2-5}
  & $P^{1}$ & $P^{2}$ & $P^{3}$ & $P^{4}$ \\ 
\hline
PP & 1.7E-14 & 8.2E-15 & 7.8E-15 & 7.1E-15\\
\hline
FS1 & 4.3E-14 & 7.3E-15 & 2.8E-15 & 4.4E-15\\
\hline
FS2 & 4.3E-14 & 7.3E-15 & 2.8E-15 & 4.4E-15\\
\hline
LPR & 6.2E-15 & 3.2E-15 & 4.1E-15 & 3.3E-15\\
\hline
RH & 2.9E-14 & 9.5E-15 & 3.9E-15 & 5.2E-15\\
\hline
PPL & 5.7E-15 & 5.5E-15 & 6.7E-15 & 7.4E-15\\
\hline
MH & 8.9E-15 & 4.9E-15 & 2.7E-15 & 3.1E-15\\
\hline
\end{tabular}
\caption{$L^{2}$ difference in density between solution obtained using CSWENO limiter and SJ indicator (shown in Figure \ref{fig:LaxCSWENO}) and the solution obtained using other indicators for the Lax problem using 200 elements }
\label{table:3New}
\end{table}

\begin{figure}[htbp]
  \centering
  \subfloat[Solution of 1D Euler equations for the Shu-Osher problem with $P^{1}$, $P^{2}$ and $P^{3}$ based DGM]{\label{fig:ShockDensityWaveSolution}\includegraphics[width=0.48\textwidth]{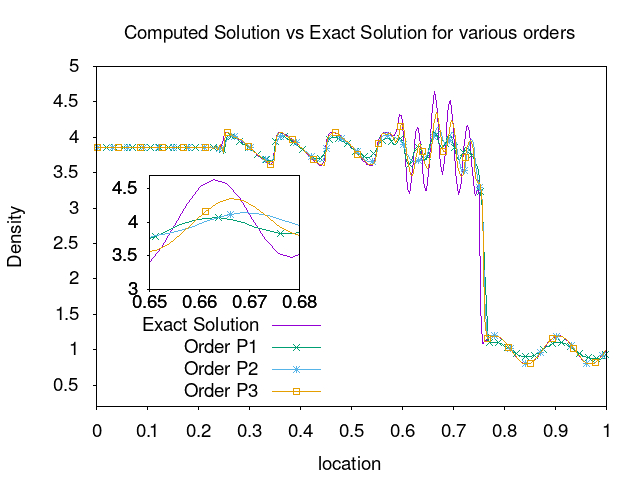}}
  \subfloat[Density error ($|\rho-\rho_{exact}|$) for $P^{1}$, $P^{2}$ and $P^{3}$ based DGM]{\label{fig:ShockDensityWaveError}\includegraphics[width=0.48\textwidth]{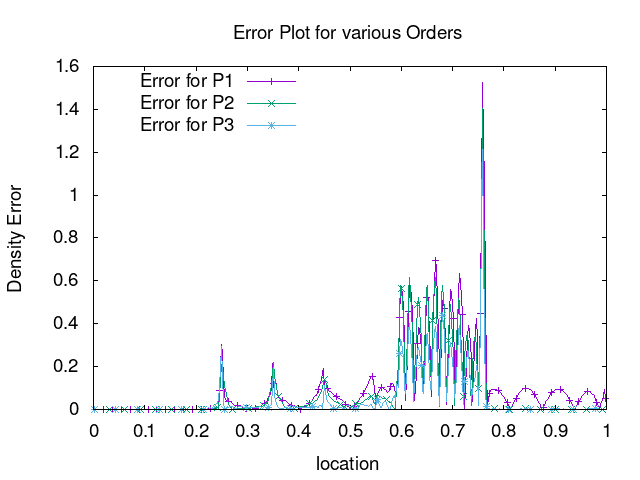}}\hfill
  \caption{Comparison of density solutions of Shu-Osher problem at $t=0.178$ using 200 elements obtained with the $P^{1}$, $P^{2}$ and $P^{3}$ based DGM and the exact solution. Density error ($|\rho-\rho_{exact}|$) for $P^{1}$, $P^{2}$ and $P^{3}$ based DGM is also plotted. Figure \ref{fig:ShockDensityWaveSolution} also includes a zoomed in portion of the solution for better comparison}
  \label{fig:ShockDensityWaveCSWENO}
\end{figure}

\begin{figure}[htbp]
  \centering
  \subfloat[PP Indicator]{\label{fig:SOPP}\includegraphics[width=0.32\textwidth]{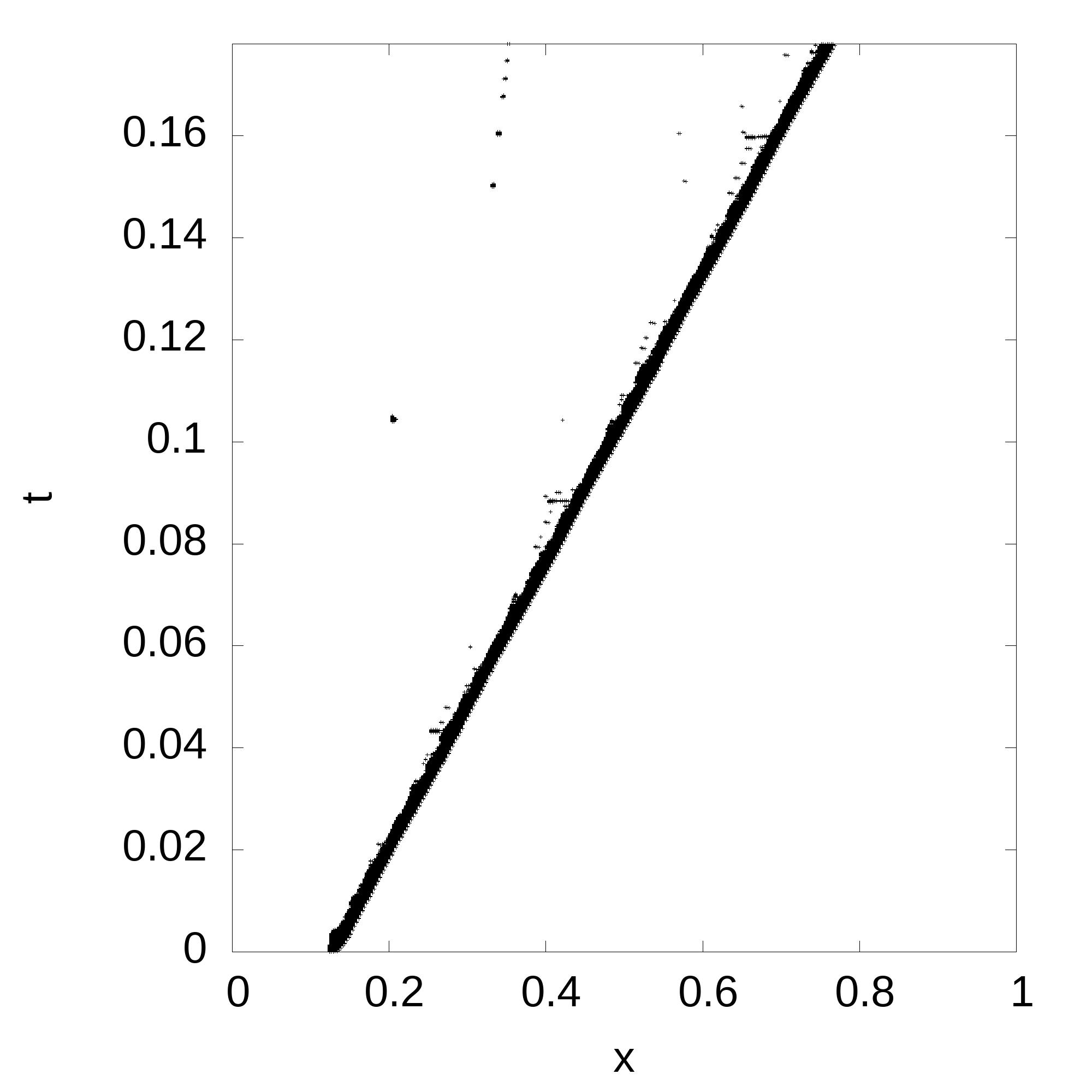}}
  \subfloat[SJ Indicator]{\label{fig:SOSJ}\includegraphics[width=0.32\textwidth]{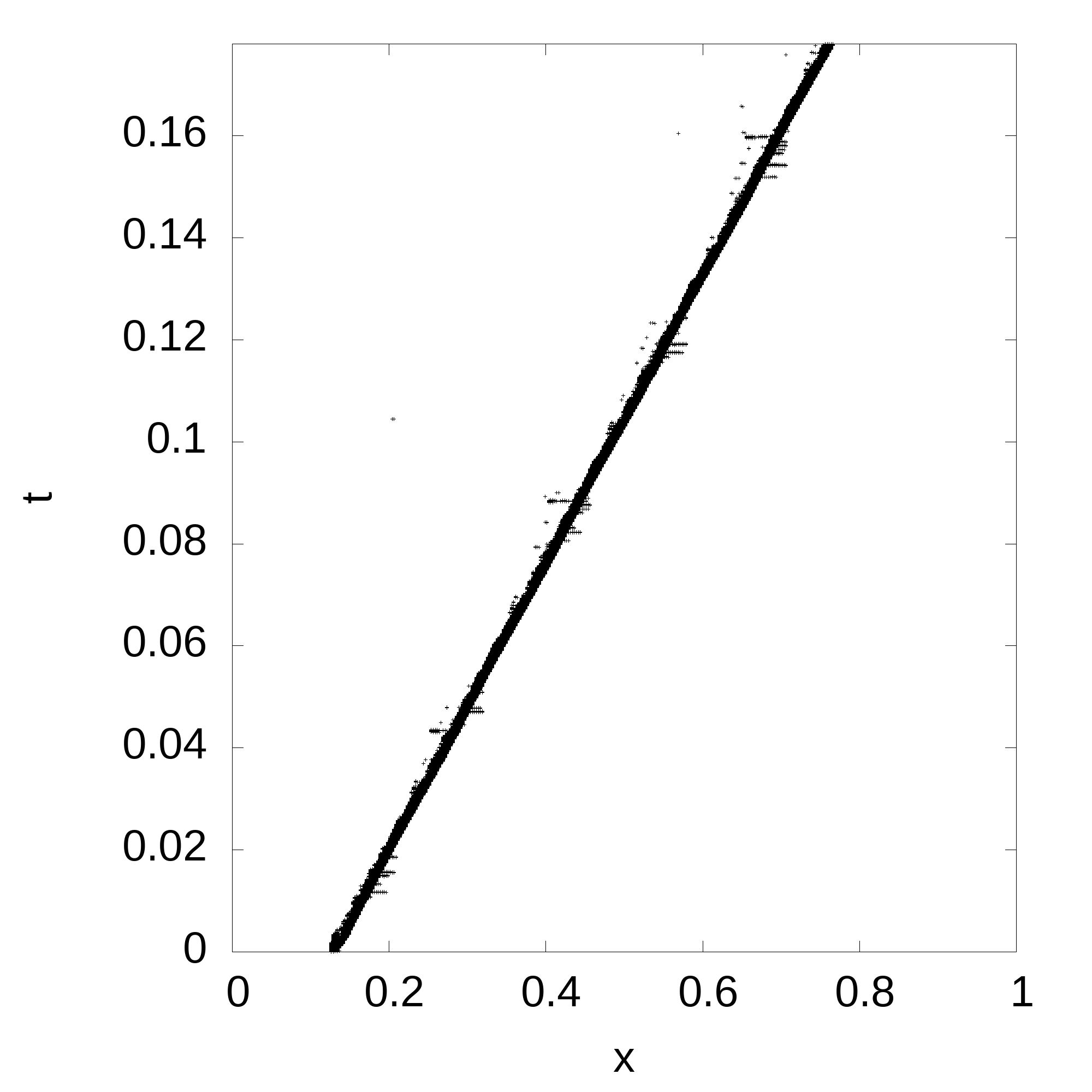}}
  \subfloat[FS1 Indicator]{\label{fig:SOFS1}\includegraphics[width=0.32\textwidth]{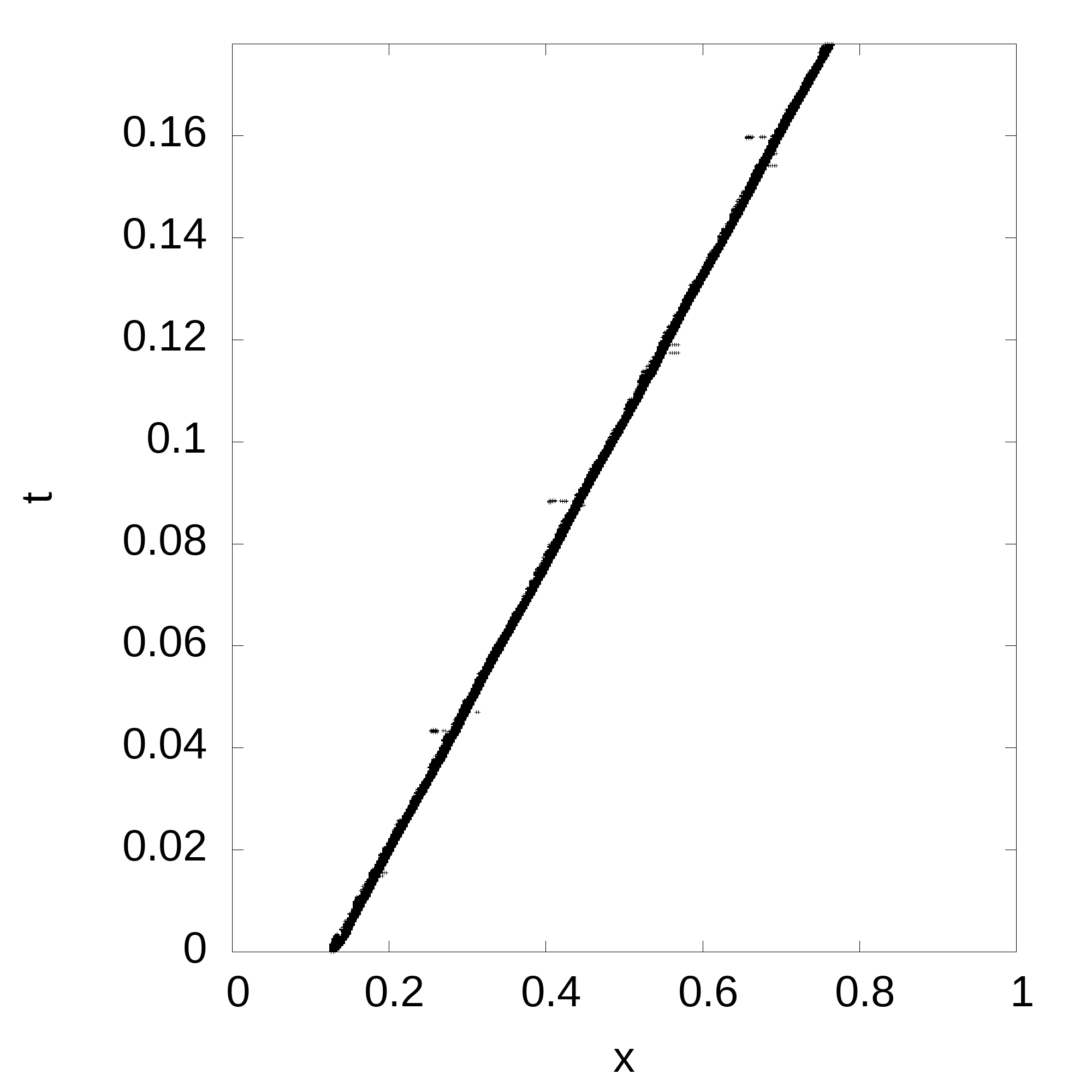}}\hfill
  \subfloat[FS2 Indicator]{\label{fig:SOFS2}\includegraphics[width=0.32\textwidth]{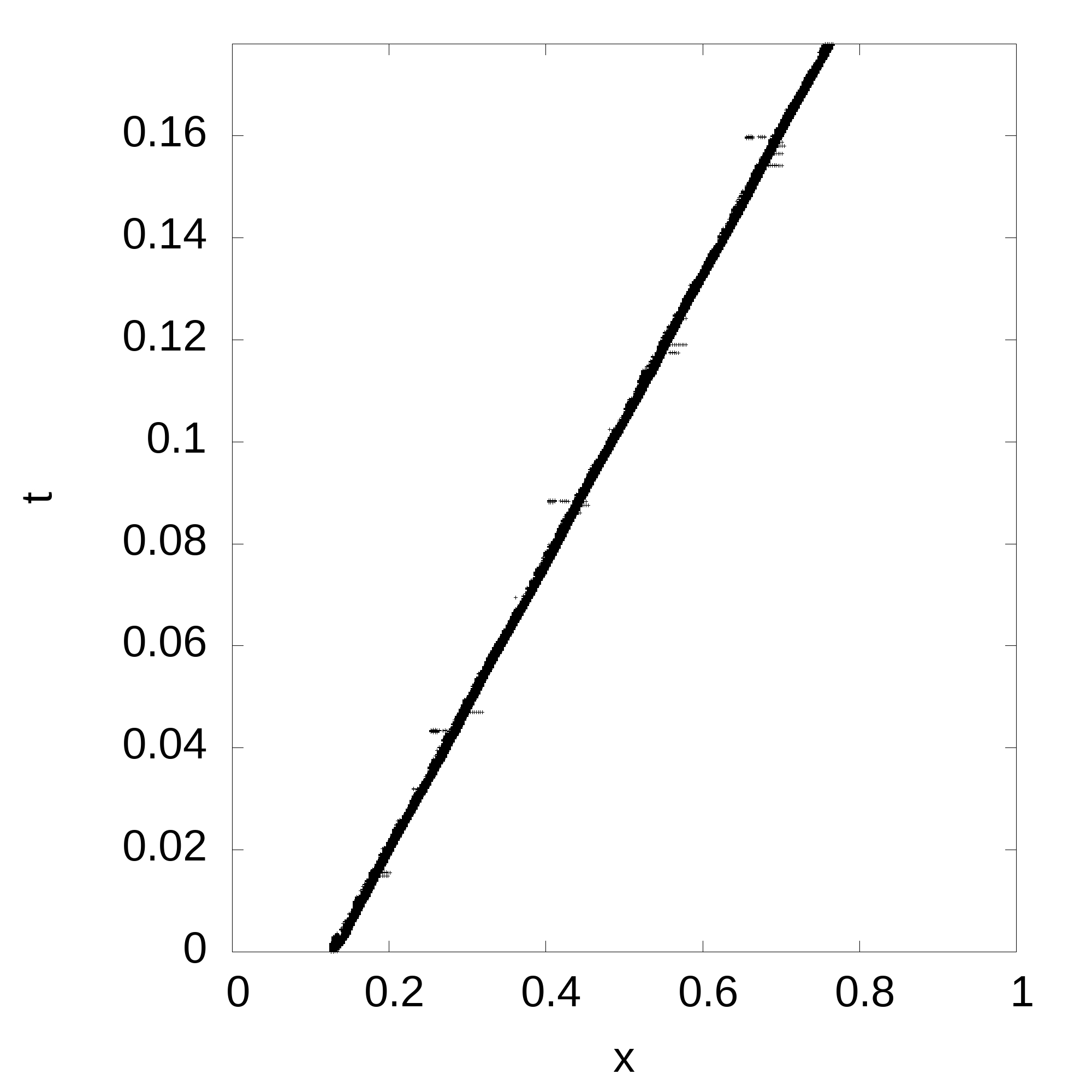}}
  \subfloat[LPR Indicator]{\label{fig:SOLPR}\includegraphics[width=0.32\textwidth]{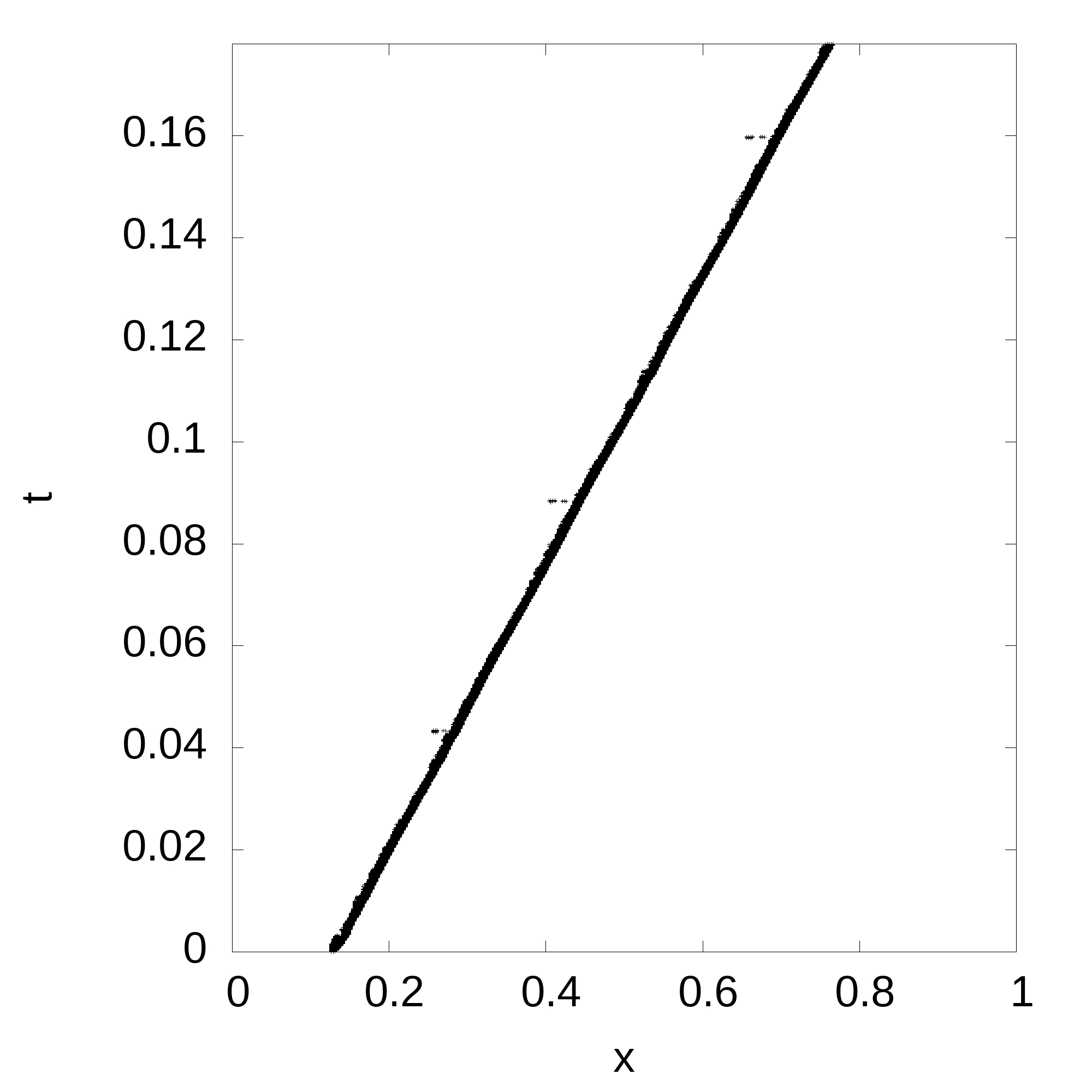}}
  \subfloat[RH Indicator]{\label{fig:SORH}\includegraphics[width=0.32\textwidth]{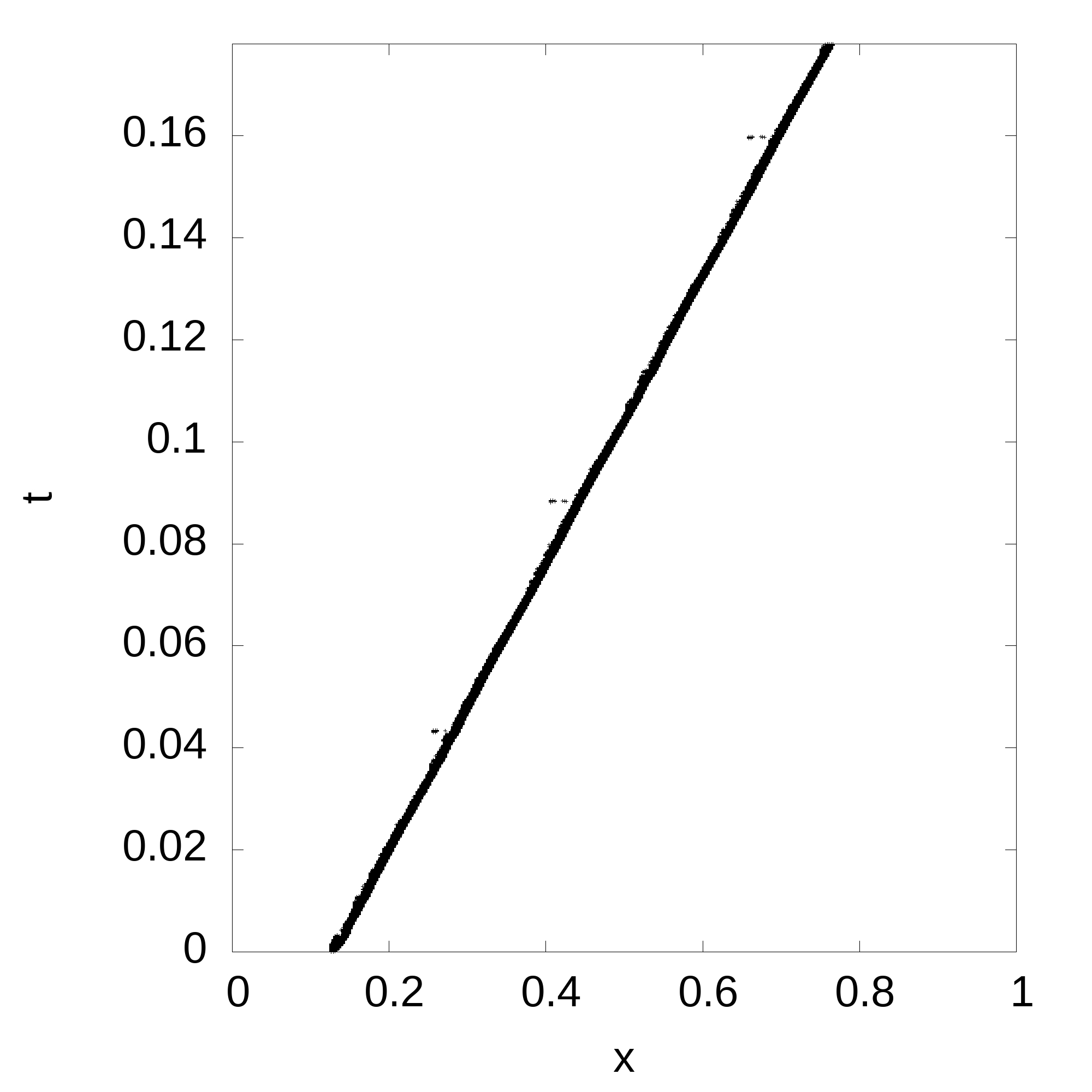}}\hfill
  \subfloat[PPL Indicator]{\label{fig:SOPPL}\includegraphics[width=0.32\textwidth]{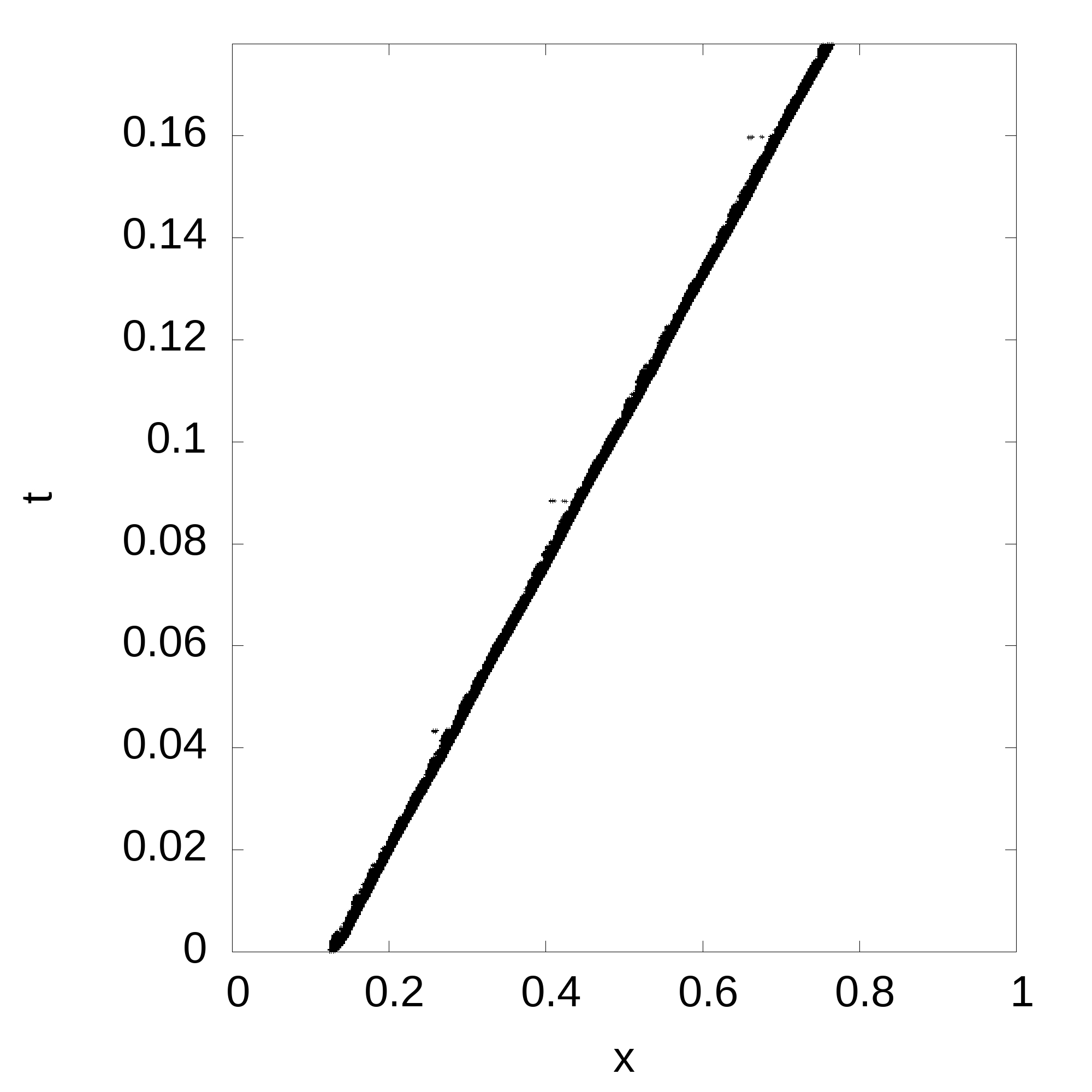}}
  \subfloat[MH Indicator]{\label{fig:SOMH}\includegraphics[width=0.32\textwidth]{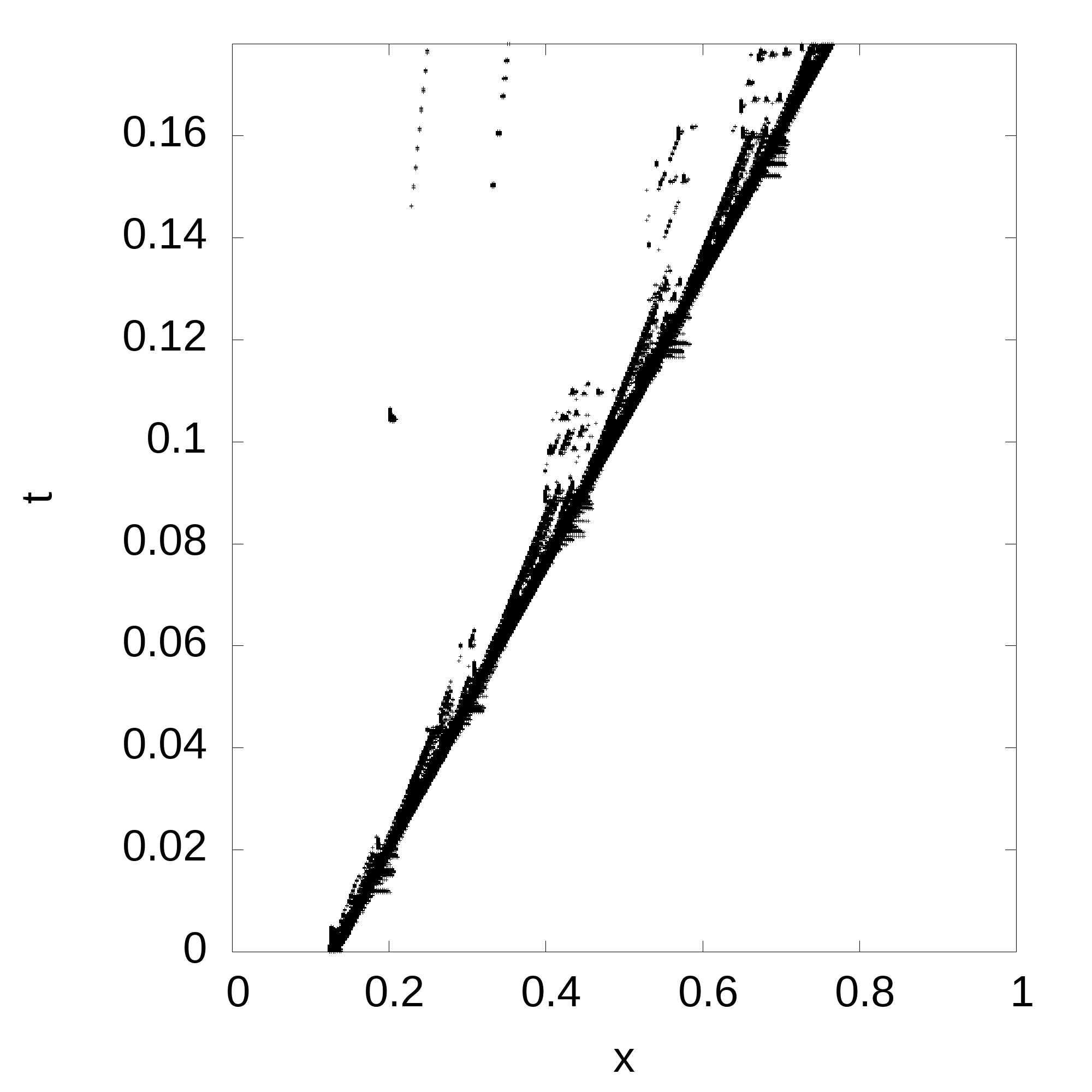}}\hfill
  \caption{The time history of flagged troubled cells of the Shu-Osher problem for the one-dimensional Euler equations, simulated until $t=0.178$ with 200 elements and $P^{1}$ based DGM}
  \label{fig:SO}
\end{figure}

\noindent \textbf{Test Problem 4 (Shu-Osher Problem)\cite{so2}:} We solve the one-dimensional Euler equations as given by \eqref{1dEulerEquations} with the initial condition

\begin{equation}
(\rho, u, p) = 
\begin{cases}
    (3.857143,2.629369,10.333333),& \text{if } x< 0.125\\
    (1.0+0.2\sin(16\pi x),0,1),              & \text{otherwise}
\end{cases}
\end{equation}

\noindent in the domain $[0,1]$. For all the troubled cell indicators, we again use the density as the detection variable. Average (over all time steps) and maximum percentages of cells being flagged as troubled cells, for the different troubled-cell indicators, are summarized in Table \ref{table:4} for two different grid sizes and various orders. The computed solution for density obtained at $t=0.178$s using 200 elements while using the SJ indicator and CSWENO limiter for $P^{1}$, $P^{2}$ and $P^{3}$ based DGM is compared and plotted against the exact solution in Figure \ref{fig:ShockDensityWaveSolution}. The error ($|\rho-\rho_{exact}|$) obtained for $P^{1}$, $P^{2}$ and $P^{3}$ based DGM is also plotted in Figure \ref{fig:ShockDensityWaveError}. We are not showing the solutions obtained using other indicators as they are quite similar to the solution obtained with the SJ indicator and stay within an $L^{2}$ norm of $\sim  10^{-14}$. We show the $L^{2}$ difference in density between those solutions in table \ref{table:4New}.
\\ 
We also show the time history of the flagged troubled cells using the eight different indicators in Figure \ref{fig:SO} using $P^{1}$ based DGM for 200 elements. For all the indicators, again, the number of flagged troubled cells slightly increase when the order increases. From the tabulated results and the figure, the troubled cell indicators PP and PPL perform in a similar fashion, while the SJ indicator is quite better and the MH indicator is quite bad comparatively. We can also say that the troubled cell indicators FS1, FS2, RH, and LPR out perform the other four indicators and perform in a similar fashion.
\\
\\
\begin{table}[htbp]
\small
\centering
\begin{tabular}{|c|c|c|c|c|c|c|c|c|c|}
%\hline
%\multicolumn{6}{|c|}{$L_{2}$ error for the Riemann problem configuration-10} \\
\hline
 \makecell{No. of \\ Cells} & \makecell{Scheme \\ Indicator} & \multicolumn{2}{|c|}{$P^{1}$} & \multicolumn{2}{|c|}{$P^{2}$} & \multicolumn{2}{|c|}{$P^{3}$} & \multicolumn{2}{|c|}{$P^{4}$} \\
\cline{3-10}
 & & Ave & Max & Ave & Max & Ave & Max & Ave & Max \\
\hline 
\multirow{8}{*}{200} & PP & 3.84 & 7.5 & 4.22 & 7.5 & 4.35 & 8.0 & 4.97 & 8.0 \\
\cline{2-10}
 & SJ & 2.96 & 4.5 & 3.24 & 5.0 & 3.84 & 5.0 & 4.17 & 5.5 \\
\cline{2-10}
 & FS1 & 2.35 & 3.5 & 2.83 & 3.5 & 3.08 & 4.0 & 3.66 & 4.0 \\
\cline{2-10}
 & FS2 & 2.42 & 3.0 & 2.87 & 3.5 & 3.07 & 4.0 & 3.69 & 4.0 \\
\cline{2-10}
 & LPR & 2.64 & 4.5 & 3.14 & 4.5 & 3.87 & 5.0 & 4.19 & 5.5\\
\cline{2-10}
 & RH & 2.53 & 3.5 & 2.99 & 4.0 & 3.17 & 4.5 & 3.62 & 5.0 \\
 \cline{2-10}
 & PPL & 3.46 & 6.5 & 4.10 & 7.5 & 4.23 & 7.5 & 4.76 & 8.0 \\
\cline{2-10}
 & MH & 7.24 & 13.5 & 7.79 & 14.0 & 8.17 & 14.5 & 9.11 & 14.5\\
\hline
\multirow{8}{*}{400} & PP & 3.64 & 7.25 & 4.05 & 7.5 & 4.27 & 7.75 & 4.78 & 8.0 \\
\cline{2-10}
 & SJ & 2.84 & 4.25 & 3.13 & 4.5 & 3.65 & 5.0 & 3.98 & 5.25 \\
\cline{2-10}
 & FS1 & 2.15 & 3.25 & 2.64 & 3.5 & 3.01 & 3.75 & 3.47 & 4.0 \\
\cline{2-10}
 & FS2 & 2.22 & 3.0 & 2.68 & 3.5 & 3.09 & 4.0 & 3.54 & 4.25 \\
\cline{2-10}
 & LPR & 2.44 & 4.25 & 3.10 & 4.75 & 3.69 & 5.0 & 3.96 & 5.25\\
\cline{2-10}
 & RH & 2.39 & 3.25 & 2.85 & 4.25 & 3.14 & 4.75 & 3.48 & 5.0 \\
 \cline{2-10}
 & PPL & 3.25 & 6.25 & 3.94 & 7.0 & 4.05 & 7.25 & 4.65 & 7.75 \\
\cline{2-10}
 & MH & 7.03 & 13.25 & 7.54 & 14.0 & 7.98 & 14.25 & 8.84 & 14.5\\
\hline
\end{tabular}
\caption{Average (marked as Ave) and maximum (marked as Max) percentages of cells flagged as troubled cells subject to different troubled-cell indicators for the Shu-Osher problem for various orders.}
\label{table:4}
\end{table}

\begin{table}[htbp]
%\small
\centering
\begin{tabular}{|c|c|c|c|c|}
%\hline
%\multicolumn{6}{|c|}{$L_{2}$ error for the Riemann problem configuration-10} \\
\hline
Scheme Indicator & \multicolumn{4}{|c|}{\makecell{$L^{2}$ difference in density in comparison with \\ SJ indicator solution for various orders}} \\ \cline{2-5}
  & $P^{1}$ & $P^{2}$ & $P^{3}$ & $P^{4}$ \\ 
\hline
PP & 3.1E-14 & 1.8E-14 & 4.9E-15 & 1.9E-15\\
\hline
FS1 & 2.2E-14 & 9.9E-14 & 5.3E-15 & 8.2E-15\\
\hline
FS2 & 2.2E-14 & 9.9E-14 & 5.3E-15 & 8.2E-15\\
\hline
LPR & 1.9E-14 & 3.5E-15 & 1.5E-15 & 3.5E-15\\
\hline
RH & 8.7E-14 & 2.8E-15 & 3.5E-15 & 9.1E-15\\
\hline
PPL & 6.2E-14 & 8.9E-15 & 7.3E-15 & 3.6E-15\\
\hline
MH & 7.6E-14 & 6.6E-15 & 6.7E-15 & 4.2E-15\\
\hline
\end{tabular}
\caption{$L^{2}$ difference in density between solution obtained using CSWENO limiter and SJ indicator (shown in Figure \ref{fig:ShockDensityWaveCSWENO}) and the solution obtained using other indicators for the Shu-Osher problem using 200 elements}
\label{table:4New}
\end{table}

\begin{figure}[htbp]
  \centering
  \subfloat[Solution of 1D Euler equations for the shock entropy wave problem with $P^{1}$, $P^{2}$ and $P^{3}$ based DGM]{\label{fig:BlastWaveSolution}\includegraphics[width=0.48\textwidth]{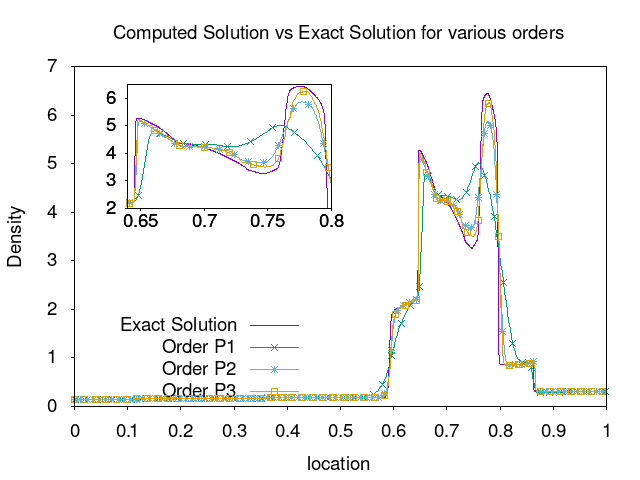}}
  \subfloat[Density error ($|\rho-\rho_{exact}|$) for $P^{1}$, $P^{2}$ and $P^{3}$ based DGM]{\label{fig:BlastWaveError}\includegraphics[width=0.48\textwidth]{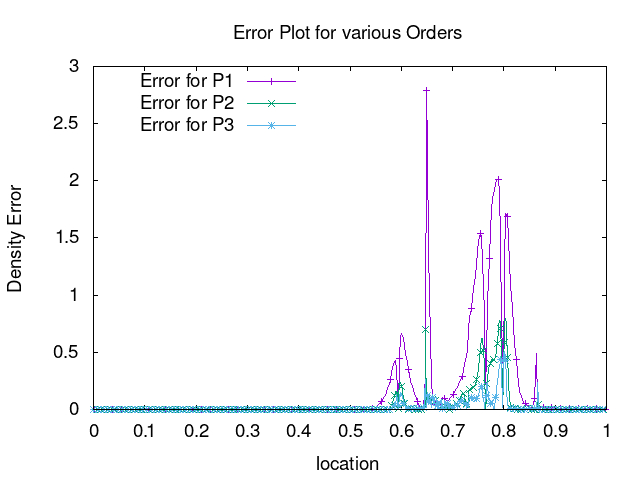}}\hfill
  \caption{Comparison of density solutions of blast wave Problem at $t=0.038$ using 400 elements obtained with the $P^{1}$, $P^{2}$ and $P^{3}$ based DGM and the exact solution. Density error ($|\rho-\rho_{exact}|$) for $P^{1}$, $P^{2}$ and $P^{3}$ based DGM is also plotted. Figure \ref{fig:BlastWaveSolution} also includes a zoomed in portion of the solution for better comparison}
  \label{fig:BlastWaveCSWENO}
\end{figure}

\begin{figure}[htbp]
  \centering
  \subfloat[PP Indicator]{\label{fig:BWPP}\includegraphics[width=0.32\textwidth]{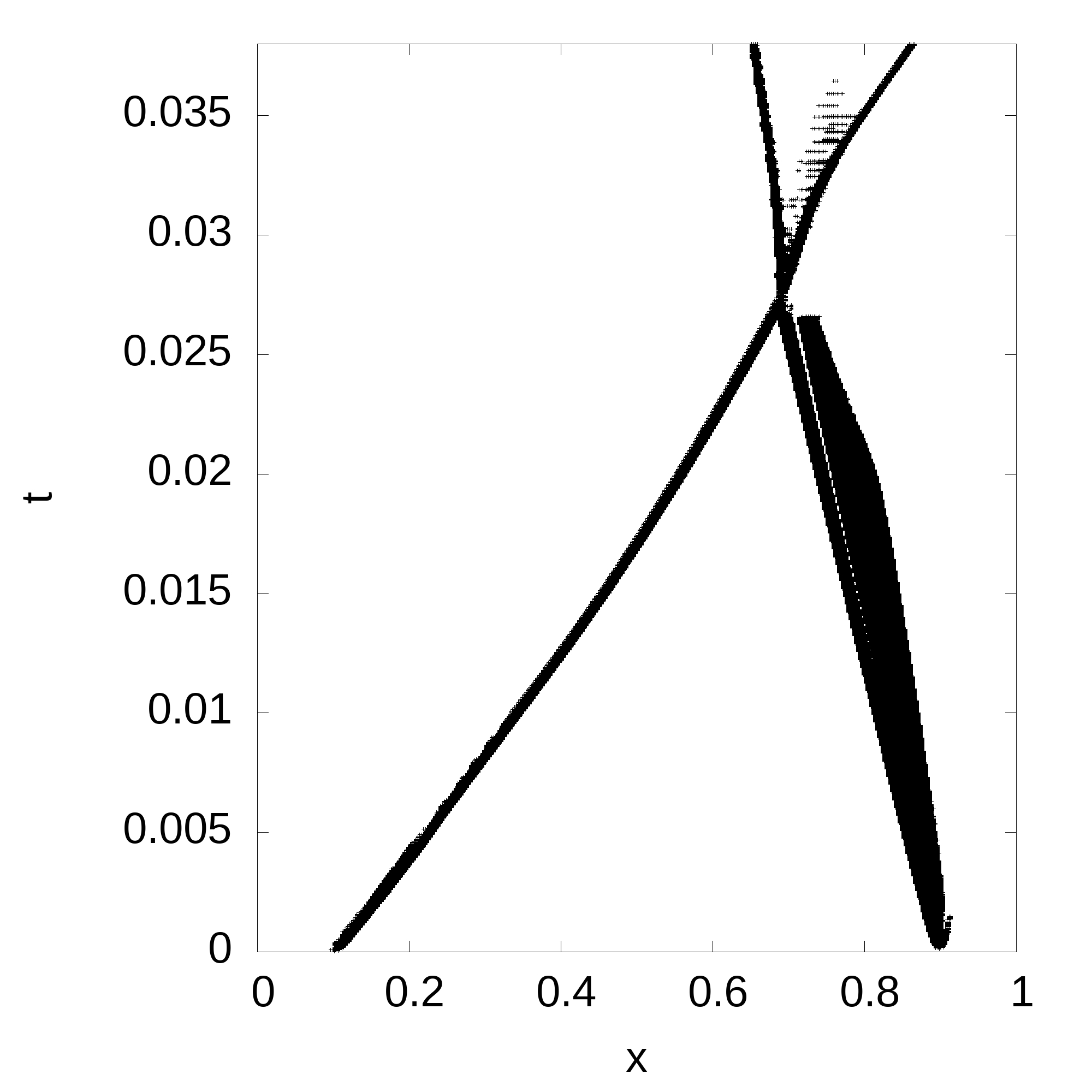}}
  \subfloat[SJ Indicator]{\label{fig:BWSJ}\includegraphics[width=0.32\textwidth]{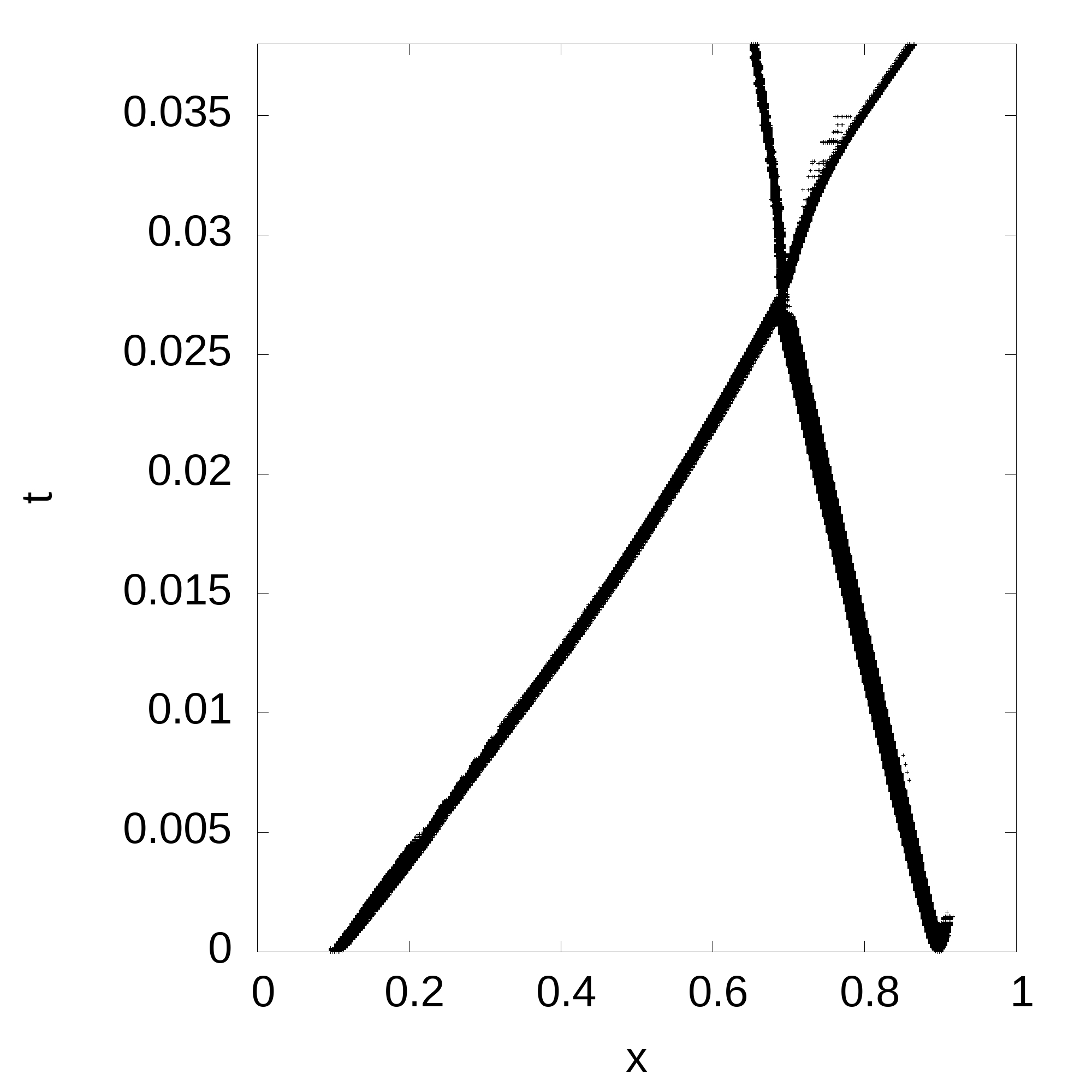}}
  \subfloat[FS1 Indicator]{\label{fig:BWFS1}\includegraphics[width=0.32\textwidth]{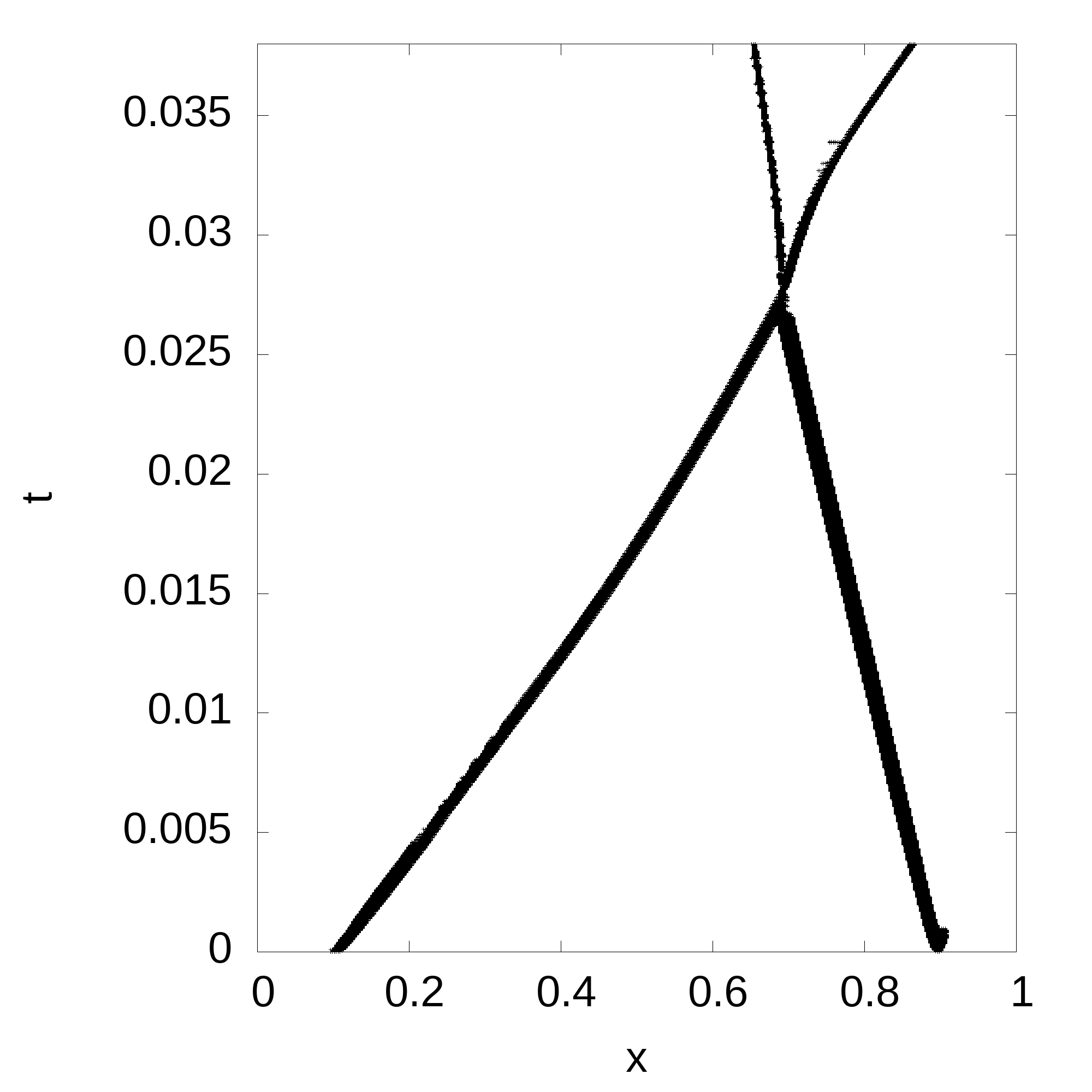}}\hfill
  \subfloat[FS2 Indicator]{\label{fig:BWFS2}\includegraphics[width=0.32\textwidth]{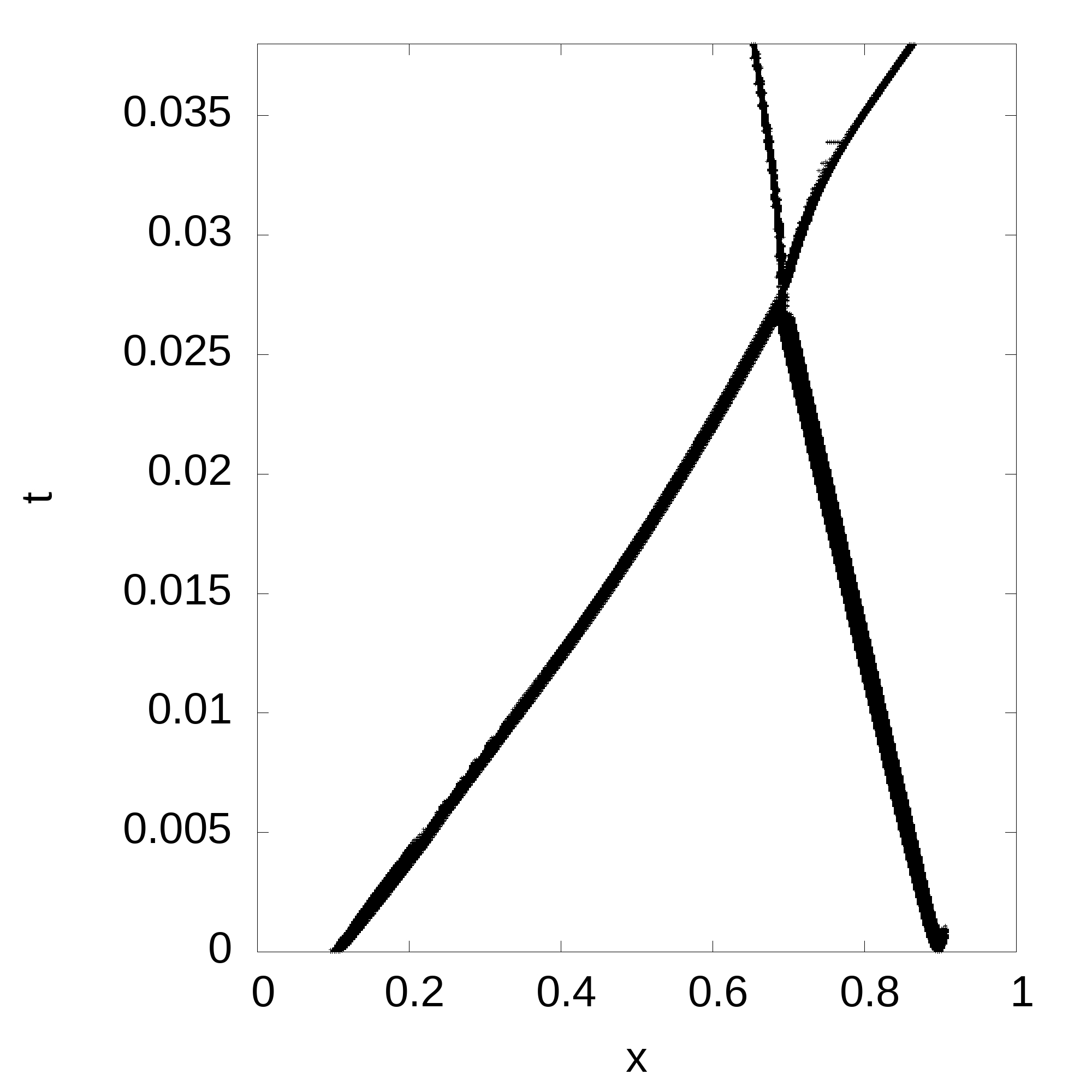}}
  \subfloat[LPR Indicator]{\label{fig:BWLPR}\includegraphics[width=0.32\textwidth]{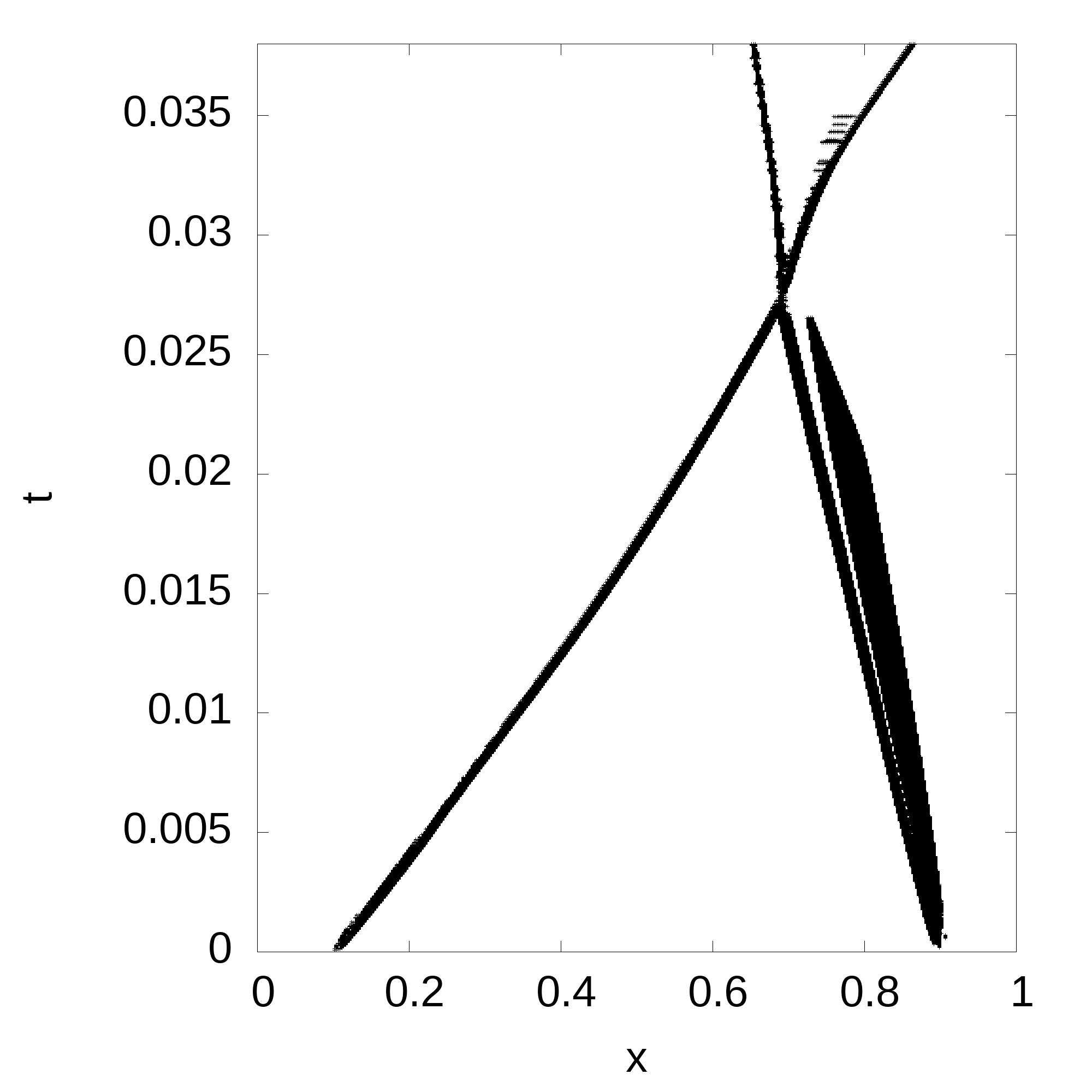}}
  \subfloat[RH Indicator]{\label{fig:BWRH}\includegraphics[width=0.32\textwidth]{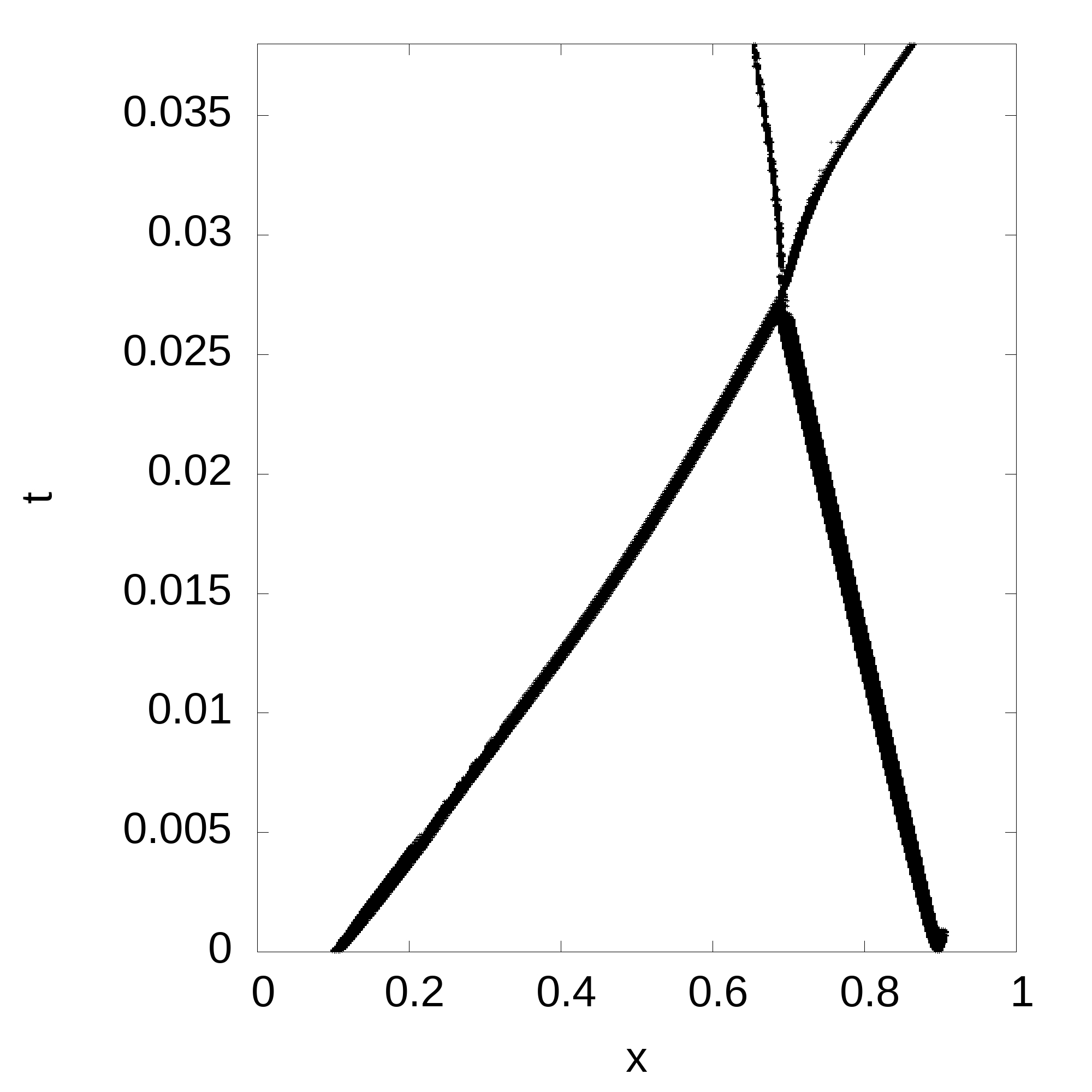}}\hfill
  \subfloat[PPL Indicator]{\label{fig:BWPPL}\includegraphics[width=0.32\textwidth]{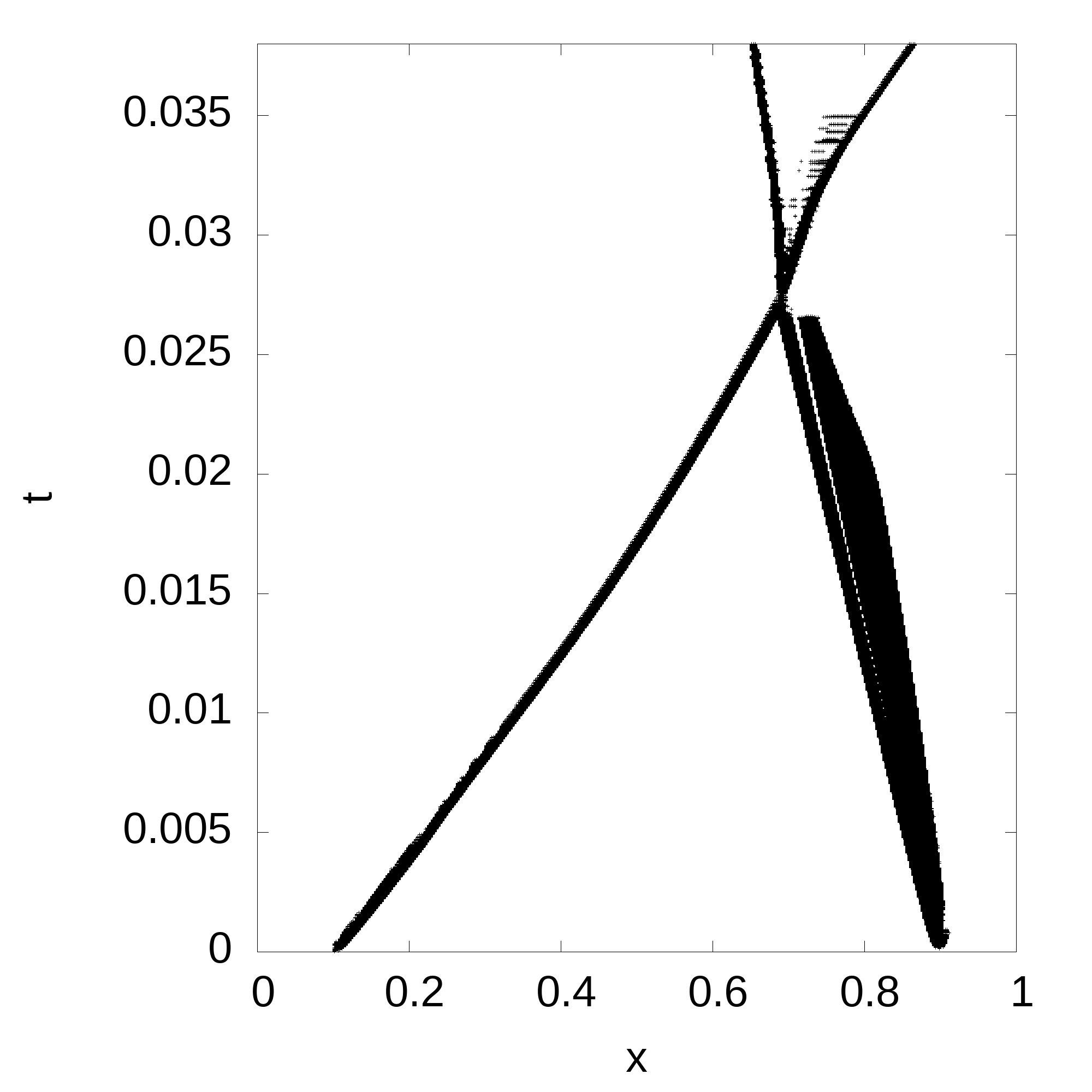}}
  \subfloat[MH Indicator]{\label{fig:BWMH}\includegraphics[width=0.32\textwidth]{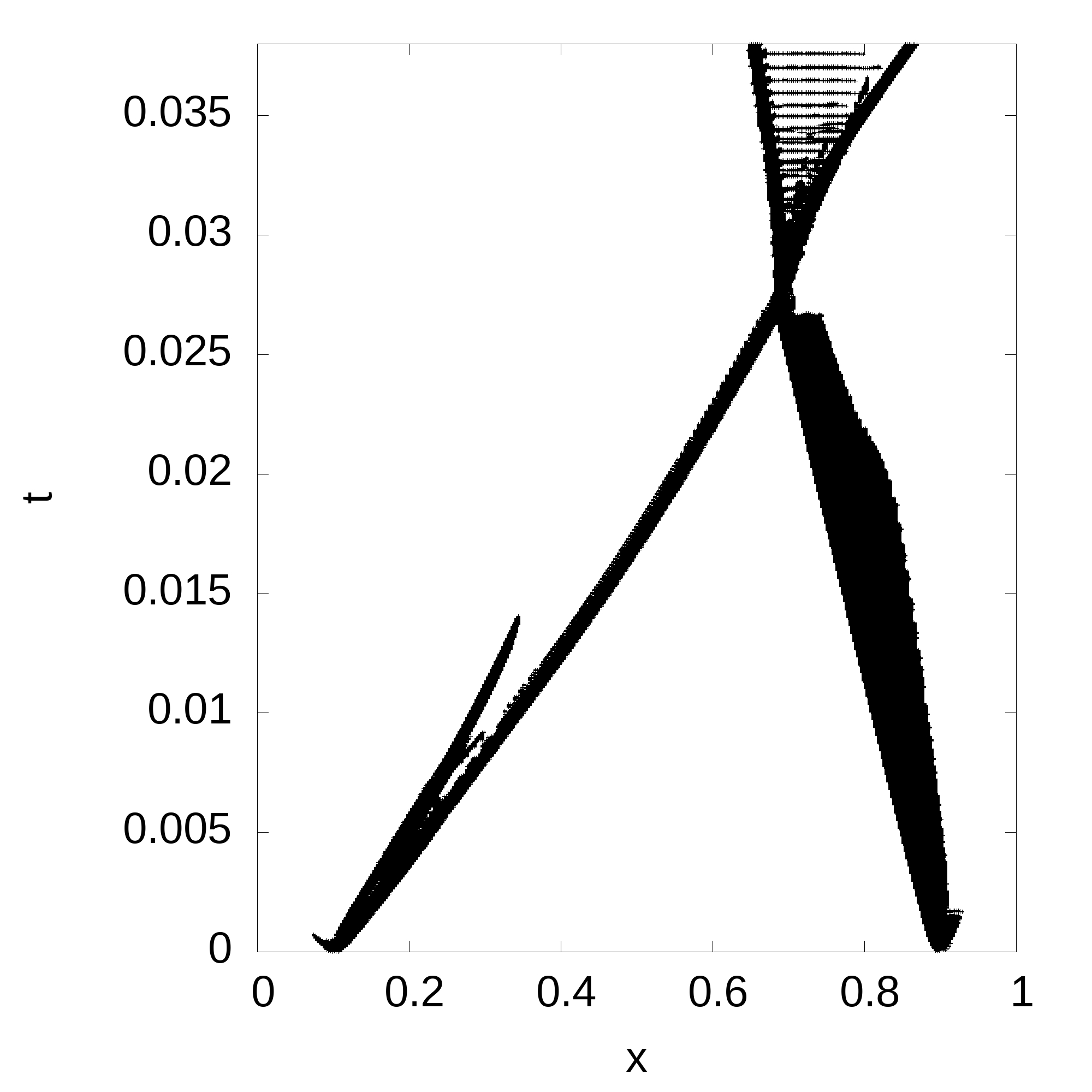}}\hfill
  \caption{The time history of flagged troubled cells of the blast wave problem for the one-dimensional Euler equations, simulated until $t=0.038$ with 200 elements and $P^{1}$ based DGM}
  \label{fig:BW}
\end{figure}

\noindent \textbf{Test Problem 5 (Blast Wave Problem)\cite{wc}:} We solve the one-dimensional Euler equations as given by \eqref{1dEulerEquations} with the initial condition

\begin{equation}
(\rho, u, p) = 
\begin{cases}
    (1,0,1000),& \text{if } x< 0.1\\
    (1,0,0.01),& \text{if } 0.1\leq x< 0.9\\
    (1,0,100),              & \text{otherwise}
\end{cases}
\end{equation}

\noindent in the domain $[0,1]$ with reflecting boundary conditions on both sides. This is a more demanding test problem. For all the troubled cell indicators, we again use the density as the detection variable. Average (over all time steps) and maximum percentages of cells being flagged as troubled cells, for the different troubled-cell indicators, are summarized in Table \ref{table:5} for two different grid sizes and various orders. The computed solution for density obtained at $t=0.038$s using 400 elements while using the SJ indicator and CSWENO limiter for $P^{1}$, $P^{2}$ and $P^{3}$ based DGM is compared and plotted against the exact solution in Figure \ref{fig:BlastWaveSolution}. The error ($|\rho-\rho_{exact}|$) obtained for $P^{1}$, $P^{2}$ and $P^{3}$ based DGM is also plotted in Figure \ref{fig:BlastWaveError}. We are not showing the solutions obtained using other indicators as they are quite similar to the solution obtained with the SJ indicator and stay within an $L^{2}$ norm of $\sim  10^{-14}$. We show the $L^{2}$ difference in density between those solutions in table \ref{table:5New}.
\\
We also show the time history of the flagged troubled cells using the eight different indicators in Figure \ref{fig:BW} using $P^{1}$ based DGM for 200 elements. For all the indicators, the number of flagged troubled cells slightly increase when the order increases. From the tabulated results and the figure, the troubled cell indicators PP and PPL perform in a similar fashion, while the LPR indicator is quite better and the MH indicator is quite bad comparatively. We can also say that, for this problem, the troubled cell indicators FS1, FS2, RH, and SJ out perform the other four indicators and perform in a similar fashion.
\\
\\
\begin{table}[htbp]
\small
\centering
\begin{tabular}{|c|c|c|c|c|c|c|c|c|c|}
%\hline
%\multicolumn{6}{|c|}{$L_{2}$ error for the Riemann problem configuration-10} \\
\hline
 \makecell{No. of \\ Cells} & \makecell{Scheme \\ Indicator} & \multicolumn{2}{|c|}{$P^{1}$} & \multicolumn{2}{|c|}{$P^{2}$} & \multicolumn{2}{|c|}{$P^{3}$} & \multicolumn{2}{|c|}{$P^{4}$} \\
\cline{3-10}
 & & Ave & Max & Ave & Max & Ave & Max & Ave & Max \\
\hline 
\multirow{8}{*}{200} & PP & 7.84 & 12.5 & 8.18 & 13.5 & 8.95 & 14.0 & 9.47 & 14.5 \\
\cline{2-10}
 & SJ & 5.25 & 9.5 & 5.57 & 9.5 & 6.11 & 10.0 & 6.84 & 10.5 \\
\cline{2-10}
 & FS1 & 4.37 & 7.5 & 4.93 & 7.5 & 5.14 & 8.5 & 5.48 & 9.0 \\
\cline{2-10}
 & FS2 & 4.41 & 8.0 & 4.95 & 8.5 & 5.19 & 9.0 & 5.51 & 9.0 \\
\cline{2-10}
 & LPR & 5.95 & 10.5 & 6.44 & 11.0 & 6.92 & 11.0 & 7.11 & 11.5\\
\cline{2-10}
 & RH & 4.14 & 7.0 & 4.72 & 7.5 & 5.05 & 8.0 & 5.39 & 8.5 \\
 \cline{2-10}
 & PPL & 7.46 & 11.5 & 8.02 & 13.0 & 8.64 & 13.5 & 9.11 & 14.0 \\
\cline{2-10}
 & MH & 10.31 & 18.5 & 10.85 & 18.5 & 11.24 & 19.0 & 11.73 & 19.5\\
\hline
\multirow{8}{*}{400} & PP & 7.26 & 12.0 & 7.83 & 12.75 & 8.11 & 13.0 & 8.85 & 13.5 \\
\cline{2-10}
 & SJ & 5.01 & 9.25 & 5.34 & 9.75 & 5.94 & 10.0 & 6.32 & 10.25 \\
\cline{2-10}
 & FS1 & 4.24 & 7.25 & 4.88 & 7.5 & 5.10 & 8.0 & 5.44 & 8.75 \\
\cline{2-10}
 & FS2 & 4.25 & 7.75 & 4.91 & 8.0 & 5.08 & 8.0 & 5.43 & 9.0 \\
\cline{2-10}
 & LPR & 5.65 & 10.25 & 6.13 & 10.5 & 6.58 & 11.0 & 6.89 & 11.25\\
\cline{2-10}
 & RH & 3.96 & 6.75 & 4.54 & 7.25 & 4.88 & 7.75 & 5.26 & 8.25 \\
 \cline{2-10}
 & PPL & 7.04 & 12.0 & 7.56 & 12.5 & 7.94 & 12.75 & 8.43 & 13.25 \\
\cline{2-10}
 & MH & 10.02 & 18.25 & 10.63 & 18.75 & 10.97 & 19.0 & 11.34 & 19.25\\
\hline
\end{tabular}
\caption{Average (marked as Ave) and maximum (marked as Max) percentages of cells flagged as troubled cells subject to different troubled-cell indicators for the blast wave problem for various orders.}
\label{table:5}
\end{table}

\begin{table}[htbp]
%\small
\centering
\begin{tabular}{|c|c|c|c|c|}
%\hline
%\multicolumn{6}{|c|}{$L_{2}$ error for the Riemann problem configuration-10} \\
\hline
Scheme Indicator & \multicolumn{4}{|c|}{\makecell{$L^{2}$ difference in density in comparison with \\ SJ indicator solution for various orders}} \\ \cline{2-5}
  & $P^{1}$ & $P^{2}$ & $P^{3}$ & $P^{4}$ \\ 
\hline
PP & 1.5E-14 & 1.6E-14 & 3.8E-14 & 5.6E-15\\
\hline
FS1 & 2.4E-14 & 2.9E-14 & 7.1E-14 & 3.7E-15\\
\hline
FS2 & 2.4E-14 & 2.9E-14 & 7.1E-14 & 3.7E-15\\
\hline
LPR & 3.8E-14 & 4.8E-14 & 4.6E-14 & 4.6E-15\\
\hline
RH & 5.6E-14 & 5.4E-14 & 6.5E-14 & 5.9E-15\\
\hline
PPL & 4.2E-14 & 6.6E-14 & 4.3E-14 & 7.1E-15\\
\hline
MH & 3.9E-14 & 7.5E-14 & 7.2E-14 & 3.9E-15\\
\hline
\end{tabular}
\caption{$L^{2}$ difference in density between solution obtained using CSWENO limiter and SJ indicator (shown in Figure \ref{fig:BlastWaveCSWENO}) and the solution obtained using other indicators for the blast wave problem using 200 elements}
\label{table:5New}
\end{table}

\begin{figure}[htbp]
  \centering
  \subfloat[Density contours for the solution at t=0.2 with $P^{1}$ based DGM]{\label{fig:P1DGMDMRCSWENOUnstructured}\includegraphics[width=0.45\textwidth]{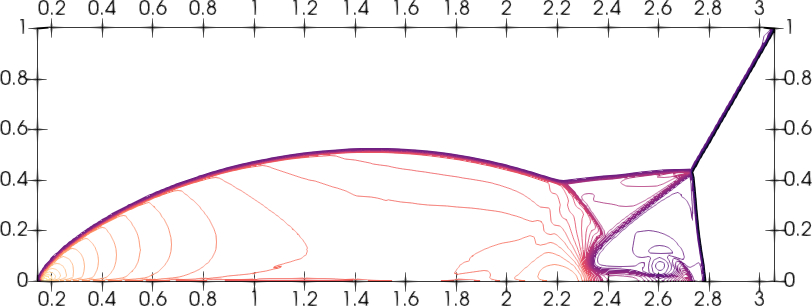}}
  \subfloat[Density contours for the solution at t=0.2 with $P^{2}$ based DGM]{\label{fig:P2DGMDMRCSWENOUnstructured}\includegraphics[width=0.45\textwidth]{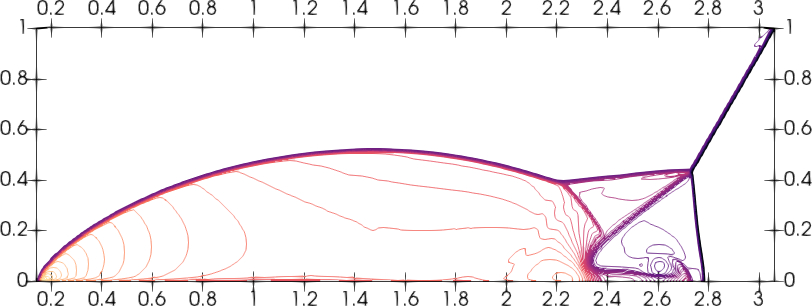}}\hfill
  \subfloat[Density contours for the solution at t=0.2 with $P^{3}$ based DGM]{\label{fig:P3DGMDMRCSWENOUnstructured}\includegraphics[width=0.45\textwidth]{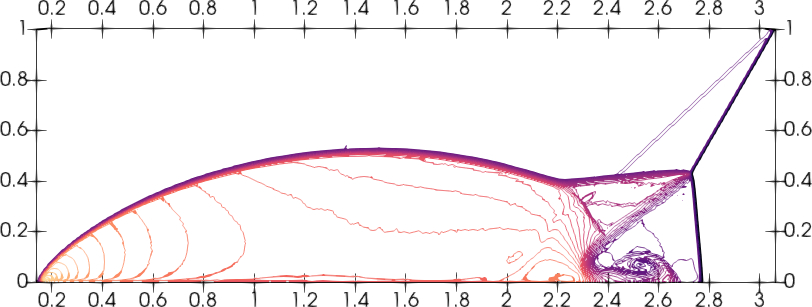}}
  \subfloat[Density contours for the solution at t=0.2 with $P^{4}$ based DGM]{\label{fig:P4DGMDMRCSWENOUnstructured}\includegraphics[width=0.45\textwidth]{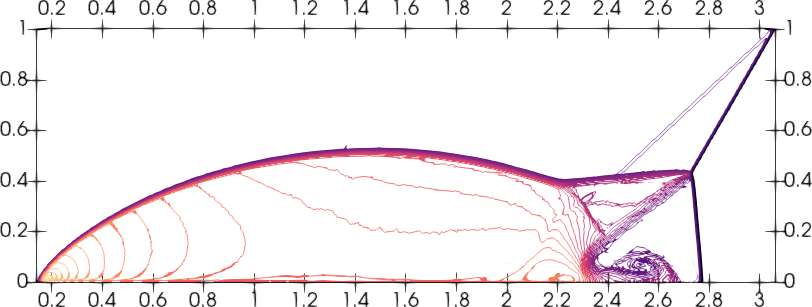}}\hfill
  \subfloat[Density Range]{\label{fig:DensityRangeDMRCSWENOUnstructured}\includegraphics[width=0.9\textwidth]{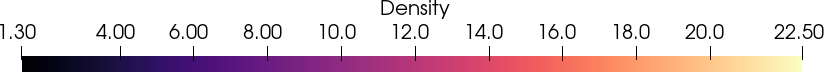}}
  \caption{50 equally spaced density contours for double Mach reflection solution using the CSWENO limiter with triangular cells for $P^{1}$, $P^{2}$, $P^{3}$ and $P^{4}$ based DGM for a mesh containing 370,046 triangles}
  \label{fig:CSWENODMRSolutionUnstructured}
\end{figure}

\begin{figure}[htbp]
  \centering
  \subfloat[PP Indicator]{\label{fig:DMRPP}\includegraphics[width=0.32\textwidth]{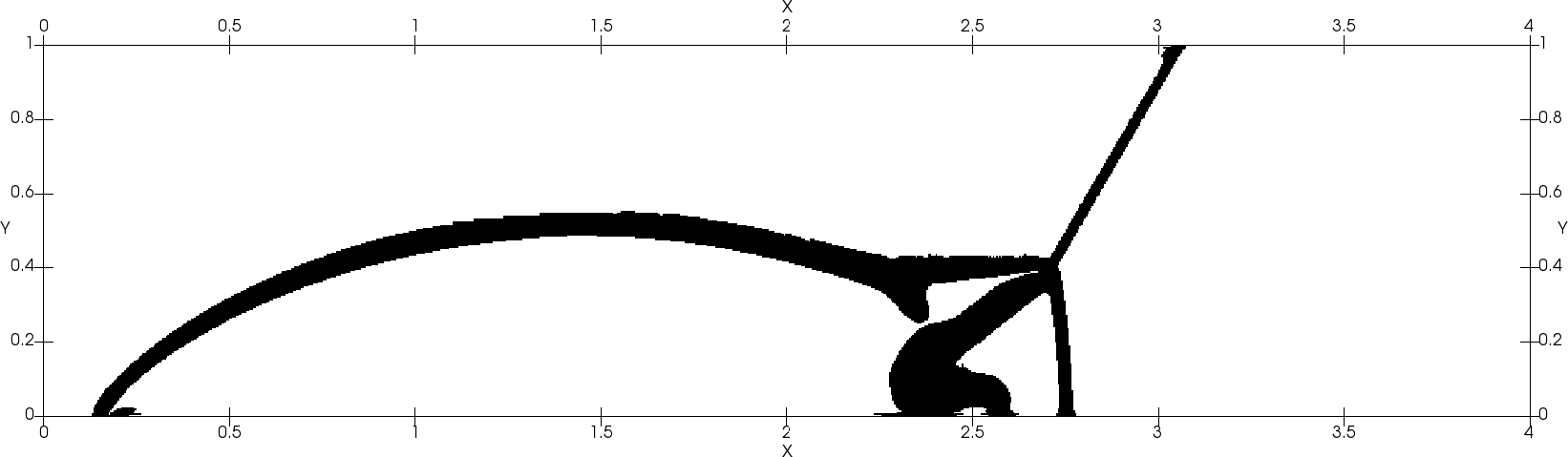}}
  \subfloat[SJ Indicator]{\label{fig:DMRSJ}\includegraphics[width=0.32\textwidth]{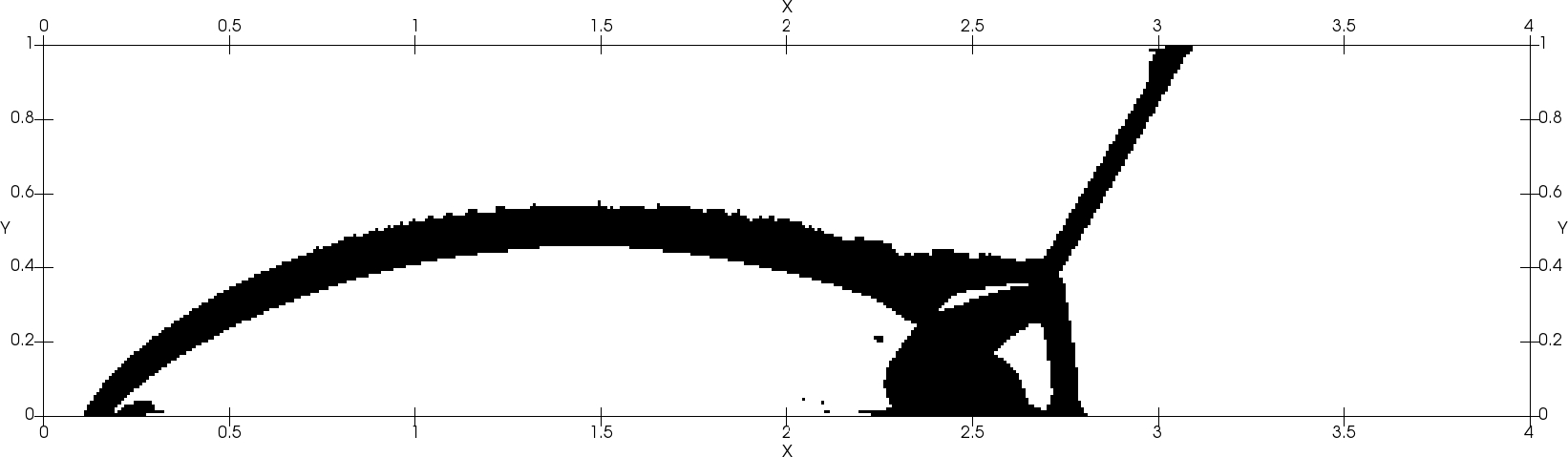}}
  \subfloat[FS1 Indicator]{\label{fig:DMRFS1}\includegraphics[width=0.32\textwidth]{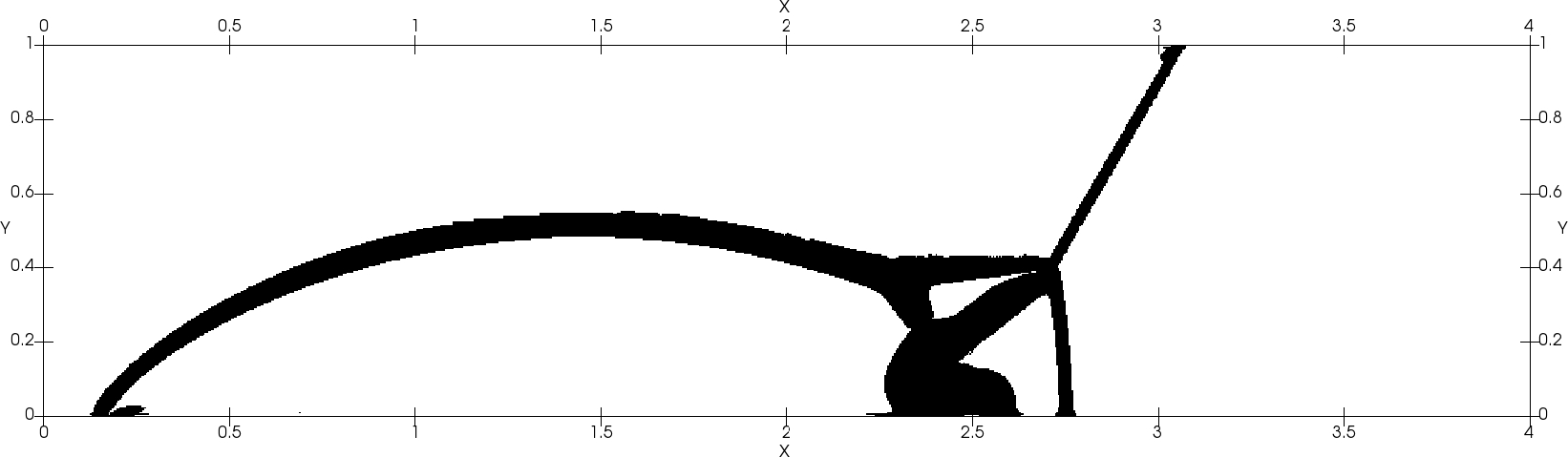}}\hfill
  \subfloat[FS2 Indicator]{\label{fig:DMRFS2}\includegraphics[width=0.32\textwidth]{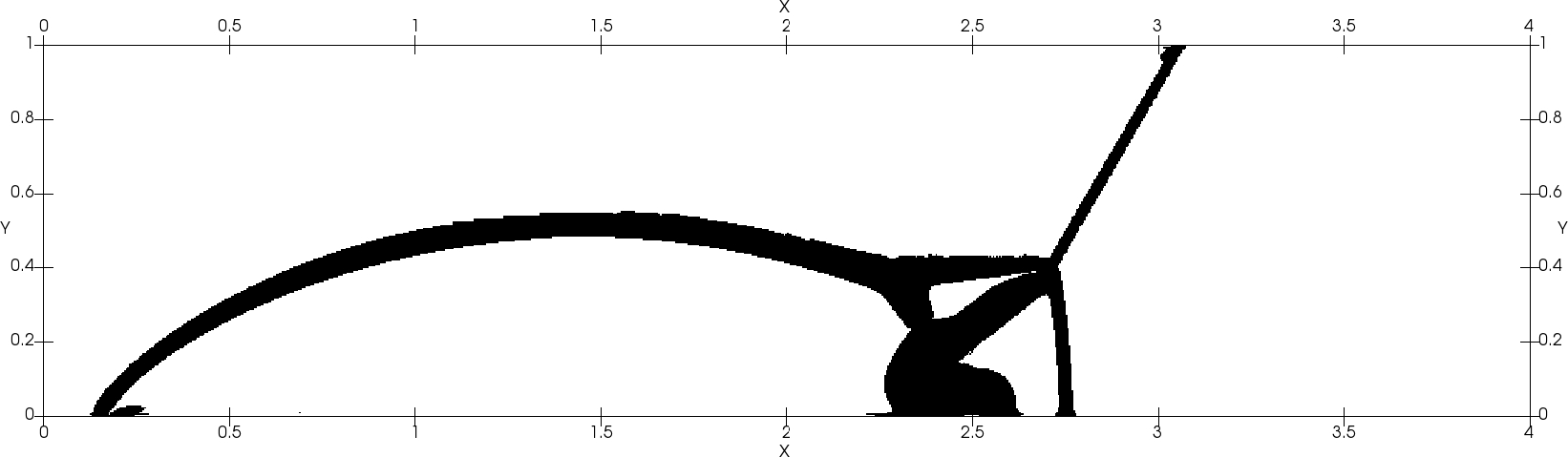}}
  \subfloat[LPR Indicator]{\label{fig:DMRLPR}\includegraphics[width=0.32\textwidth]{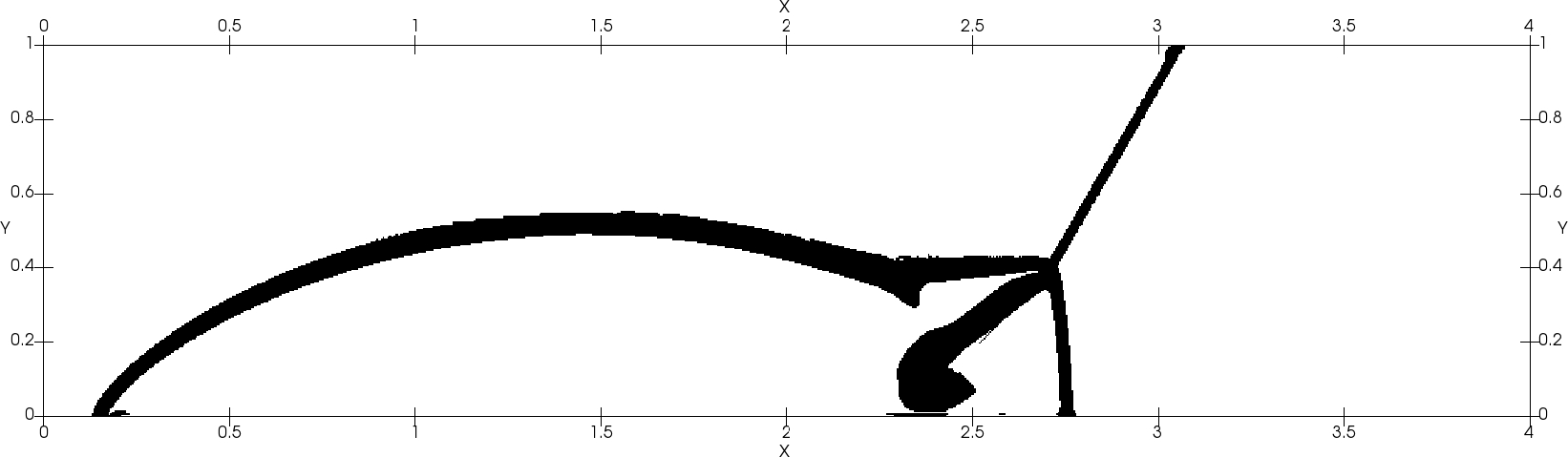}}
  \subfloat[RH Indicator]{\label{fig:DMRRH}\includegraphics[width=0.32\textwidth]{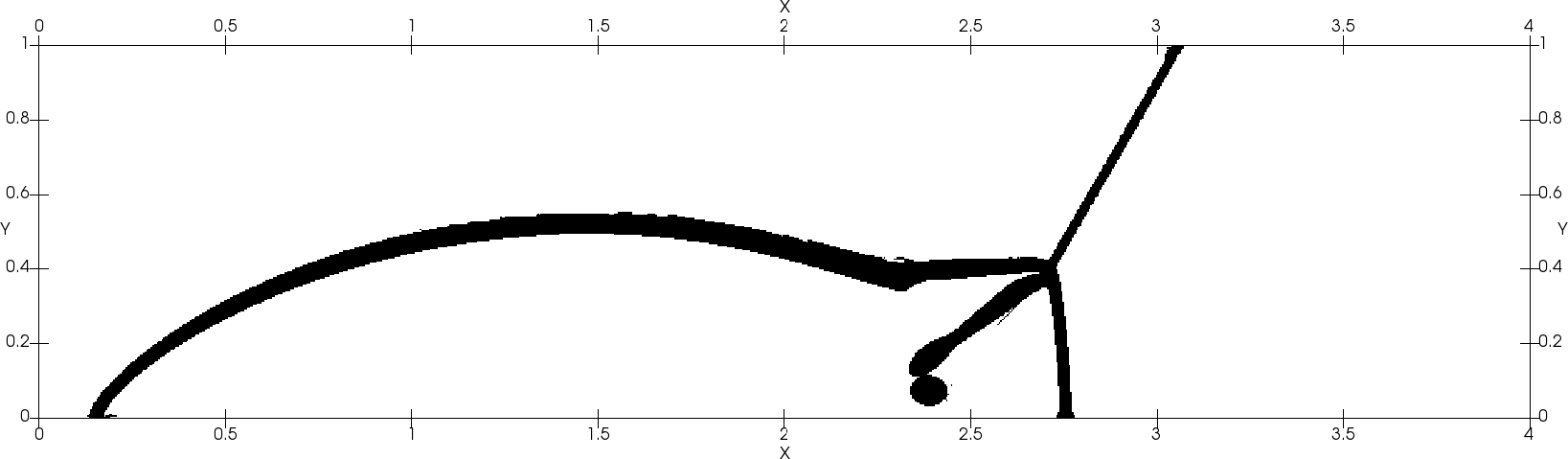}}\hfill
  \subfloat[PPL Indicator]{\label{fig:DMRPPL}\includegraphics[width=0.32\textwidth]{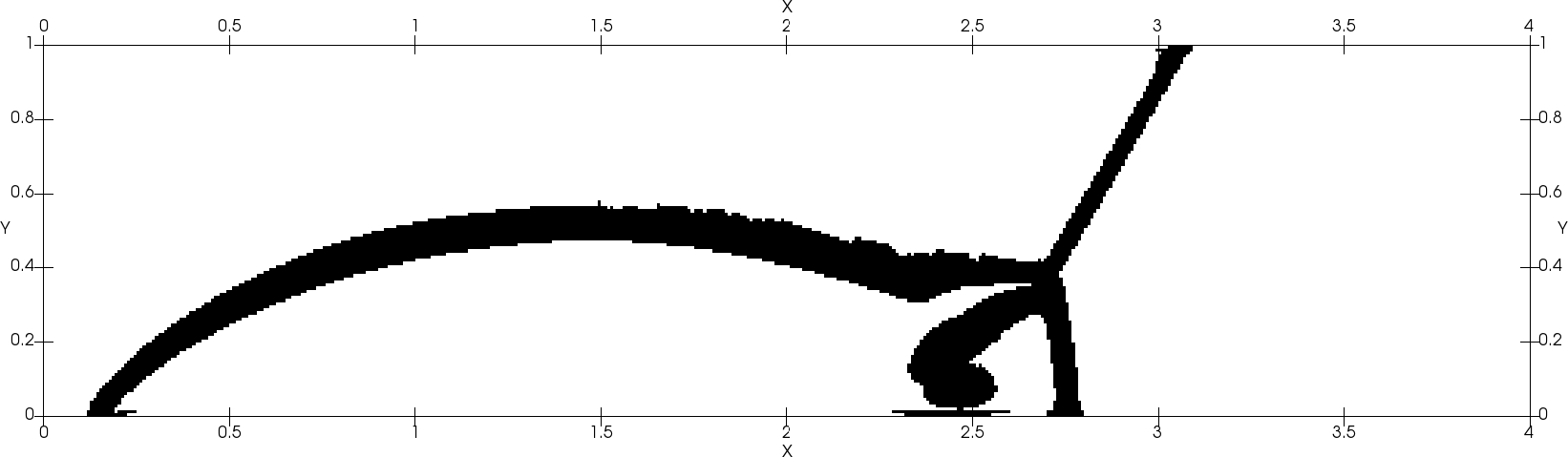}}
  \subfloat[MH Indicator]{\label{fig:DMRMH}\includegraphics[width=0.32\textwidth]{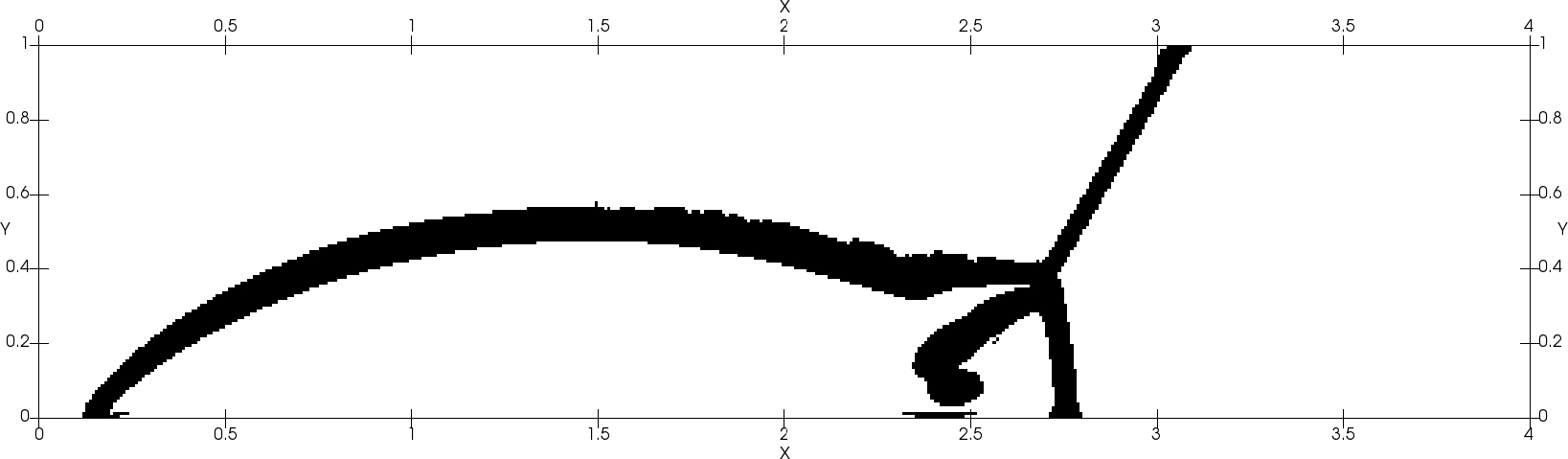}}\hfill
  \caption{The time history of flagged troubled cells of the double Mach reflection problem for the two-dimensional Euler equations, simulated until $t=0.2$ with $480\times 120$ elements and $P^{1}$ based DGM}
  \label{fig:DMR}
\end{figure}

\noindent \textbf{Test Problem 6 (Double Mach reflection)\cite{wc}:} We solve the two-dimensional Euler equations for an ideal gas given by

\begin{equation}\label{2dEulerEquations}
\frac{\partial Q}{\partial t} + \frac{\partial F(Q)}{\partial x} + \frac{\partial G(Q)}{\partial y} = 0
\end{equation}

\noindent in the domain $[0,4]\times [0,1]$ where $Q = (\rho, \rho u, \rho v, E)^{T}$, $F(Q) = uQ + (0, p, 0, pu)^{T}$, $G(Q) = vQ + (0, 0, p, pv)^{T}$, $p = (\gamma -1)(E-\frac{1}{2}\rho (u^{2}+v^{2}))$ and $\gamma = 1.4$. Here, $\rho$ is the density, $(u,v)$ is the velocity, $E$ is the total energy and $p$ is the pressure. Initially, a right moving Mach 10 shock is positioned at $x = 1/6$, $y = 0$ and it makes an angle $60^{\circ}$ with the $x$-axis. For the bottom boundary, we impose the exact post shock conditions from $x = 0$ to $x = 1/6$ and for the rest of the $x$-axis, we use reflective boundary conditions. For the top boundary, we set conditions to describe the exact motion of a Mach 10 shock. We compute the solution upto time $t = 0.2$ for both structured and unstructured grids. For structured grids, we have used two different uniform meshes with $480\times 120$ and $960\times 240$ elements. For unstructured grids (generated using Gmsh 4.6.0 software \cite{GMshGR}), we have used meshes with 185023 and 370046 triangles. Average (over all time steps) and maximum percentages of cells being flagged as troubled cells, for the different troubled-cell indicators, are summarized for structured and unstructured grids respectively in Tables \ref{table:6} and \ref{table:7} for the two different grid sizes and various orders. The density contours for the solution obtained using the SJ indicator and CSWENO limiter for $P^{1}$, $P^{2}$, $P^{3}$ and $P^{4}$ based DGM using 370046 triangles are shown in Figure \ref{fig:CSWENODMRSolutionUnstructured}. We are not showing the solutions obtained using other indicators as they are quite similar to the solution obtained with the SJ indicator and stay within an $L^{2}$ norm of $\sim  10^{-14}$. We show the $L^{2}$ difference in density between those solutions in table \ref{table:7New}.
\\
We also show the troubled cell profile using the eight different indicators in Figure \ref{fig:DMR} using $P^{1}$ based DGM for $480\times 120$ elements at $t=0.2$. For all the indicators, the number of flagged troubled cells slightly increase when the order increases. From the tabulated results and the figure, the troubled cell indicators PP, SJ, FS1, FS2, MH, and PPL perform in a similar fashion, while the LPR indicator is quite better and the RH indicator is the best comparatively. From the tables, we also observe that the performance of the troubled-cell indicators for structured and unstructured grids is quite similar.
\\
\\
\begin{table}[htbp]
\small
\centering
\begin{tabular}{|c|c|c|c|c|c|c|c|c|c|}
%\hline
%\multicolumn{6}{|c|}{$L_{2}$ error for the Riemann problem configuration-10} \\
\hline
 \makecell{No. of \\ Cells} & \makecell{Scheme \\ Indicator} & \multicolumn{2}{|c|}{$P^{1}$} & \multicolumn{2}{|c|}{$P^{2}$} & \multicolumn{2}{|c|}{$P^{3}$} & \multicolumn{2}{|c|}{$P^{4}$} \\
\cline{3-10}
 & & Ave & Max & Ave & Max & Ave & Max & Ave & Max \\
\hline 
\multirow{8}{*}{$480\times 120$} & PP & 4.04 & 6.75 & 4.43 & 6.97 & 4.87 & 7.16 & 5.13 & 7.22 \\
\cline{2-10}
 & SJ & 4.12 & 6.81 & 4.89 & 6.92 & 5.14 & 7.24 & 5.37 & 7.35 \\
\cline{2-10}
 & FS1 & 4.03 & 6.72 & 4.35 & 6.99 & 4.74 & 7.08 & 4.98 & 7.18 \\
\cline{2-10}
 & FS2 & 4.03 & 6.72 & 4.35 & 6.99 & 4.74 & 7.08 & 4.98 & 7.18 \\
\cline{2-10}
 & LPR & 3.65 & 5.99 & 3.89 & 6.29 & 4.28 & 6.74 & 4.64 & 6.87 \\
\cline{2-10}
 & RH & 3.23 & 5.82 & 3.65 & 5.94 & 3.89 & 6.09 & 4.12 & 6.32 \\
 \cline{2-10}
 & PPL & 3.89 & 6.27 & 4.02 & 6.75 & 4.59 & 7.02 & 4.78 & 7.10 \\
\cline{2-10}
 & MH & 3.86 & 6.37 & 4.05 & 6.68 & 4.56 & 6.99 & 4.76 & 7.08 \\
\hline
\multirow{8}{*}{$960\times 240$} & PP & 3.96 & 6.54 & 4.24 & 6.88 & 4.56 & 6.90 & 5.02 & 6.97 \\
\cline{2-10}
 & SJ & 3.94 & 6.73 & 4.33 & 6.88 & 5.11 & 7.09 & 5.22 & 7.22 \\
\cline{2-10}
 & FS1 & 3.82 & 6.53 & 4.14 & 6.73 & 4.53 & 6.96 & 4.76 & 7.12 \\
\cline{2-10}
 & FS2 & 3.82 & 6.53 & 4.14 & 6.73 & 4.53 & 6.96 & 4.76 & 7.12 \\
\cline{2-10}
 & LPR & 3.54 & 5.86 & 3.76 & 6.18 & 4.04 & 6.57 & 4.53 & 6.81 \\
\cline{2-10}
 & RH & 3.09 & 5.75 & 3.58 & 5.91 & 3.77 & 6.03 & 3.98 & 6.17 \\
 \cline{2-10}
 & PPL & 3.73 & 6.14 & 3.94 & 6.70 & 4.16 & 6.89 & 4.54 & 7.00 \\
\cline{2-10}
 & MH & 3.69 & 6.08 & 3.84 & 6.40 & 4.08 & 6.76 & 4.58 & 7.06 \\
\hline
\end{tabular}
\caption{Average (marked as Ave) and maximum (marked as Max) percentages of cells flagged as troubled cells subject to different troubled-cell indicators for the double Mach reflection problem for various orders using structured grids.}
\label{table:6}
\end{table}

\begin{table}[htbp]
\small
\centering
\begin{tabular}{|c|c|c|c|c|c|c|c|c|c|}
%\hline
%\multicolumn{6}{|c|}{$L_{2}$ error for the Riemann problem configuration-10} \\
\hline
 \makecell{No. of \\ Cells} & \makecell{Scheme \\ Indicator} & \multicolumn{2}{|c|}{$P^{1}$} & \multicolumn{2}{|c|}{$P^{2}$} & \multicolumn{2}{|c|}{$P^{3}$} & \multicolumn{2}{|c|}{$P^{4}$} \\
\cline{3-10}
 & & Ave & Max & Ave & Max & Ave & Max & Ave & Max \\
\hline 
\multirow{8}{*}{185023} & PP & 4.02 & 6.71 & 4.47 & 6.99 & 4.84 & 7.12 & 5.18 & 7.24 \\
\cline{2-10}
 & SJ & 4.16 & 6.88 & 4.87 & 6.90 & 5.19 & 7.22 & 5.34 & 7.32 \\
\cline{2-10}
 & FS1 & 4.06 & 6.76 & 4.37 & 6.96 & 4.77 & 7.04 & 4.93 & 7.15 \\
\cline{2-10}
 & FS2 & 4.06 & 6.76 & 4.37 & 6.96 & 4.77 & 7.04 & 4.93 & 7.15 \\
\cline{2-10}
 & LPR & 3.67 & 5.94 & 3.88 & 6.23 & 4.24 & 6.72 & 4.66 & 6.83 \\
\cline{2-10}
 & RH & 3.26 & 5.86 & 3.62 & 5.97 & 3.84 & 6.02 & 4.17 & 6.39 \\
 \cline{2-10}
 & PPL & 3.83 & 6.24 & 4.08 & 6.79 & 4.55 & 7.04 & 4.75 & 7.15 \\
\cline{2-10}
 & MH & 3.89 & 6.34 & 4.03 & 6.65 & 4.59 & 6.96 & 4.74 & 7.06 \\
\hline
\multirow{8}{*}{370046} & PP & 3.94 & 6.59 & 4.28 & 6.83 & 4.54 & 6.96 & 5.06 & 6.94 \\
\cline{2-10}
 & SJ & 3.92 & 6.78 & 4.36 & 6.83 & 5.15 & 7.04 & 5.24 & 7.27 \\
\cline{2-10}
 & FS1 & 3.85 & 6.58 & 4.18 & 6.77 & 4.58 & 6.99 & 4.74 & 7.23 \\
\cline{2-10}
 & FS2 & 3.85 & 6.58 & 4.18 & 6.77 & 4.58 & 6.99 & 4.74 & 7.23 \\
\cline{2-10}
 & LPR & 3.51 & 5.83 & 3.78 & 6.13 & 4.08 & 6.54 & 4.52 & 6.87 \\
\cline{2-10}
 & RH & 3.03 & 5.73 & 3.57 & 5.93 & 3.71 & 6.04 & 3.95 & 6.13 \\
 \cline{2-10}
 & PPL & 3.72 & 6.11 & 3.97 & 6.73 & 4.12 & 6.92 & 4.52 & 7.11 \\
\cline{2-10}
 & MH & 3.66 & 6.03 & 3.87 & 6.45 & 4.04 & 6.73 & 4.54 & 7.08 \\
\hline
\end{tabular}
\caption{Average (marked as Ave) and maximum (marked as Max) percentages of cells flagged as troubled cells subject to different troubled-cell indicators for the double Mach reflection problem for various orders using unstructured triangles.}
\label{table:7}
\end{table}

\begin{table}[htbp]
%\small
\centering
\begin{tabular}{|c|c|c|c|c|}
%\hline
%\multicolumn{6}{|c|}{$L_{2}$ error for the Riemann problem configuration-10} \\
\hline
Scheme Indicator & \multicolumn{4}{|c|}{\makecell{$L^{2}$ difference in density in comparison with \\ SJ indicator solution for various orders}} \\ \cline{2-5}
  & $P^{1}$ & $P^{2}$ & $P^{3}$ & $P^{4}$ \\ 
\hline
PP & 2.4E-14 & 2.9E-14 & 4.8E-14 & 6.1E-15\\
\hline
FS1 & 3.3E-14 & 3.4E-14 & 4.2E-14 & 6.8E-15\\
\hline
FS2 & 3.3E-14 & 3.4E-14 & 4.2E-14 & 6.8E-15\\
\hline
LPR & 1.7E-14 & 5.4E-14 & 5.7E-14 & 7.6E-15\\
\hline
RH & 4.9E-14 & 1.7E-14 & 7.4E-14 & 6.3E-15\\
\hline
PPL & 6.6E-14 & 3.9E-14 & 8.2E-14 & 3.2E-15\\
\hline
MH & 7.2E-14 & 6.5E-14 & 9.5E-14 & 5.5E-15\\
\hline
\end{tabular}
\caption{$L^{2}$ difference in density between solution obtained using CSWENO limiter and SJ indicator (shown in Figure \ref{fig:CSWENODMRSolutionUnstructured}) and the solution obtained using other indicators for the double Mach reflection problem using 370046 triangles}
\label{table:7New}
\end{table}

\noindent \textbf{Test Problem 7 (2D Riemann Problem)\cite{ll}:} We solve the two-dimensional Euler equations for an ideal gas given by \eqref{2dEulerEquations} in the domain $[0,1]\times [0,1]$ for the 2D Riemann problem configurations 4 and 12 as given by the nomenclature in \cite{ll}. The initial conditions for configurations 4 and 12 are given respectively as

\begin{flalign}\label{2DRiemannInitialConfig04}
 (\rho,u,v,p)(x,y,0) = \begin{cases}
                            (1.1,0,0,1.1) \quad \text{if $x \geq 0.5$ and $y \geq 0.5$} \\
                            (0.5065,0.8939,0,0.35) \quad \text{if $x<0.5$ and $y \geq 0.5$} \\
                            (1.1,0.8939,0.8939,1.1) \quad \text{if $x<0.5$ and $y<0.5$} \\
                            (0.5065,0,0.8939,0.35) \quad \text{otherwise}
                        \end{cases}
\end{flalign}

\begin{flalign}\label{2DRiemannInitialConfig12}
 (\rho,u,v,p)(x,y,0) = \begin{cases}
                            (0.5313,0,0,0.4) \quad \text{if $x \geq 0.5$ and $y \geq 0.5$} \\
                            (1,0.7276,0,1) \quad \text{if $x<0.5$ and $y \geq 0.5$} \\
                            (0.8,0,0,1) \quad \text{if $x<0.5$ and $y<0.5$} \\
                            (1,0,0.7276,1) \quad \text{otherwise}
                        \end{cases}
\end{flalign}

\begin{figure}[htbp]
  \centering
  \subfloat[Density contours for the solution at t=0.25 with $P^{1}$ based DGM]{\label{fig:P1DGM2DRiemannConfig04CSWENOUnstructured}\includegraphics[width=0.45\textwidth]{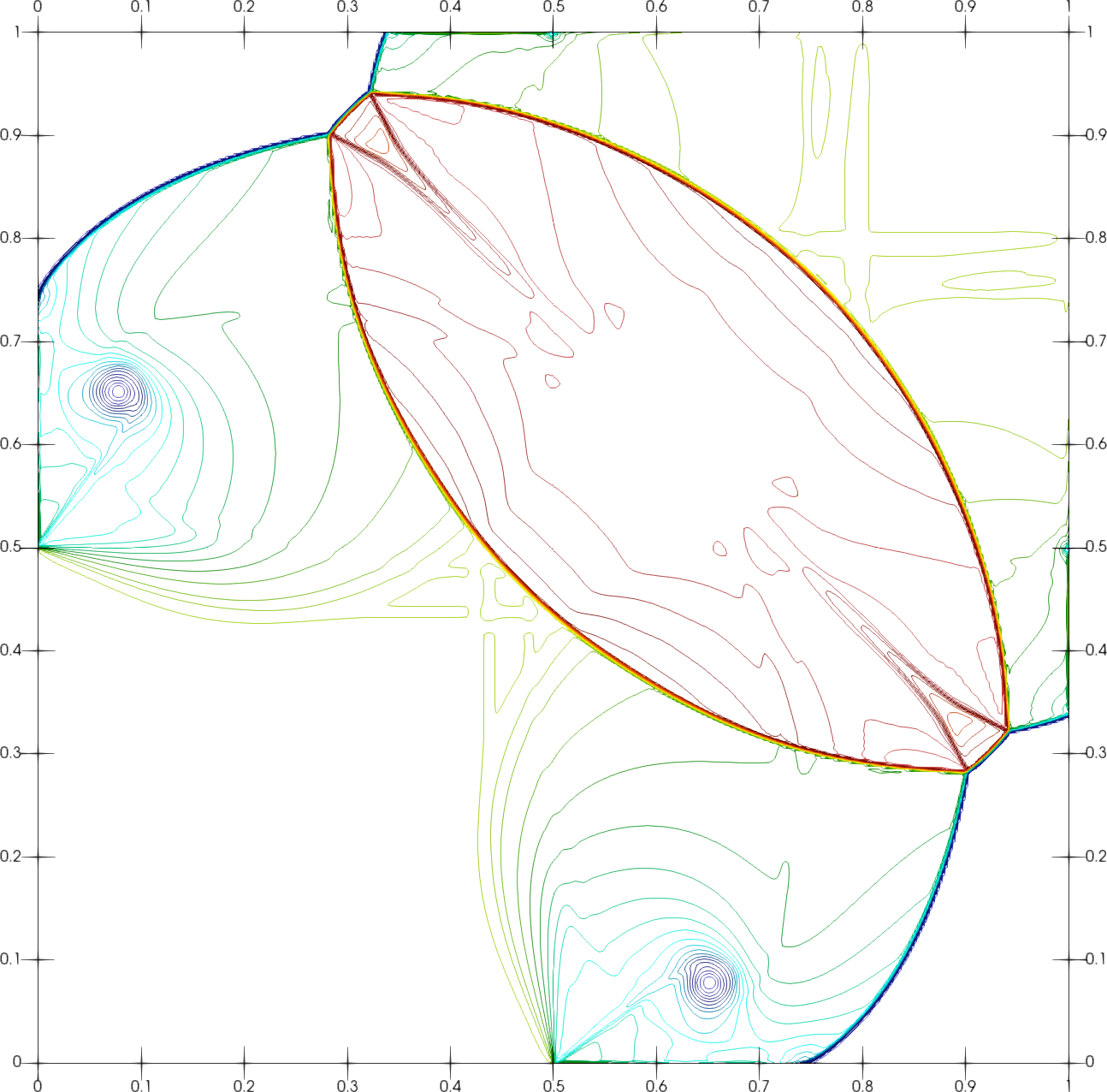}}
  \subfloat[Density contours for the solution at t=0.25 with $P^{2}$ based DGM]{\label{fig:P2DGM2DRiemannConfig04CSWENOUnstructured}\includegraphics[width=0.45\textwidth]{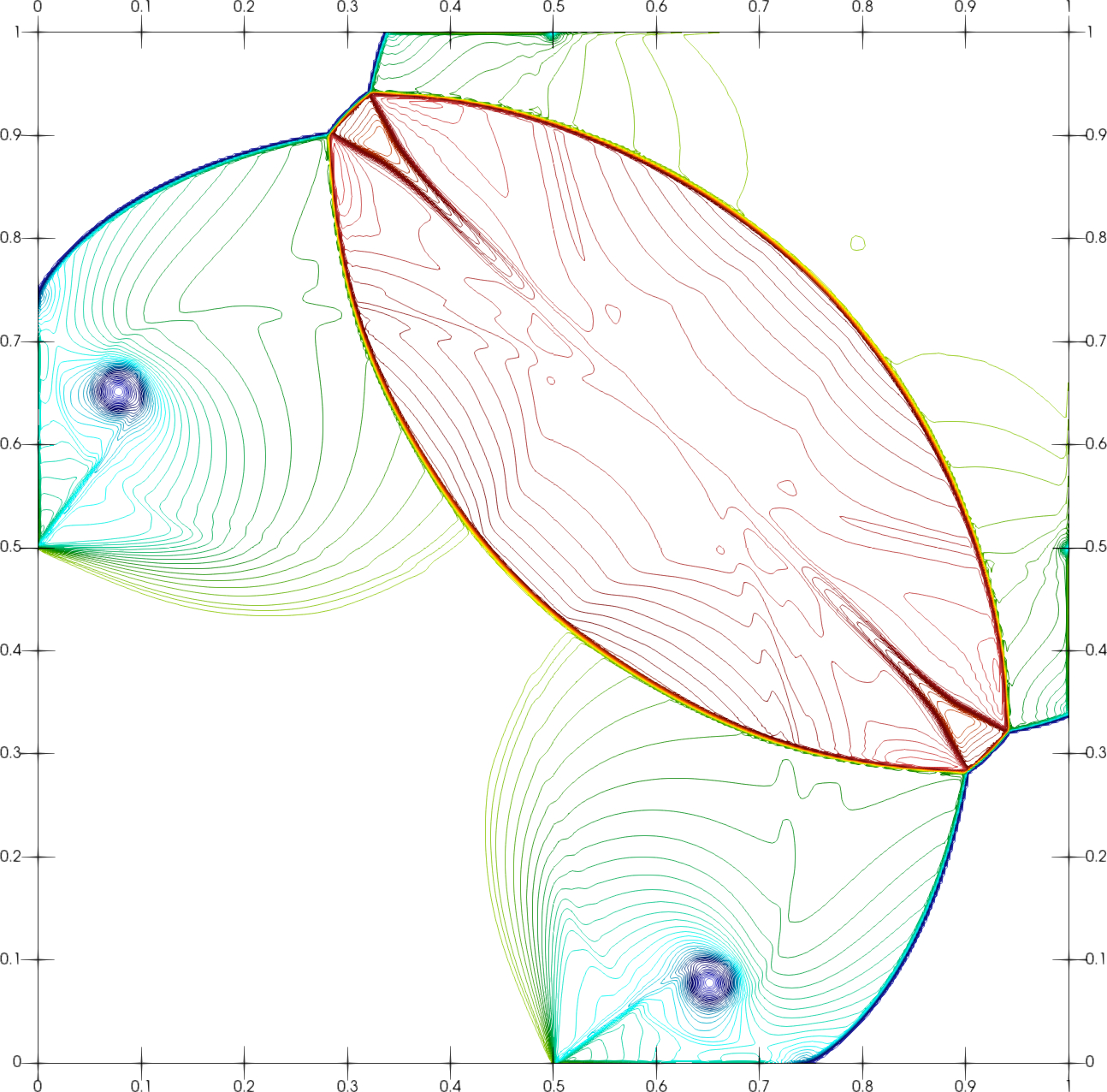}}\hfill
  \subfloat[Density contours for the solution at t=0.25 with $P^{3}$ based DGM]{\label{fig:P3DGM2DRiemannConfig04CSWENOUnstructured}\includegraphics[width=0.45\textwidth]{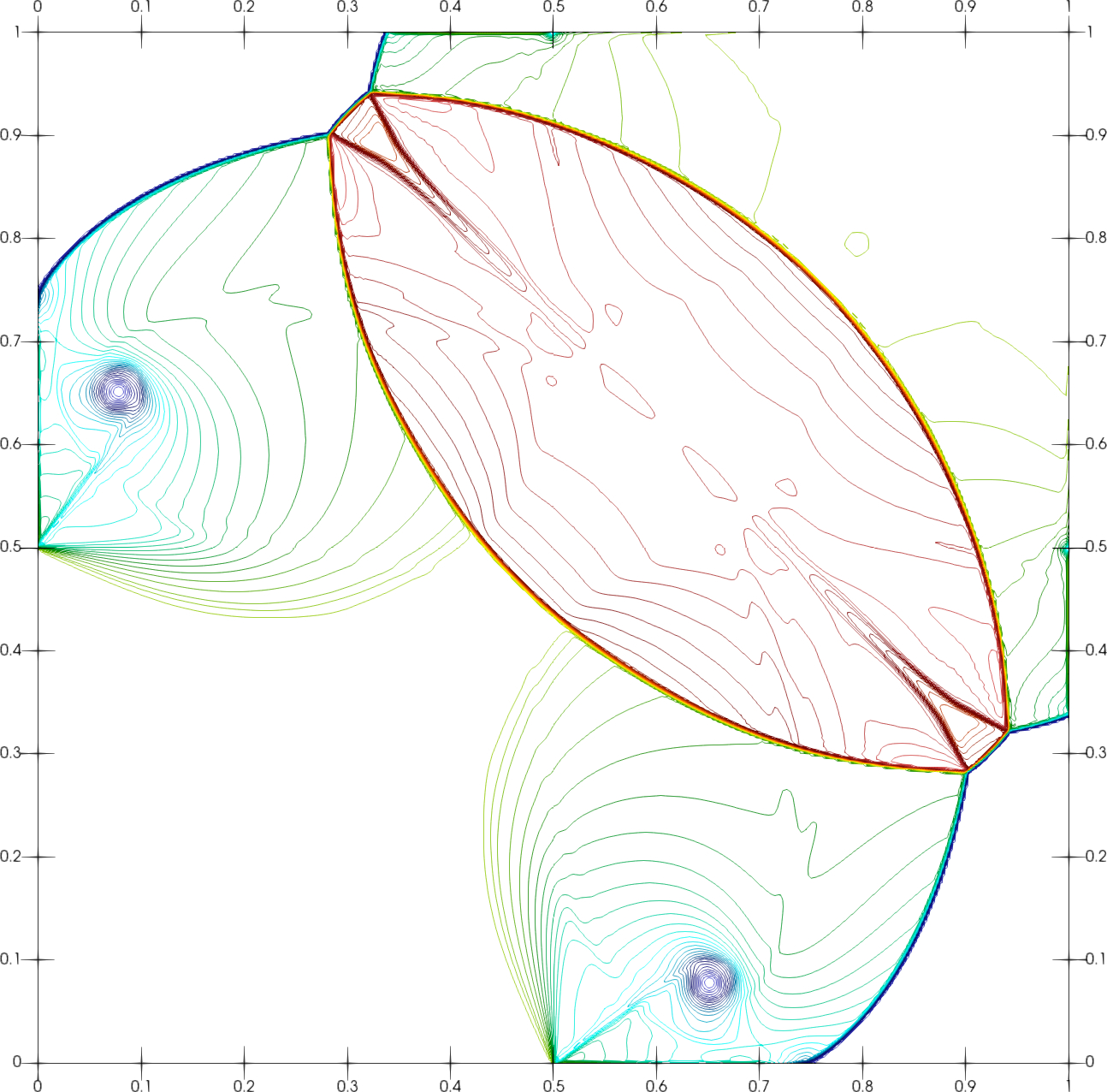}}
  \subfloat[Density contours for the solution at t=0.25 with $P^{4}$ based DGM]{\label{fig:P4DGM2DRiemannConfig04CSWENOUnstructured}\includegraphics[width=0.45\textwidth]{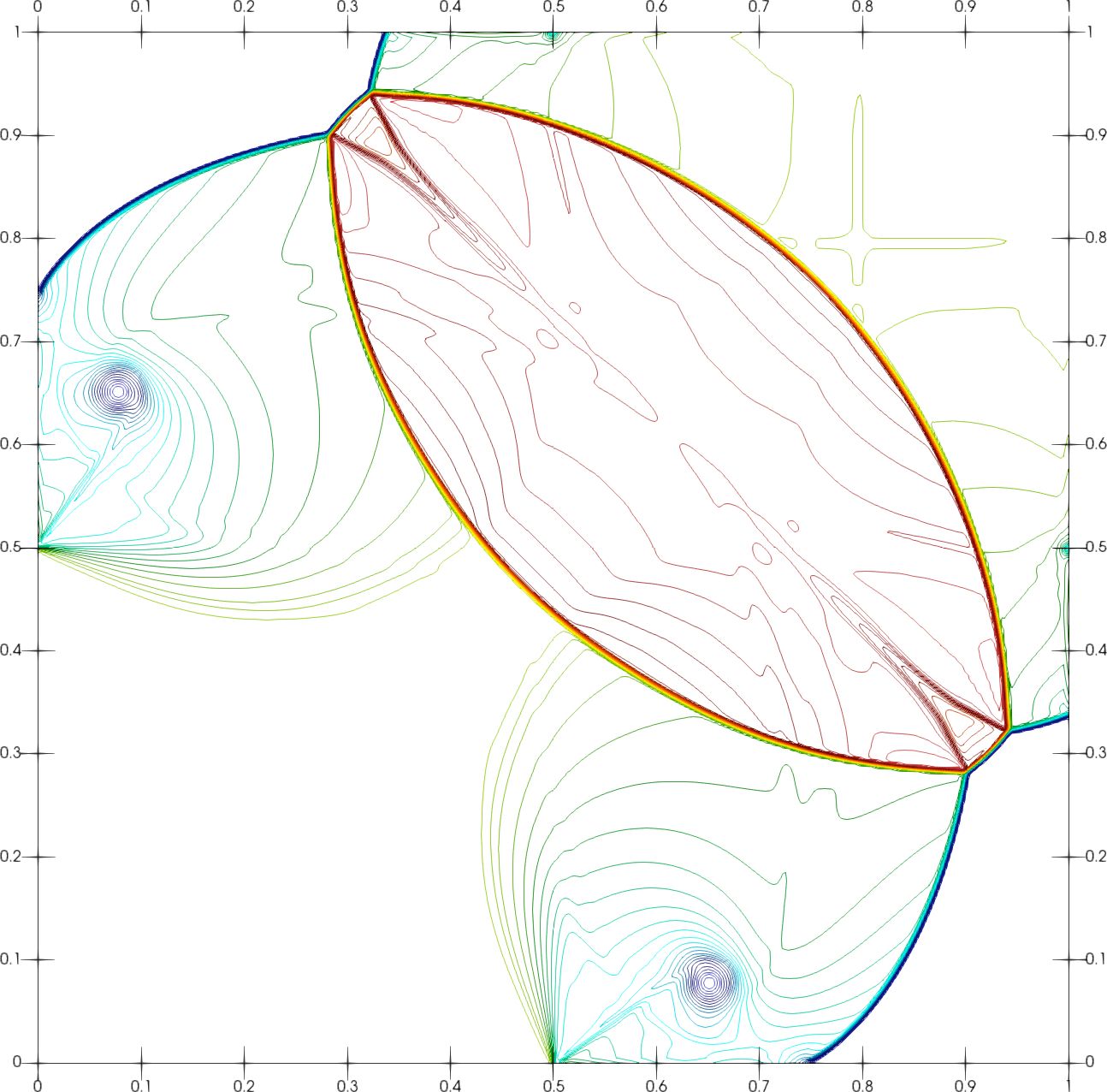}}\hfill
  \subfloat[Density Range]{\label{fig:DensityRange2DRiemannConfig04CSWENOUnstructured}\includegraphics[width=0.9\textwidth]{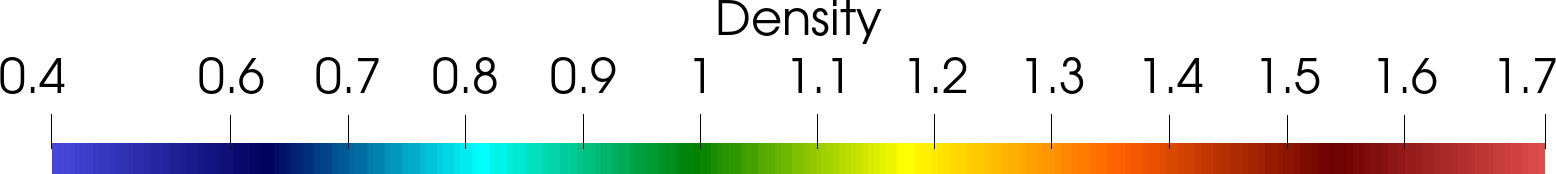}}
  \caption{50 equally spaced density contours for solution at $t=0.25$ for 2D Riemann problem configuration 4 using the CSWENO limiter for $P^{1}$, $P^{2}$, $P^{3}$, and $P^{4}$ based DGM for a mesh containing 92552 triangles}
  \label{fig:CSWENO2DRiemannConfig04SolutionUnstructured}
\end{figure}

\begin{figure}[htbp]
  \centering
  \subfloat[PP Indicator]{\label{fig:2DRiemannConfig04PP}\includegraphics[width=0.32\textwidth]{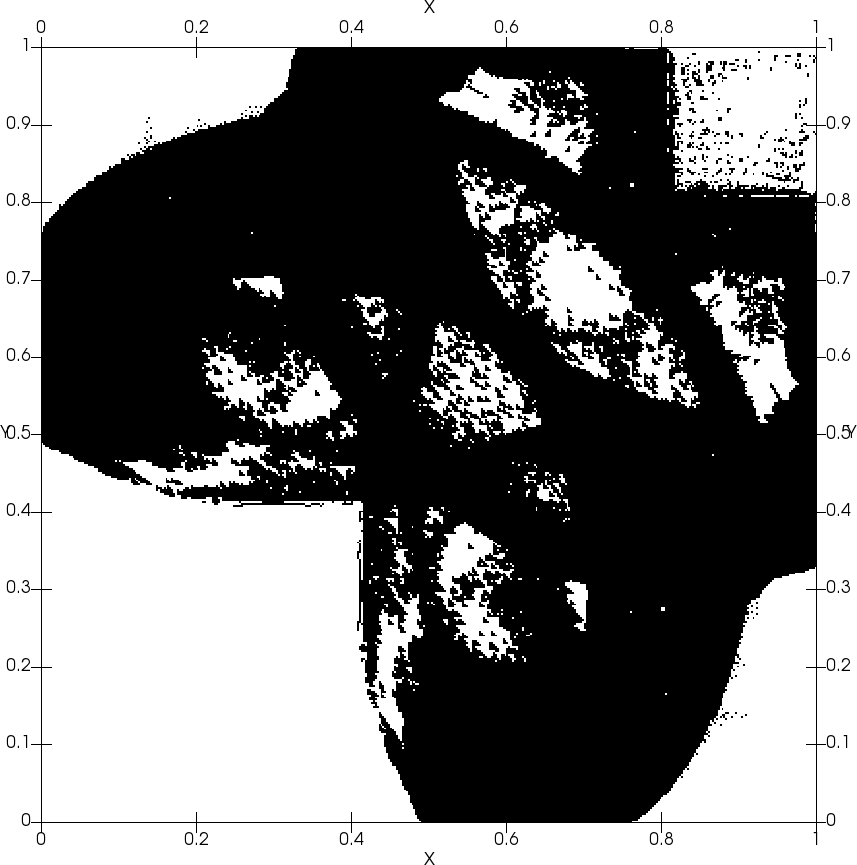}}
  \subfloat[SJ Indicator]{\label{fig:2DRiemannConfig04SJ}\includegraphics[width=0.32\textwidth]{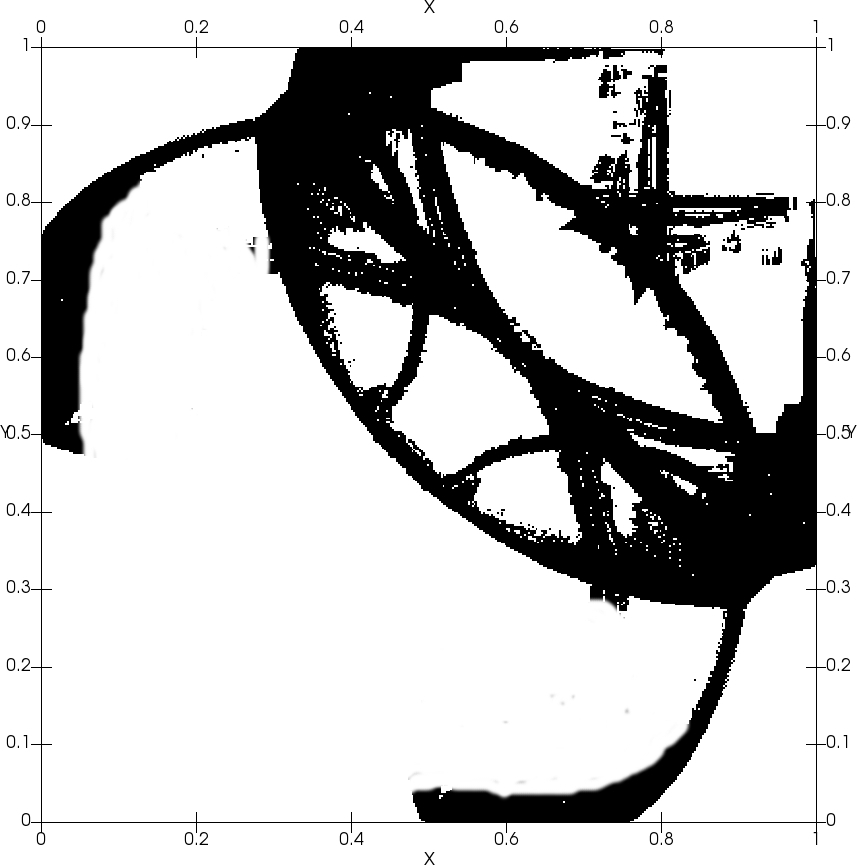}}
  \subfloat[FS1 Indicator]{\label{fig:2DRiemannConfig04FS1}\includegraphics[width=0.32\textwidth]{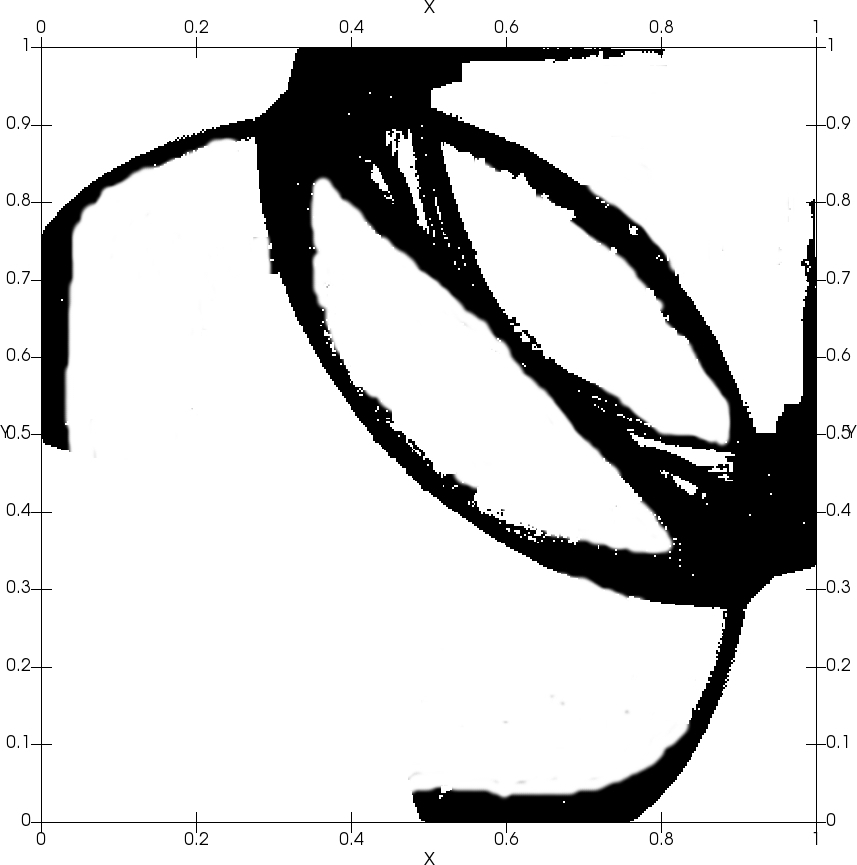}}\hfill
  \subfloat[FS2 Indicator]{\label{fig:2DRiemannConfig04FS2}\includegraphics[width=0.32\textwidth]{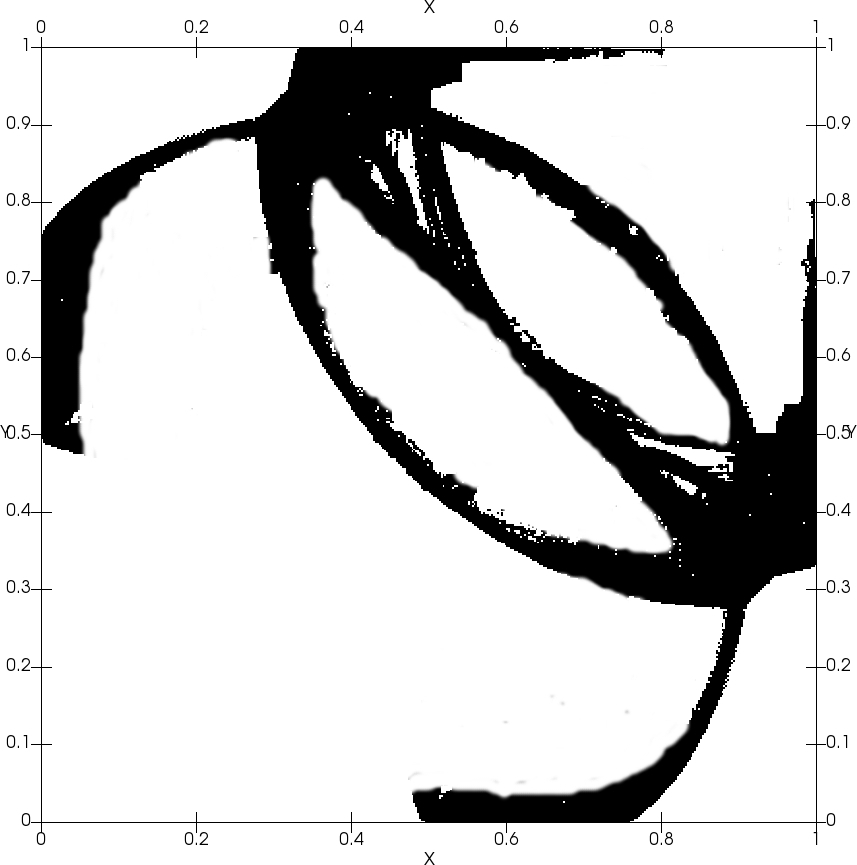}}
  \subfloat[LPR Indicator]{\label{fig:2DRiemannConfig04LPR}\includegraphics[width=0.32\textwidth]{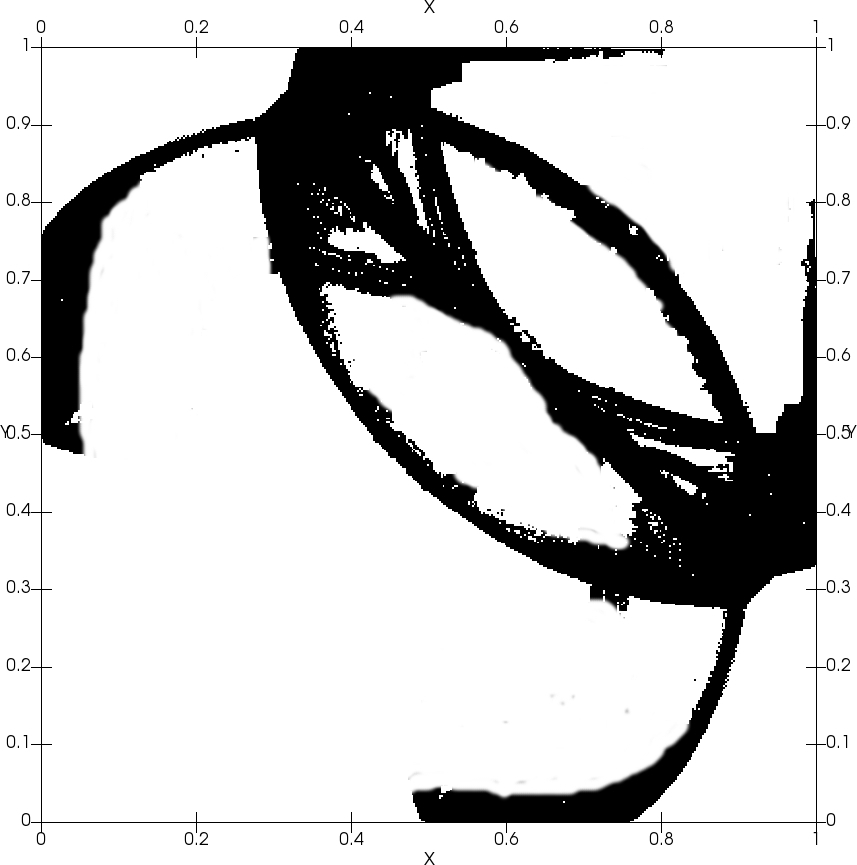}}
  \subfloat[RH Indicator]{\label{fig:2DRiemannConfig04RH}\includegraphics[width=0.32\textwidth]{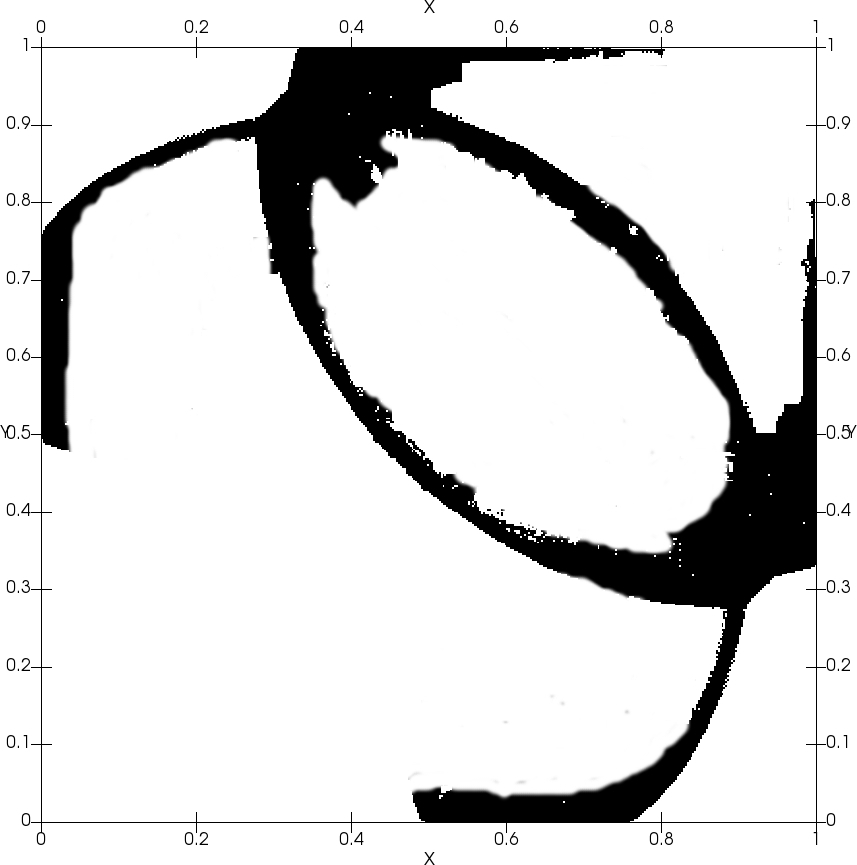}}\hfill
  \subfloat[PPL Indicator]{\label{fig:2DRiemannConfig04PPL}\includegraphics[width=0.32\textwidth]{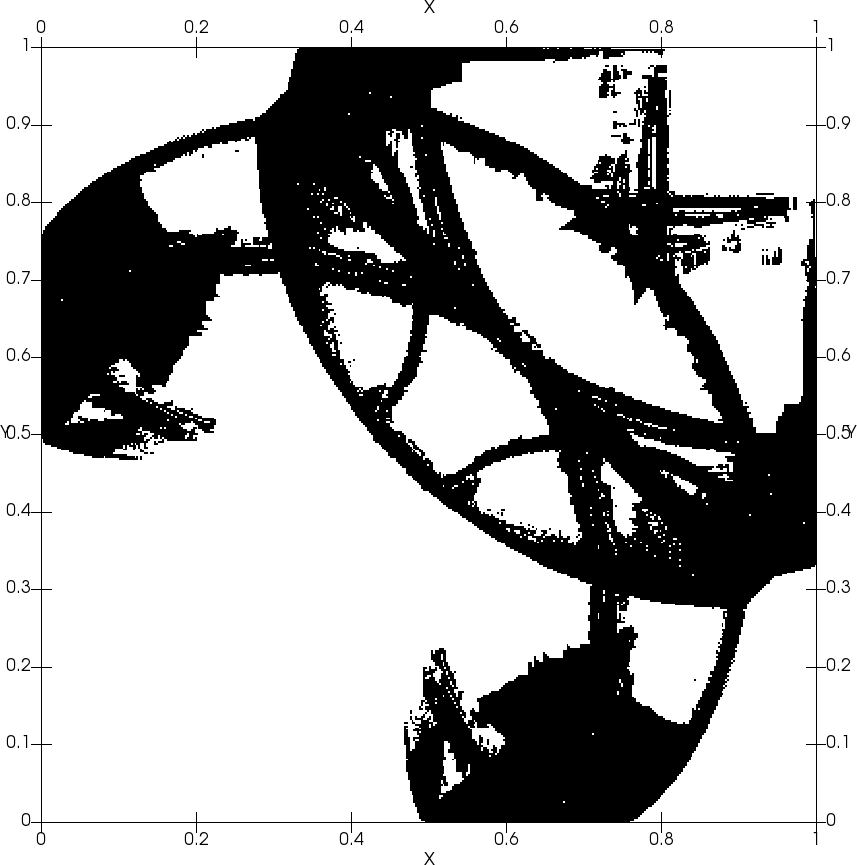}}
  \subfloat[MH Indicator]{\label{fig:2DRiemannConfig04MH}\includegraphics[width=0.32\textwidth]{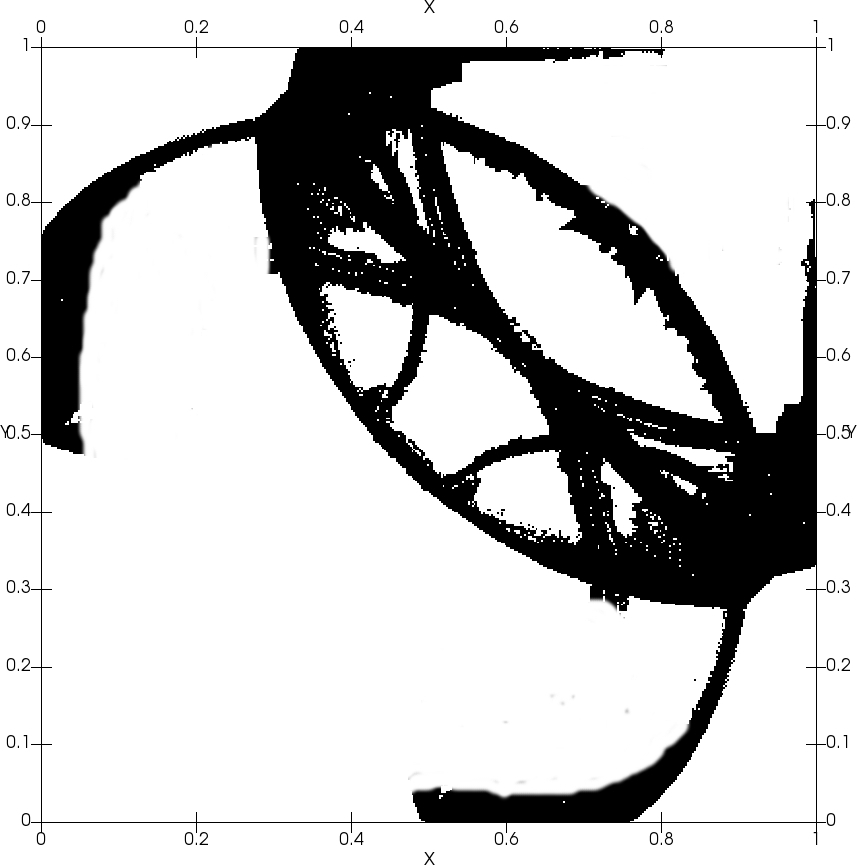}}\hfill
  \caption{The time history of flagged troubled cells of the configuration 4 of the 2D Riemann problem for the two-dimensional Euler equations, simulated until $t=0.25$ with $200\times 200$ elements and $P^{1}$ based DGM}
  \label{fig:2DRiemannConfig04}
\end{figure}

\begin{figure}[htbp]
  \centering
  \subfloat[Density contours for the solution at t=0.25 with $P^{1}$ based DGM]{\label{fig:P1DGM2DRiemannConfig12CSWENOUnstructured}\includegraphics[width=0.45\textwidth]{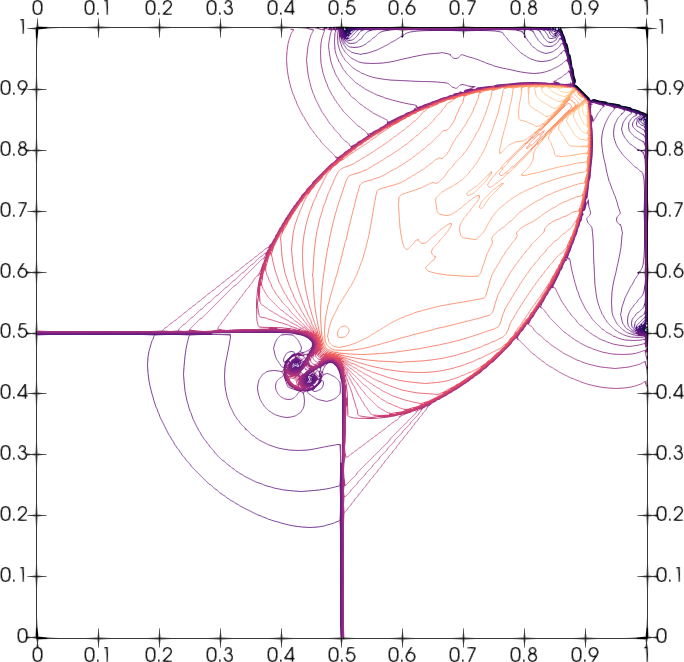}}
  \subfloat[Density contours for the solution at t=0.25 with $P^{2}$ based DGM]{\label{fig:P2DGM2DRiemannConfig12CSWENOUnstructured}\includegraphics[width=0.45\textwidth]{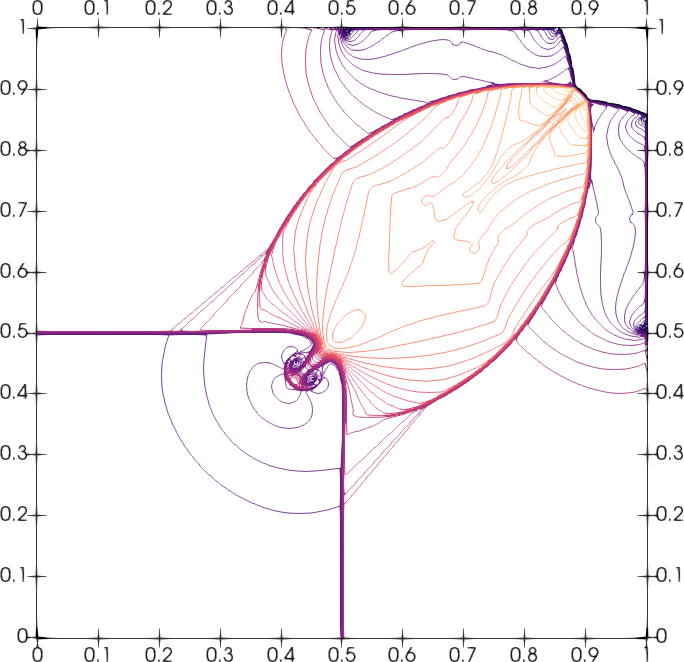}}\hfill
  \subfloat[Density contours for the solution at t=0.25 with $P^{3}$ based DGM]{\label{fig:P3DGM2DRiemannConfig12CSWENOUnstructured}\includegraphics[width=0.45\textwidth]{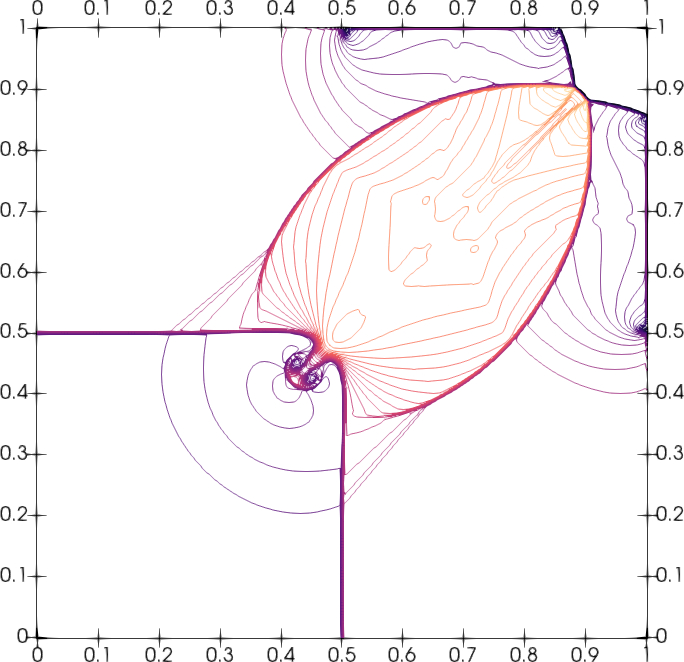}}
  \subfloat[Density contours for the solution at t=0.25 with $P^{4}$ based DGM]{\label{fig:P4DGM2DRiemannConfig12CSWENOUnstructured}\includegraphics[width=0.45\textwidth]{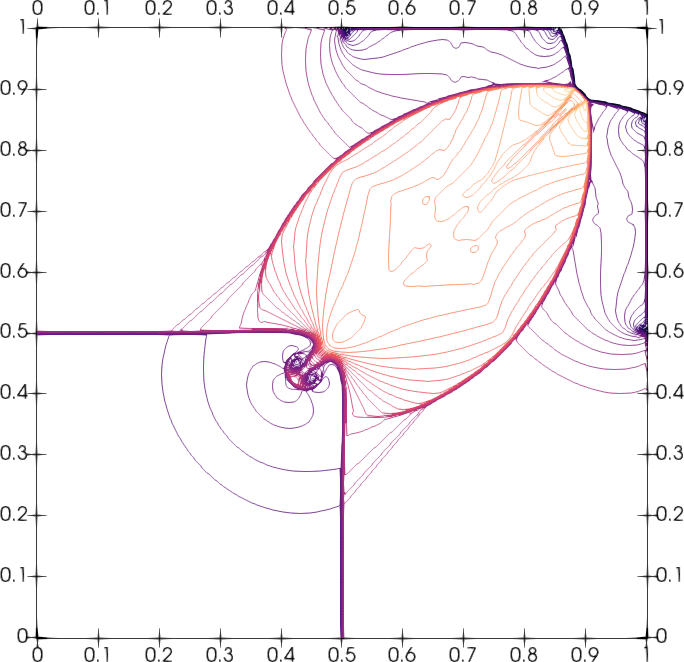}}\hfill
  \subfloat[Density Range]{\label{fig:DensityRange2DRiemannConfig12CSWENOUnstructured}\includegraphics[width=0.9\textwidth]{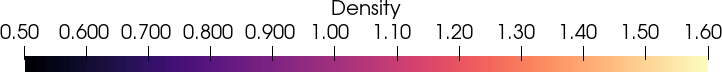}}
  \caption{50 equally spaced density contours for solution at $t=0.25$ for 2D Riemann problem configuration 12 using the CSWENO limiter for $P^{1}$, $P^{2}$, $P^{3}$, and $P^{4}$ based DGM for a mesh containing 92552 triangles}
  \label{fig:CSWENO2DRiemannConfig12SolutionUnstructured}
\end{figure}

\begin{figure}[htbp]
  \centering
  \subfloat[PP Indicator]{\label{fig:2DRiemannConfig12PP}\includegraphics[width=0.32\textwidth]{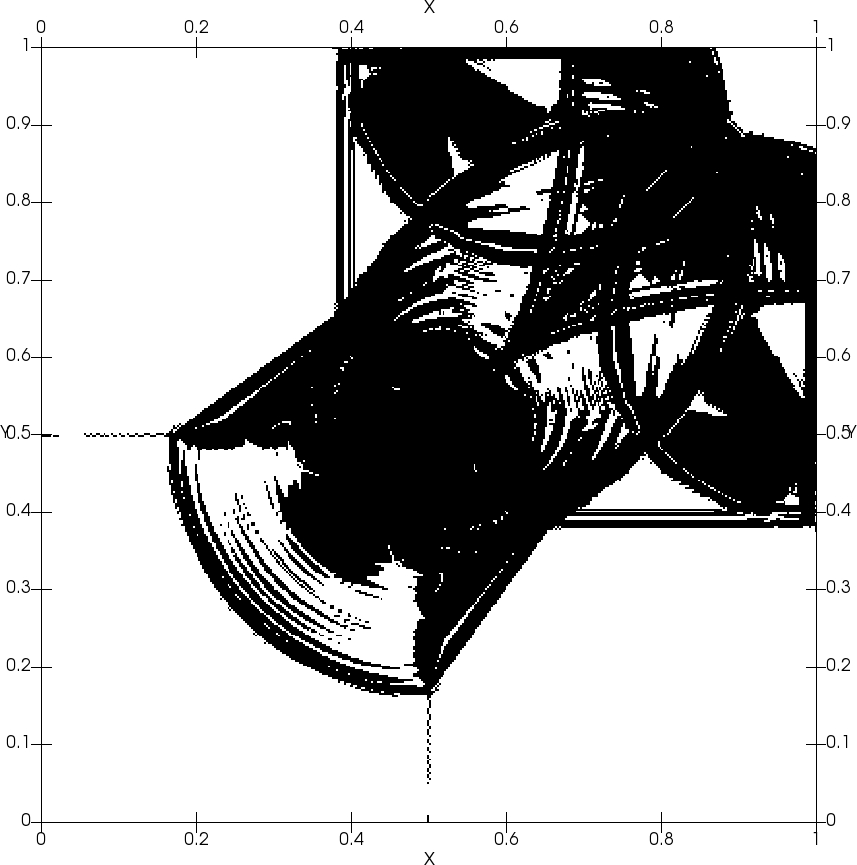}}
  \subfloat[SJ Indicator]{\label{fig:2DRiemannConfig12SJ}\includegraphics[width=0.32\textwidth]{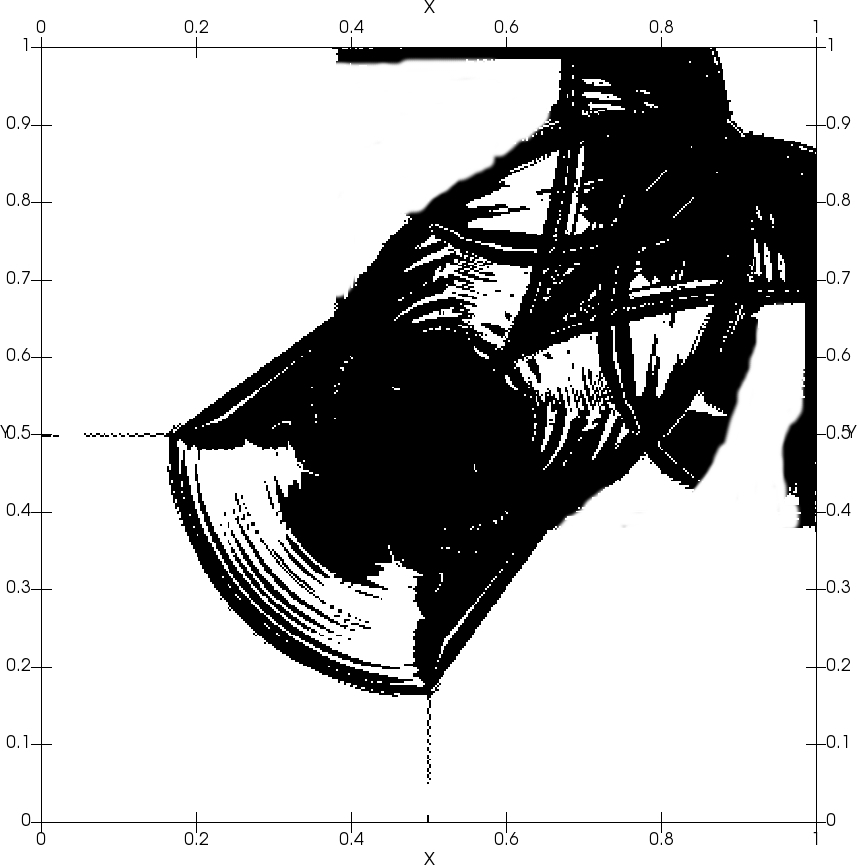}}
  \subfloat[FS1 Indicator]{\label{fig:2DRiemannConfig12FS1}\includegraphics[width=0.32\textwidth]{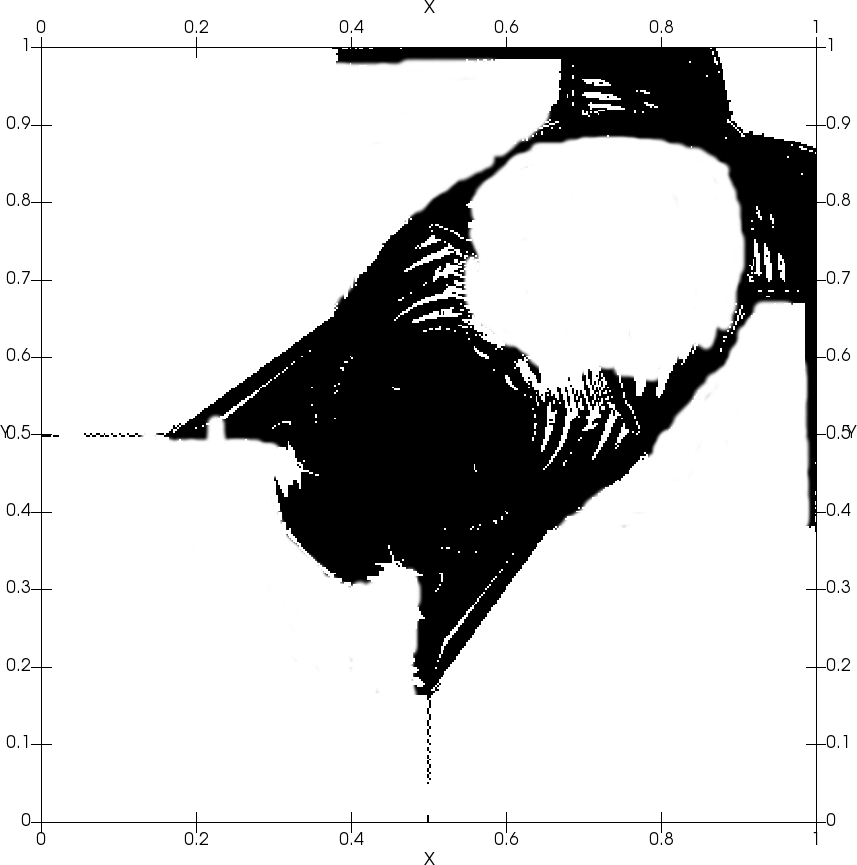}}\hfill
  \subfloat[FS2 Indicator]{\label{fig:2DRiemannConfig12FS2}\includegraphics[width=0.32\textwidth]{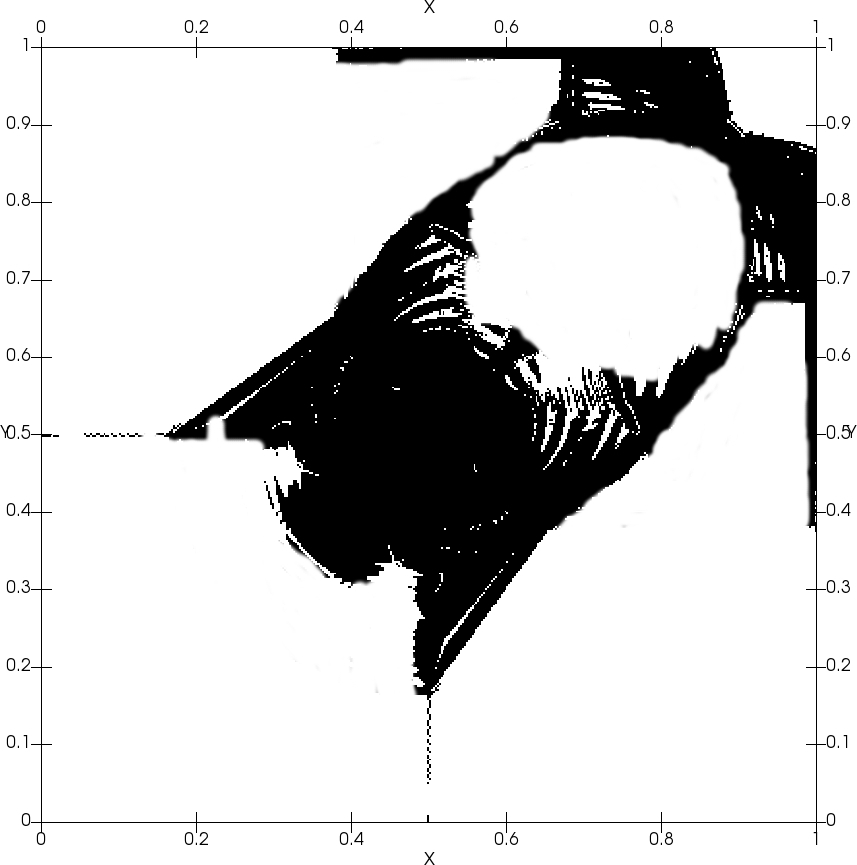}}
  \subfloat[LPR Indicator]{\label{fig:2DRiemannConfig12LPR}\includegraphics[width=0.32\textwidth]{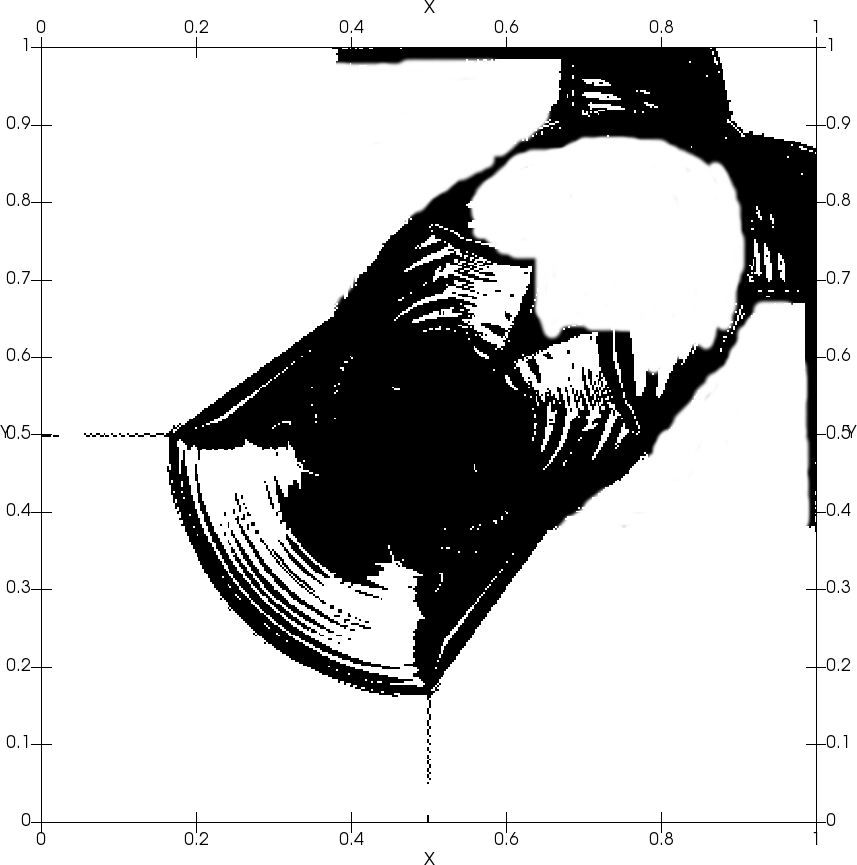}}
  \subfloat[RH Indicator]{\label{fig:2DRiemannConfig12RH}\includegraphics[width=0.32\textwidth]{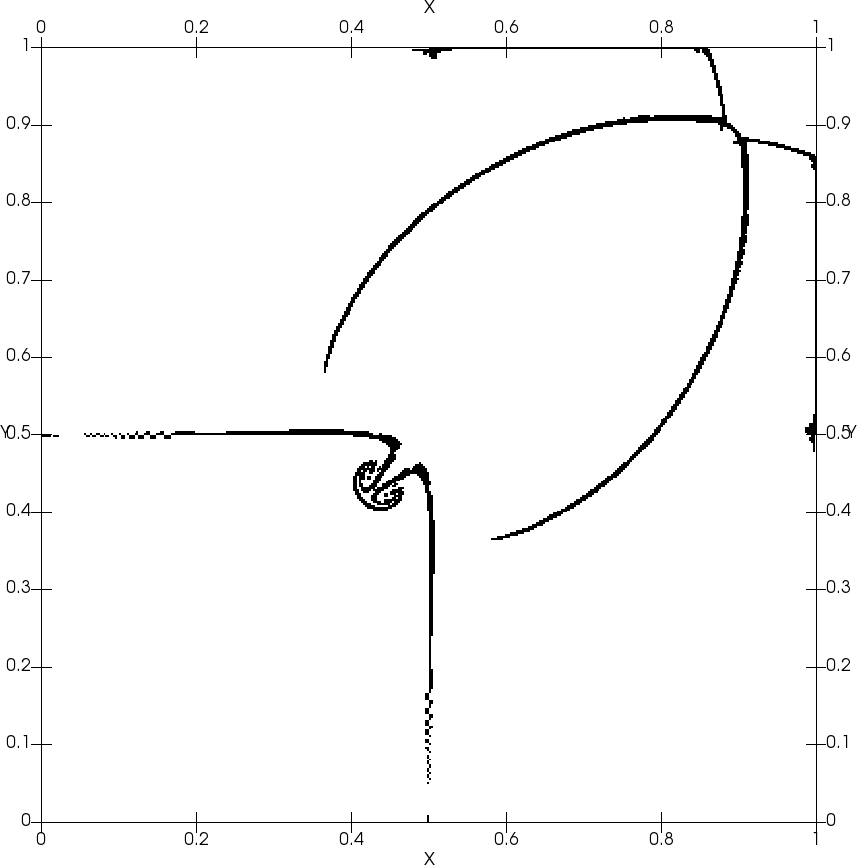}}\hfill
  \subfloat[PPL Indicator]{\label{fig:2DRiemannConfig12PPL}\includegraphics[width=0.32\textwidth]{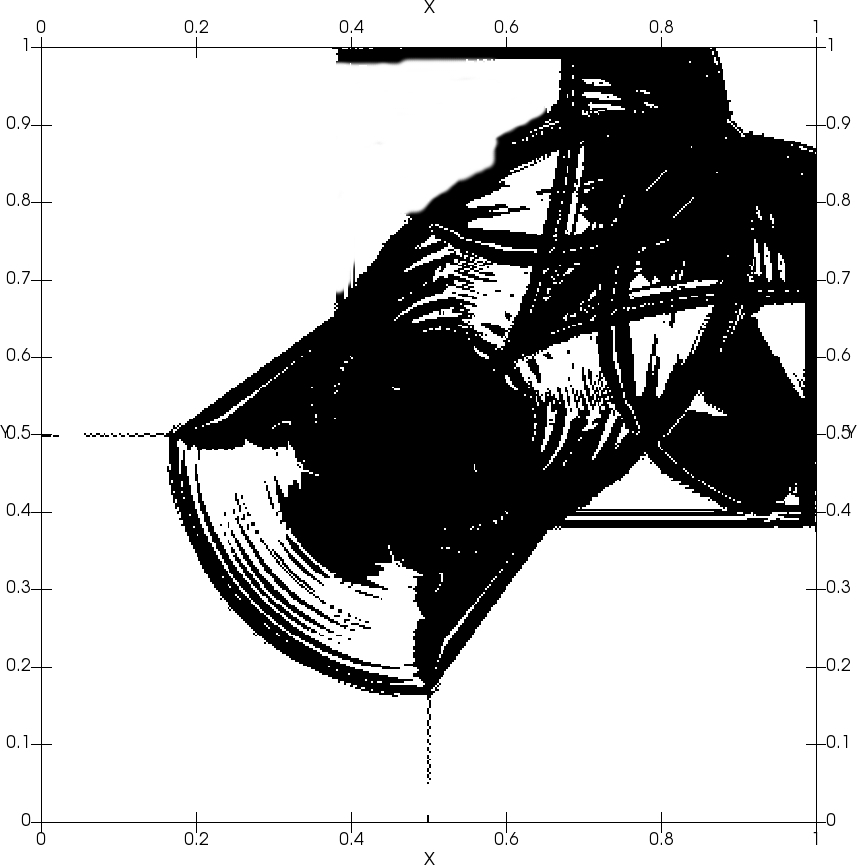}}
  \subfloat[MH Indicator]{\label{fig:2DRiemannConfig12MH}\includegraphics[width=0.32\textwidth]{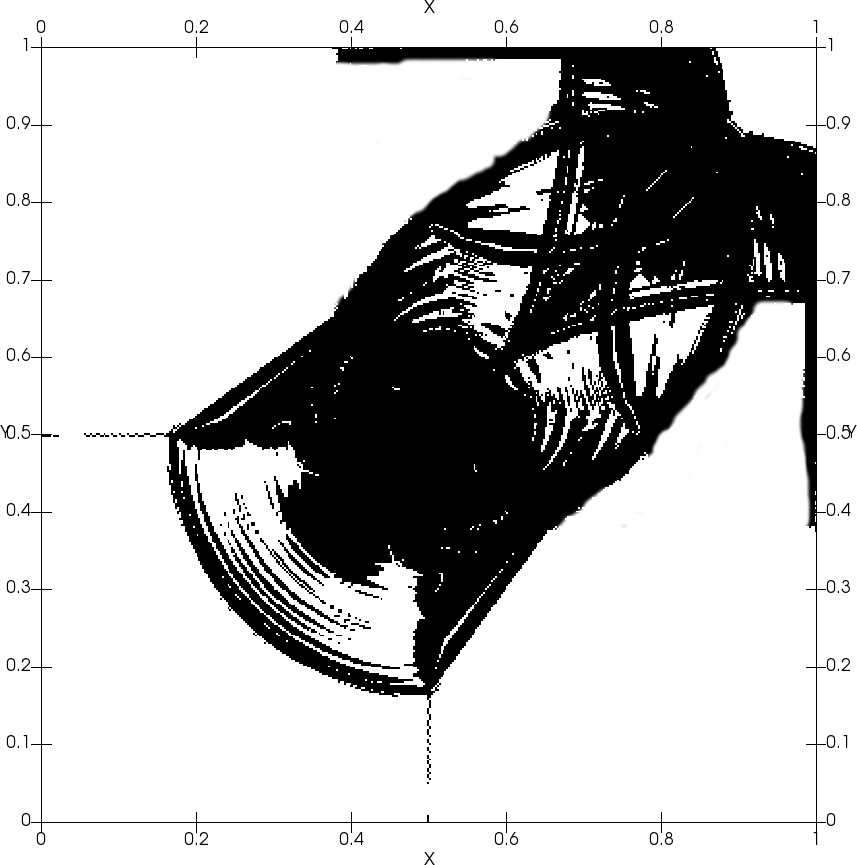}}\hfill
  \caption{The time history of flagged troubled cells of the configuration 12 of the 2D Riemann problem for the two-dimensional Euler equations, simulated until $t=0.25$ with $200\times 200$ elements and $P^{1}$ based DGM}
  \label{fig:2DRiemannConfig12}
\end{figure}

We compute the solution upto time $t=0.25$ for both configurations for structured and unstructured grids. For structured grids, we have used two different uniform meshes with $200\times 200$ and $400\times 400$ elements in each mesh. For unstructured grids (generated using Gmsh 4.6.0 software \cite{GMshGR}), we have used meshes with 92552 and 185104 triangles. Average (over all time steps) and maximum percentages of cells being flagged as troubled cells, for the different troubled-cell indicators, are summarized for structured and unstructured grids respectively in Tables \ref{table:8} and \ref{table:9} for configuration 4 and in Tables \ref{table:10} and \ref{table:11} for configuration 12. The density contours for the solution obtained using the SJ indicator and CSWENO limiter for $P^{1}$, $P^{2}$, $P^{3}$ and $P^{4}$ based DGM using 92552 triangles for configuration 4 and 12 are shown in Figures \ref{fig:CSWENO2DRiemannConfig04SolutionUnstructured} and \ref{fig:CSWENO2DRiemannConfig12SolutionUnstructured} respectively. We are not showing the solutions obtained using other indicators as they are quite similar to the solution obtained with the SJ indicator and stay within an $L^{2}$ norm of $\sim  10^{-14}$. We show the $L^{2}$ difference in density between those solutions for configurations 4 and 12 in tables \ref{table:9New} and \ref{table:11New} respectively.
\\
We also show the troubled cell profile using the eight different indicators in Figure \ref{fig:2DRiemannConfig04} using $P^{1}$ based DGM for $200\times 200$ elements for configuration 4 and in Figure \ref{fig:2DRiemannConfig12} for configuration 12 at $t=0.25$. For all the indicators, the number of flagged troubled cells slightly increase when the order increases. From the tabulated results and the figures, for both configurations, the troubled cell indicators PPL, SJ, MH, and LPR perform in a similar fashion, while the FS1 and FS2 indicators are quite better and the PP indicator is quite bad comparatively. We can also say that, for this problem, the troubled cell indicator RH out performs all the other indicators. From the tables, we also observe that the performance of the troubled-cell indicators for structured and unstructured grids is quite similar.
\\
\\
\begin{table}[htbp]
\small
\centering
\begin{tabular}{|c|c|c|c|c|c|c|c|c|c|}
%\hline
%\multicolumn{6}{|c|}{$L_{2}$ error for the Riemann problem configuration-10} \\
\hline
 \makecell{No. of \\ Cells} & \makecell{Scheme \\ Indicator} & \multicolumn{2}{|c|}{$P^{1}$} & \multicolumn{2}{|c|}{$P^{2}$} & \multicolumn{2}{|c|}{$P^{3}$} & \multicolumn{2}{|c|}{$P^{4}$} \\
\cline{3-10}
 & & Ave & Max & Ave & Max & Ave & Max & Ave & Max \\
\hline 
\multirow{8}{*}{$200\times 200$} & PP & 42.16 & 48.93 & 43.23 & 49.63 & 44.14 & 50.32 & 44.75 & 50.36 \\
\cline{2-10}
 & SJ & 24.23 & 27.46 & 24.93 & 27.99 & 25.26 & 28.39 & 25.74 & 28.88 \\
\cline{2-10}
 & FS1 & 14.15 & 18.28 & 14.65 & 18.76 & 14.98 & 18.98 & 15.14 & 19.15 \\
\cline{2-10}
 & FS2 & 14.26 & 18.29 & 14.66 & 18.78 & 14.99 & 18.94 & 15.16 & 19.16 \\
\cline{2-10}
 & LPR & 17.42 & 21.72 & 17.83 & 21.96 & 18.21 & 22.15 & 18.57 & 22.47 \\
\cline{2-10}
 & RH & 9.17 & 12.76 & 9.42 & 12.98 & 9.88 & 13.16 & 10.24 & 13.45 \\
 \cline{2-10}
 & PPL & 31.65 & 36.32 & 31.92 & 36.54 & 32.26 & 36.94 & 32.56 & 37.16 \\
\cline{2-10}
 & MH & 20.13 & 25.44 & 20.54 & 25.88 & 20.86 & 26.04 & 21.07 & 26.38 \\
\hline
\multirow{8}{*}{$400\times 400$} & PP & 41.24 & 47.86 & 41.57 & 48.13 & 41.98 & 48.59 & 42.32 & 49.02 \\
\cline{2-10}
 & SJ & 22.15 & 26.42 & 22.58 & 26.93 & 22.85 & 27.27 & 23.07 & 27.63 \\
\cline{2-10}
 & FS1 & 13.66 & 17.13 & 13.93 & 17.53 & 14.26 & 17.85 & 14.59 & 18.12 \\
\cline{2-10}
 & FS2 & 13.64 & 17.13 & 13.92 & 17.54 & 14.28 & 17.88 & 14.60 & 18.14 \\
\cline{2-10}
 & LPR & 16.32 & 20.11 & 16.75 & 20.48 & 16.99 & 20.76 & 17.29 & 21.07 \\
\cline{2-10}
 & RH & 8.54 & 11.85 & 8.93 & 12.04 & 9.23 & 12.54 & 9.58 & 12.76 \\
 \cline{2-10}
 & PPL & 30.42 & 35.11 & 30.85 & 35.63 & 31.18 & 35.95 & 31.36 & 32.37 \\
\cline{2-10}
 & MH & 19.78 & 24.76 & 19.96 & 25.02 & 20.18 & 25.45 & 20.58 & 25.76 \\
\hline
\end{tabular}
\caption{Average (marked as Ave) and maximum (marked as Max) percentages of cells flagged as troubled cells subject to different troubled-cell indicators for the 2D Riemann problem - configuration 4, for various orders using structured grids.}
\label{table:8}
\end{table}

\begin{table}[htbp]
\small
\centering
\begin{tabular}{|c|c|c|c|c|c|c|c|c|c|}
%\hline
%\multicolumn{6}{|c|}{$L_{2}$ error for the Riemann problem configuration-10} \\
\hline
 \makecell{No. of \\ Cells} & \makecell{Scheme \\ Indicator} & \multicolumn{2}{|c|}{$P^{1}$} & \multicolumn{2}{|c|}{$P^{2}$} & \multicolumn{2}{|c|}{$P^{3}$} & \multicolumn{2}{|c|}{$P^{4}$} \\
\cline{3-10}
 & & Ave & Max & Ave & Max & Ave & Max & Ave & Max \\
\hline 
\multirow{8}{*}{92552} & PP & 42.13 & 48.88 & 43.26 & 49.69 & 44.18 & 50.39 & 44.72 & 50.33 \\
\cline{2-10}
 & SJ & 24.25 & 27.48 & 24.99 & 27.95 & 25.23 & 28.35 & 25.79 & 28.82 \\
\cline{2-10}
 & FS1 & 14.12 & 18.24 & 14.68 & 18.71 & 14.93 & 18.97 & 15.14 & 19.18 \\
\cline{2-10}
 & FS2 & 14.13 & 18.24 & 14.68 & 18.72 & 14.93 & 18.97 & 15.12 & 19.17 \\
\cline{2-10}
 & LPR & 17.45 & 21.78 & 17.81 & 21.93 & 18.28 & 22.13 & 18.51 & 22.49 \\
\cline{2-10}
 & RH & 9.14 & 12.73 & 9.48 & 12.95 & 9.83 & 13.19 & 10.22 & 13.41 \\
 \cline{2-10}
 & PPL & 31.60 & 36.37 & 31.96 & 36.59 & 32.24 & 36.92 & 32.58 & 37.14 \\
\cline{2-10}
 & MH & 20.17 & 25.47 & 20.53 & 25.80 & 20.84 & 26.09 & 21.04 & 26.32 \\
 \hline
\multirow{8}{*}{185104} & PP & 41.22 & 47.88 & 41.54 & 48.11 & 41.95 & 48.55 & 42.39 & 49.05 \\
\cline{2-10}
 & SJ & 22.18 & 26.45 & 22.53 & 26.92 & 22.89 & 27.23 & 23.07 & 27.64 \\
\cline{2-10}
 & FS1 & 13.62 & 17.14 & 13.98 & 17.56 & 14.26 & 17.86 & 14.55 & 18.15 \\
\cline{2-10}
 & FS2 & 13.64 & 17.14 & 13.98 & 17.55 & 14.26 & 17.86 & 14.56 & 18.14 \\
\cline{2-10}
 & LPR & 16.32 & 20.15 & 16.77 & 20.49 & 16.95 & 20.74 & 17.24 & 21.01 \\
\cline{2-10}
 & RH & 8.59 & 11.83 & 8.99 & 12.02 & 9.28 & 12.54 & 9.58 & 12.74 \\
 \cline{2-10}
 & PPL & 30.42 & 35.18 & 30.88 & 35.65 & 31.13 & 35.93 & 31.34 & 32.36 \\
\cline{2-10}
 & MH & 19.72 & 24.75 & 19.92 & 25.07 & 20.13 & 25.44 & 20.57 & 25.72 \\
\hline
\end{tabular}
\caption{Average (marked as Ave) and maximum (marked as Max) percentages of cells flagged as troubled cells subject to different troubled-cell indicators for the 2D Riemann problem - configuration 4, for various orders using unstructured triangles.}
\label{table:9}
\end{table}

\begin{table}[htbp]
%\small
\centering
\begin{tabular}{|c|c|c|c|c|}
%\hline
%\multicolumn{6}{|c|}{$L_{2}$ error for the Riemann problem configuration-10} \\
\hline
Scheme Indicator & \multicolumn{4}{|c|}{\makecell{$L^{2}$ difference in density in comparison with \\ SJ indicator solution for various orders}} \\ \cline{2-5}
  & $P^{1}$ & $P^{2}$ & $P^{3}$ & $P^{4}$ \\ 
\hline
PP & 7.5E-14 & 1.9E-14 & 6.6E-14 & 5.6E-14\\
\hline
FS1 & 6.2E-14 & 3.0E-14 & 2.1E-14 & 4.9E-14\\
\hline
FS2 & 6.2E-14 & 3.0E-14 & 2.1E-14 & 4.9E-14\\
\hline
LPR & 7.4E-14 & 4.2E-14 & 6.2E-14 & 5.4E-14\\
\hline
RH & 4.3E-14 & 5.2E-14 & 7.3E-14 & 3.2E-14\\
\hline
PPL & 5.5E-14 & 2.9E-14 & 9.1E-14 & 4.0E-14\\
\hline
MH & 6.5E-14 & 6.3E-14 & 9.4E-14 & 5.1E-14\\
\hline
\end{tabular}
\caption{$L^{2}$ difference in density between solution obtained using CSWENO limiter and SJ indicator (shown in Figure \ref{fig:CSWENO2DRiemannConfig04SolutionUnstructured}) and the solution obtained using other indicators for the 2D Riemann problem - configuration 4 using 92552 triangles}
\label{table:9New}
\end{table}

\begin{table}[htbp]
\small
\centering
\begin{tabular}{|c|c|c|c|c|c|c|c|c|c|}
%\hline
%\multicolumn{6}{|c|}{$L_{2}$ error for the Riemann problem configuration-10} \\
\hline
 \makecell{No. of \\ Cells} & \makecell{Scheme \\ Indicator} & \multicolumn{2}{|c|}{$P^{1}$} & \multicolumn{2}{|c|}{$P^{2}$} & \multicolumn{2}{|c|}{$P^{3}$} & \multicolumn{2}{|c|}{$P^{4}$} \\
\cline{3-10}
 & & Ave & Max & Ave & Max & Ave & Max & Ave & Max \\
\hline 
\multirow{8}{*}{$200\times 200$} & PP & 38.16 & 42.27 & 38.65 & 42.35 & 38.95 & 42.69 & 39.23 & 43.25 \\
\cline{2-10}
 & SJ & 28.43 & 32.11 & 28.84 & 32.65 & 29.13 & 32.87 & 29.42 & 33.21 \\
\cline{2-10}
 & FS1 & 16.25 & 18.34 & 16.66 & 18.89 & 17.03 & 19.55 & 17.35 & 19.87 \\
\cline{2-10}
 & FS2 & 16.56 & 18.44 & 16.74 & 18.92 & 17.06 & 19.49 & 17.46 & 19.92 \\
\cline{2-10}
 & LPR & 20.32 & 24.18 & 20.64 & 24.64 & 20.91 & 25.10 & 21.32 & 25.54 \\
\cline{2-10}
 & RH & 7.16 & 11.37 & 7.43 & 11.85 & 7.75 & 12.14 & 7.90 & 12.54 \\
 \cline{2-10}
 & PPL & 32.36 & 36.76 & 32.89 & 37.25 & 33.23 & 37.59 & 33.68 & 37.67 \\
\cline{2-10}
 & MH & 26.54 & 29.16 & 27.11 & 29.59 & 27.36 & 29.81 & 27.62 & 30.03 \\
\hline
\multirow{8}{*}{$400\times 400$} & PP & 36.83 & 41.11 & 37.18 & 41.68 & 37.73 & 41.91 & 37.96 & 42.32 \\
\cline{2-10}
 & SJ & 27.17 & 31.03 & 27.54 & 31.55 & 27.94 & 31.83 & 28.32 & 32.89 \\
\cline{2-10}
 & FS1 & 15.15 & 17.32 & 15.87 & 17.66 & 16.11 & 17.89 & 16.43 & 18.18 \\
\cline{2-10}
 & FS2 & 15.15 & 17.32 & 15.87 & 17.66 & 16.11 & 17.89 & 16.43 & 18.18 \\
\cline{2-10}
 & LPR & 19.18 & 23.00 & 19.54 & 23.47 & 19.84 & 23.78 & 20.15 & 24.16 \\
\cline{2-10}
 & RH & 6.54 & 10.11 & 6.88 & 10.38 & 7.23 & 10.69 & 7.64 & 10.95 \\
 \cline{2-10}
 & PPL & 31.13 & 35.65 & 31.65 & 35.89 & 31.86 & 36.24 & 32.06 & 36.57 \\
\cline{2-10}
 & MH & 25.15 & 28.01 & 25.57 & 28.43 & 25.92 & 28.64 & 26.25 & 28.91 \\
\hline
\end{tabular}
\caption{Average (marked as Ave) and maximum (marked as Max) percentages of cells flagged as troubled cells subject to different troubled-cell indicators for the 2D Riemann problem - configuration 12, for various orders using structured grids.}
\label{table:10}
\end{table}

\begin{table}[htbp]
\small
\centering
\begin{tabular}{|c|c|c|c|c|c|c|c|c|c|}
%\hline
%\multicolumn{6}{|c|}{$L_{2}$ error for the Riemann problem configuration-10} \\
\hline
 \makecell{No. of \\ Cells} & \makecell{Scheme \\ Indicator} & \multicolumn{2}{|c|}{$P^{1}$} & \multicolumn{2}{|c|}{$P^{2}$} & \multicolumn{2}{|c|}{$P^{3}$} & \multicolumn{2}{|c|}{$P^{4}$} \\
\cline{3-10}
 & & Ave & Max & Ave & Max & Ave & Max & Ave & Max \\
\hline 
\multirow{8}{*}{92552} & PP & 38.12 & 42.21 & 38.63 & 42.32 & 38.97 & 42.63 & 39.29 & 43.20 \\
\cline{2-10}
 & SJ & 28.41 & 32.17 & 28.82 & 32.60 & 29.14 & 32.81 & 29.43 & 33.27 \\
\cline{2-10}
 & FS1 & 16.28 & 18.31 & 16.64 & 18.84 & 17.00 & 19.53 & 17.31 & 19.84 \\
\cline{2-10}
 & FS2 & 16.28 & 18.31 & 16.64 & 18.84 & 17.00 & 19.53 & 17.31 & 19.84 \\
\cline{2-10}
 & LPR & 20.34 & 24.19 & 20.63 & 24.61 & 20.96 & 25.18 & 21.34 & 25.51 \\
\cline{2-10}
 & RH & 7.12 & 11.39 & 7.41 & 11.89 & 7.72 & 12.13 & 7.98 & 12.54 \\
 \cline{2-10}
 & PPL & 32.39 & 36.72 & 32.82 & 37.27 & 33.23 & 37.58 & 33.63 & 37.67 \\
\cline{2-10}
 & MH & 26.55 & 29.13 & 27.15 & 29.51 & 27.39 & 29.87 & 27.64 & 30.01 \\
\hline
\multirow{8}{*}{185104} & PP & 36.88 & 41.13 & 37.15 & 41.63 & 37.79 & 41.99 & 37.91 & 42.34 \\
\cline{2-10}
 & SJ & 27.12 & 31.08 & 27.52 & 31.58 & 27.91 & 31.82 & 28.35 & 32.83 \\
\cline{2-10}
 & FS1 & 15.14 & 17.33 & 15.84 & 17.63 & 16.12 & 17.84 & 16.46 & 18.12 \\
\cline{2-10}
 & FS2 & 15.14 & 17.33 & 15.84 & 17.63 & 16.12 & 17.84 & 16.46 & 18.12 \\
\cline{2-10}
 & LPR & 19.12 & 23.05 & 19.59 & 23.42 & 19.88 & 23.73 & 20.12 & 24.11 \\
\cline{2-10}
 & RH & 6.51 & 10.18 & 6.81 & 10.32 & 7.28 & 10.61 & 7.66 & 10.93 \\
 \cline{2-10}
 & PPL & 31.14 & 35.64 & 31.66 & 35.81 & 31.89 & 36.21 & 32.09 & 36.53 \\
\cline{2-10}
 & MH & 25.13 & 28.08 & 25.52 & 28.47 & 25.90 & 28.61 & 26.21 & 28.87 \\
\hline
\end{tabular}
\caption{Average (marked as Ave) and maximum (marked as Max) percentages of cells flagged as troubled cells subject to different troubled-cell indicators for the 2D Riemann problem - configuration 12, for various orders using unstructured triangles.}
\label{table:11}
\end{table}

\begin{table}[htbp]
%\small
\centering
\begin{tabular}{|c|c|c|c|c|}
%\hline
%\multicolumn{6}{|c|}{$L_{2}$ error for the Riemann problem configuration-10} \\
\hline
Scheme Indicator & \multicolumn{4}{|c|}{\makecell{$L^{2}$ difference in density in comparison with \\ SJ indicator solution for various orders}} \\ \cline{2-5}
  & $P^{1}$ & $P^{2}$ & $P^{3}$ & $P^{4}$ \\ 
\hline
PP & 3.3E-14 & 8.3E-14 & 5.2E-14 & 5.4E-14\\
\hline
FS1 & 2.1E-14 & 4.1E-14 & 6.5E-14 & 2.9E-14\\
\hline
FS2 & 2.1E-14 & 4.1E-14 & 6.5E-14 & 2.9E-14\\
\hline
LPR & 6.4E-14 & 4.9E-14 & 8.2E-14 & 7.5E-14\\
\hline
RH & 5.3E-14 & 9.1E-14 & 8.5E-14 & 2.8E-14\\
\hline
PPL & 7.7E-14 & 8.2E-14 & 3.7E-14 & 4.1E-14\\
\hline
MH & 8.1E-14 & 7.3E-14 & 4.3E-14 & 7.0E-14\\
\hline
\end{tabular}
\caption{$L^{2}$ difference in density between solution obtained using CSWENO limiter and SJ indicator (shown in Figure \ref{fig:CSWENO2DRiemannConfig12SolutionUnstructured}) and the solution obtained using other indicators for the 2D Riemann problem - configuration 12 using 92552 triangles}
\label{table:11New}
\end{table}

\section{Conclusion}\label{sec:4}

\noindent In this paper, we have studied eight troubled cell indicators for the discontinuous Galerkin method using a WENO limiter on both structured and unstructured grids. We used various one-dimensional and two-dimensional problems for various orders on the hyperbolic system of Euler equations to compare these troubled cell indicators. We evaluated them based on the percentage of cells flagged as troubled cells for various orders and various grid sizes. We have observed that for each of the selected troubled cell indicators, the CPU time taken to test a single cell for discontinuity is approximately the same (excluding the time taken to train a neural network for the RH indicator). We have also observed that as long as we use an indicator which identifies all the troubled cells, if we use a good limiter, the accuracy of the solution remains the same. But the additional troubled cells add to the computational time and reduce the efficiency of the scheme. For one-dimensional problems, the performance of Fu and Shu indicator \cite{fs1} and the modified KXRCF indicator \cite{lpr} is better than other indicators. For two-dimensional problems, the performance of the artificial neural network (ANN) indicator of Ray and Hesthaven \cite{rh1} is quite good and the Fu and Shu and the modified KXRCF indicators are also good. We can conclude that these three indicators are suitable candidates for applications of DGM using WENO limiters though it should be noted that the ANN indicator is quite expensive and requires a lot of training.
\\
\\




\noindent \textbf{Acknowledgements:} We thank Dr. P.A. Ramakrishna for his valuable support.

\end{document}